\documentclass[twoside,reqno,a4paper]{amsart}

\usepackage{graphicx}

\usepackage{xcolor}
\definecolor{webgreen}{rgb}{0,.5,0}
\definecolor{webbrown}{rgb}{.6,0,0}
\definecolor{RoyalBlue}{cmyk}{1, 0.50, 0, 0}
\usepackage{epsfig, graphicx, subfigure,tikz}
\usepackage{verbatim, setspace}
\usepackage{amsmath, amssymb}

\usepackage{style}

\usepackage{pdfpages}
\usepackage{chngcntr} 
\counterwithin{section}{part}

\newcommand{\ddd }  {\stackrel{\rm def}{=}}
\newcommand{\n}{{\vec n}}

\newcommand{\HN}{\mathrm{HN}}

\newcommand{\RS}{\boldsymbol{\mathfrak{R}}}
\newcommand{\z}	{{\boldsymbol z}}

\newcommand{\cws}{\stackrel{*}{\to}}

\newcommand{\mjo}{\mathfrak{J}}
\newcommand{\jo}{\mathcal{J}}

\newcommand{\jack}{\jo_{\vec\kappa}}

\newcommand{\jox}{\jo_{[X]}}
\newcommand{\jackn}{\jo_{\vec\kappa,\vec N}}

\newcommand{\jnt}{\mathrm{Joint}(E)}
\newcommand{\jnts}{\mathrm{Joint}^*(E)}

\begin{document}
\title[ Spectral theory of Jacobi matrices on trees whose coefficients are generated by multiple orthogonality]{Spectral theory of Jacobi matrices on trees  whose coefficients are generated by multiple orthogonality}

\thanks{ 
The research of the first author was supported by the
grant NSF DMS-1764245, the grant from the Moscow Center for Fundamental and Applied Mathematics, project No. 20-03-01,  and Van Vleck Professorship Research Award. 
The second author's research was supported by the Moscow Center for Fundamental and Applied Mathematics, project No. 20-03-01, and grant CGM-354538 from the Simons Foundation.}

\author[S. Denisov]{Sergey A. Denisov}
\address{Department of Mathematics, University of Wisconsin-Madison, 480 Lincoln Dr., Madison, WI 53706, USA}
\address{Keldysh Institute of Applied Mathematics, Russian Academy of Science, Miusskaya Pl. 4, Moscow, 125047 Russian Federation}
\email{\href{mailto:denissov@math.wisc.edu}{denissov@math.wisc.edu}}

\author[M. Yattselev]{Maxim L. Yattselev}
\address{Department of Mathematical Sciences, Indiana University-Purdue University Indianapolis, 402~North Blackford Street, Indianapolis, IN 46202, USA}
\address{Keldysh Institute of Applied Mathematics, Russian Academy of Science, Miusskaya Pl. 4, Moscow, 125047 Russian Federation}
\email{\href{mailto:maxyatts@iupui.edu}{maxyatts@iupui.edu}}


\begin{abstract} We study  Jacobi matrices on trees whose coefficients are  generated by multiple orthogonal polynomials.  Hilbert space decomposition into an orthogonal sum of cyclic subspaces is obtained. For each subspace, we find generators and the generalized eigenfunctions  written in terms of the orthogonal polynomials. The spectrum and its spectral type are studied for large classes of orthogonality measures.

\end{abstract}

\maketitle

\setcounter{tocdepth}{1}

\tableofcontents


\part*{Introduction}

This paper is the third in the sequence of works  \cite{ady1,ady2} that study the connection between  Jacobi matrices on  trees and the theory of multiple orthogonal polynomials (MOPs). In \cite{ady1}, we have described a large class of MOPs that generate bounded and self-adjoint Jacobi matrices on rooted homogeneous trees and established some basic facts explaining this connection. In particular, we constructed a bijection between MOPs of the first type and a class of such Jacobi matrices. In the follow-up paper \cite{ady2}, we performed a case study of the Angelesco systems generated by two measures of orthogonality with  analytic densities. We used Riemann-Hilbert analysis to obtain asymptotics of MOPs and their recurrence coefficients.  That led to a complete description of all the ``right limits''  of these Jacobi matrices and allowed us to find their essential spectrum. In the current paper, we study the spectrum and spectral decomposition in a more general situation. We focus on the case of two measures only and  address several questions that were left open in \cite{ady1}. 

The rest of the paper is organized as follows. In the remaining part of the introduction, we emphasize the importance of Jacobi matrices, outline their connection to orthogonal polynomials, provide a general definition of Jacobi matrices on graphs, and state some of the properties of multiple orthogonal polynomials on the real line that we need to study the Jacobi matrices we are interested in. After that, we focus exclusively on the study of spectral properties of Jacobi operators on trees generated by MOPs on the real line. In Part~\ref{part1}, we provide a full Spectral Theorem for finite Jacobi matrices. In Part~\ref{part2}, we define Jacobi matrices on a 2-homogeneous infinite rooted Cayley tree and discuss some of their basic properties. In Part~\ref{part3}, we study Jacobi matrices generated by Angelesco systems and describe cyclic subspaces, generalized eigenfunctions, and the corresponding spectral measures. Part~\ref{part4} contains the spectral decomposition for Jacobi matrices on rooted trees with periodic coefficients. That complements the construction in Part~\ref{part3}.

\section*{Orthogonal decomposition and spectrum}

We recall some basic facts from  the spectral theory of bounded self-adjoin operators (see, \cite{akh1,akh2} and \cite[Section VII.2]{rs1}). Let $\mathfrak{H}$ be a Hilbert space and $\mathfrak{A}$ be a bounded self-adjoint operator acting on it. We can study the spectrum of this operator by obtaining a decomposition of $\mathfrak{H}$ into an orthogonal sum of cyclic subspaces of  $\mathfrak{A}$. That is, take any $\mathfrak{g}_1\in \mathfrak{H}$ with unit norm, i.e., $\|\mathfrak{g}_1\|=1$, and generate the cyclic subspace 
\[
\mathfrak{C}_1 \ddd \overline{\text{span}\{\mathfrak{A}^m\mathfrak{g}_1:~m=0,1,\ldots\}}.
\]
We shall call $\mathfrak{g}_1$ the first generator and $\mathfrak{C}_1$ the first cyclic subspace. One can show that $\mathfrak{C}_1$ is invariant with respect to $\mathfrak{A}$. If  $\mathfrak{C}_1\subset \mathfrak{H}$, we take $\mathfrak{g}_2\in \mathfrak{H}$, that satisfies $\|\mathfrak{g}_2\|=1$ and $\mathfrak{g}_2\perp \mathfrak{C}_1$. We denote by $\mathfrak{C}_2$ the cyclic space generated by $\mathfrak{g}_2$. It is also invariant under $\mathfrak{A}$ and satisfies $\mathfrak{C}_1\perp \mathfrak{C}_2$. Continuing that way, we obtain the following representation of \( \mathfrak H \) as a sum of orthogonal cyclic subspaces:
\begin{equation}
\label{f1}
\mathfrak{H}=\oplus_{m=1}^N\mathfrak{C}_m,
\end{equation}
where $N\in \mathbb{N}\cup \infty$ is called the multiplicity of the spectrum. Since $\mathfrak{A}$ is self-adjoint,  the operator $(\mathfrak A-z)^{-1}$ is bounded on $\mathfrak H$ for every $z\in \C_+$, the upper half-plane. For each $f\in \mathfrak H$, the function $\langle(\mathfrak A-z)^{-1}f,f\rangle$ is in Herglotz-Nevanlinna class, i.e., it is analytic in $\C_+$ and has non-negative imaginary part there (we discuss this class below, see \eqref{df1}).  Moreover, since $\mathfrak A$ is bounded, we have an integral representation
\begin{equation}
\label{spectral-m}
\langle (\mathfrak A-z)^{-1}f,f\rangle=\int_{\R}\frac{d\rho_f(x)}{x-z}, \quad z\in \C_+\,,
\end{equation}
where the measure $\rho_f$ is called the \emph{spectral measure} of $f$. Then, the following result holds.

\begin{Thm}
\label{thm:spectrum}
Let \( \mathfrak A \) be a bounded self-adjoint operator on a Hilbert space \( \mathfrak H \) and let \( \sigma(\mathfrak A) \) denote its spectrum. It holds that
\[
\sigma(\mathfrak A) = \bigcup_{m=1}^\infty \supp\,\rho_{\mathfrak g_m}, 
\]
where \( \rho_{\mathfrak g_m} \) is the spectral measure of the generator \( \mathfrak g_m \) for the cyclic subspace \( \mathfrak C_m \) from  decomposition \eqref{f1}.
\end{Thm}

Decomposition \eqref{f1} can be used as follows. Fix $\mathfrak{C}_m$. Taking a sequence of vectors 
\[
\{\mathfrak{g}_m,\mathfrak{A}\mathfrak{g}_m,\mathfrak{A}^2\mathfrak{g}_m, \ldots\}
\]
and running Gramm-Schmidt orthogonalization procedure gives the orthonormal basis in $\mathfrak{C}_m$ in which the restriction of $\mathfrak{A}$ to $\mathfrak{C}_m$ takes the form of either an infinite or a finite (depending on \( \dim \mathfrak{C}_m \)) one-sided  Jacobi matrix, see \eqref{1.20} and \eqref{sd_a31}, further below. It turns out that these matrices are  related to orthogonal polynomials, a connection that is central to our interested in the subject.

\section*{Classical Jacobi matrices}

Let $\{a_j\},\{b_j\}\in \ell^\infty(\Z_+)$ and $a_j>0, b_j\in \R$, hereafter \( \Z_+\ddd \{0,1,2,\ldots\} \) and \( \mathbb N \ddd \{1,2,\ldots\} \). The infinite one-sided  Jacobi matrix is a matrix of the form
 \begin{equation}
 \label{1.20}
\mjo\ddd   \left[
\begin{array}{ccccc}
b_0 & \sqrt{a_0} & 0& 0 &\ldots\\
\sqrt{a_0} & b_1 &\sqrt{a_1}&0&\ldots \\
0 & \sqrt{a_1}& b_2& \sqrt{a_2}&\ldots\\
0 & 0 & \sqrt{a_2} & b_3&\ldots\\
\ldots &\ldots&\ldots&\ldots&\ldots\\
\end{array}
\right]\,,
\end{equation}
and \( N \)--dimensional Jacobi matrix is the upper-left $N\times N$ corner of \eqref{1.20}, see \eqref{sd_a31} further below.  We define two sets of measures on the real line 
\[
\mathfrak M \ddd \big\{\mu:~\supp \mu\subset[-R_\mu,R_\mu],~R_\mu<\infty,~\text{and}~~\#\supp \mu=\infty\big\} \quad \text{and} \quad \mathfrak{M}_1 \ddd \big\{\mu\in\mathfrak M:~\mu(\R)=1\big\},
\]
where the cardinality of a set \( S \) is denoted by \( \#S \). One-sided infinite Jacobi matrices with uniformly bounded entries are known to be in one-to-one correspondence with \( \mathfrak{M}_1 \),  the set  of probability measures on $\R$ whose support is compact and has infinite cardinality. This bijection is realized via polynomials orthogonal on the real line. On one hand, since $\mjo$  defines a bounded self-adjoint operator on the Hilbert space $\ell^2(\Z_+)$, we can consider the spectral measure  of the vector $(1,0,0,\ldots)$, see \eqref{spectral-m}.  We will call it \( \rho(\mjo) \). On the other hand, given $\mu\in \mathfrak{M}_1$, one can produce a Jacobi matrix in the following way. Let $p_n(x,\mu)$ be the \( n \)-th orthonormal polynomial with respect to~$\mu$, i.e., $p_n(x,\mu)$ is a polynomial of degree \( n \) such that
\[
\int_\R p_n(x,\mu)x^md\mu(x)=0,\quad  m=0,\ldots,n-1\,,
\]
that is normalized so that
\[
 {\rm coeff}_n p_n>0, \quad \int_\R p^2_n(x,\mu)d\mu(x)=1\,,
\]
where ${\rm coeff}_n Q$ is the coefficient in  front of  $x^n$ of the polynomial $Q(x)$. It is known that polynomials $p_n(x,\mu)$ satisfy the three-term recurrence relations
\begin{equation}
\label{1.19}
xp_n(x,\mu)=\sqrt{a_n}p_{n+1}(x,\mu)+b_np_n(x,\mu)+\sqrt{a_{n-1}}p_{n-1}(x,\mu), \quad n=0,1,\ldots\,,
\end{equation}
where $a_n>0, \, b_n\in \R$ and $p_{-1}\ddd  0, a_{-1}\ddd  0$. The coefficients $\{a_n\}, \{b_n\}$ are defined uniquely by $\mu$ and one can show that 
\[
\{a_n\},\{b_n\}\in \ell^\infty(\Z_+)\,.
\]
Let \( \mjo \) be defined via \eqref{1.20} with these coefficients. It is a general fact of the theory \cite{akh1,akh2} that
\begin{equation}
\label{sd_21}
\rho(\mjo) = \mu \quad \text{and therefore} \quad \sigma(\mjo) = \supp \mu\,.
\end{equation}

The above correspondence is one-to-one: one can start with a bounded self-adjoint Jacobi matrix \eqref{1.20}, compute \( \rho(\mjo) \), the spectral measure of $(1,0,0,\ldots)$, via \eqref{spectral-m}, take $\rho(\mjo)$ as a measure of orthogonality $\mu$ and, finally, define the orthogonal polynomials whose recurrence coefficients will give rise to the same $\mjo$. 

It follows from \eqref{1.19} that the sequence $\{p_n(x,\mu)\}$, with \( \mu=\rho(\mjo) \), represents the generalized eigenfunction of $\mjo$. That can be made explicit by the following statement, see \cite{akh1,akh2}, which, together with \eqref{1.19}, can be taken as a definition of a generalized eigenfunction.

\begin{Prop}
\label{lem1} 
Suppose $\mu \in \mathfrak{M}_1$. The map
\[
\alpha(x) \mapsto \widehat \alpha=\big\{\widehat\alpha(n)\big\}_{n\in\Z_+}, \quad \widehat\alpha(n) \ddd \int \alpha(x)p_n(x,\mu)d\mu(x),
\]
is a unitary map from $L^2(\mu)$ onto $\ell^2(\Z_+)$ such that
\[
\|\alpha\|_{L^2(\mu)}^2=\|\widehat\alpha\|_{\ell^2(\Z_+)}^2\,.
\]
This map establishes unitary equivalence of the operator \( \mjo \) on $\ell^2(\Z_+)$ and the operator of multiplication by \( x \) on $L^2(\mu)$. In particular,
\[
x\alpha(x) \mapsto \mjo\widehat \alpha.
\]
\end{Prop}

Finite Jacobi matrices can also be studied via polynomials orthogonal on the real line although the measure of orthogonality giving rise to a particular 
\begin{equation}
\label{sd_a31}
\mjo_N\ddd  \left[
\begin{array}{cccccc}
b_0 & \sqrt{a_0} & 0&                  \ldots &       \ldots &0\\
\sqrt{a_0} & b_1 & \sqrt{a_1}&                   \ldots &       \ldots &0 \\
0 & \sqrt{a_1}& b_2&      \ldots &      \ldots &  0\\
\ldots & \ldots & \ldots &        \ldots &     \ldots&\ldots\\
0 &0&0&                           \ldots &     \sqrt{a_{N-1}} & b_{N}\\
\end{array}
\right]\,
\end{equation}
is not unique, which has to do with multiple solutions to a moment problem, see \cite{akh1}. Let \( \mu \) be any measure of orthogonality such that \( \mjo_N \) is upper-left $(N+1)\times (N+1)$ corner of \( \mjo \) generated by the orthogonal polynomials \( \{p_n(x,\mu)\} \). If $\vec{p}_{N}\ddd (p_0,\ldots,p_N)$, we get
\begin{equation}
\label{1.21a}
(\mjo_N-x)\vec{p}_N(x) = -\sqrt{a_{N}}p_{N+1}(x)\delta^{(N)}, \quad \delta^{(N)} \ddd (0,\ldots,0,1)\,.
\end{equation}
The last identity provides, in particular, the characterization of the spectrum of $\mjo_N$:
\begin{equation}\label{sar_1}
\sigma(\mjo_N)=\{E:~p_{N+1}(E,\mu)=0\}\,.
\end{equation}

\section*{Jacobi matrices on graphs}

We are interested in the generalizations of the above notion of a Jacobi matrix to the case when underlying Hilbert space is realized not as \( \ell^2(\Z_+) \), but as a space of square-integrable functions on vertices of a tree. 

Let $\mathcal{G}=(\mathcal V,\mathcal E)$ be an infinite graph, where $\mathcal{V}$ and \( \mathcal E \) stand for the sets of its vertices and edges, respectively. The set of directed edges will be denoted by $\vec{\mathcal E}$. For $Y\in \mathcal{V}$, the symbol $\delta^{(Y)}$ indicates the Kronecker symbol at $Y$, i.e., the function which is equal to $1$ at $Y$ and zero otherwise. Given two vertices $V_1,V_2\in \mathcal{V}$, we shall write $V_1\sim V_2$ if they are connected by an edge and also use this notation to denote the edge itself. The edge directed from $V_1$ to $V_2$ will be denoted by $[V_1,V_2]$.

A connected graph that has no loops is called a \emph{tree}, in which case we shall use the symbol $\mathcal{T}$ instead of \( \mathcal G \). If every vertex in a tree has the same number of neighbors, this tree is called \emph{homogeneous}. We can construct a rooted homogeneous tree of degree $d+1$  as follows. One starts with the root $O$ and connects it to $d$ ``children'' that we name $O_{(ch),j}, j=1,\ldots,d$. Then, we connect each $O_{(ch),j}$ to $d$ new vertices. Continuing this process generation by generation, we obtain an infinite rooted tree in which $O$ has $d$ neighbors, and any other vertex has $d+1$ neighbors. For each $Y\neq O$, the vertex $Y_{(p)}$ indicates its unique parent and $Y_{(ch),j}, j=1,\ldots,d$, its children. If $d=2$ and $Y\neq O$, we can define its unique \emph{sibling} $Y_{(s)}$ as the other child of $Y_{(p)}$. Given functions \( f \) and \( F \) on $\mathcal{V}$ and \( \mathcal E \), respectively, we shall denote by \( f_Y \) the value of \( f \) at \( Y \)  and by $F_{Z,Y}(=F_{Y,Z})$ the value of \( F\) at an edge \( Z\sim Y \).

Given a graph $\mathcal{G}=(\mathcal V,\mathcal E)$, let $V$, $W$, and $\sigma$ be functions on $\mathcal{V}$,  \( \mathcal E \), and $\vec{\mathcal{E}}$, respectively.  Assume that $V$ and $W$ are both bounded, $W>0$, and $\sigma$ takes value in $\{0,1\}$. By definition \( W_{Y,Z}=W_{Z,Y} \) while \( \sigma_{[Y,Z]} \) and \( \sigma_{[Z,Y]} \) might not be equal to each other. If there is a constant $C$ such that each vertex has at most $C$ neighbors, we can define an operator, a \emph{generalized Jacobi matrix on the graph $\mathcal{G}$}, by
\begin{equation}
\label{sad1}
(\jo f)_Y\ddd  V_Yf_Y+\sum_{Z\sim Y}(-1)^{\sigma_{[Y,Z]}}W_{Y,Z}^{1/2}f_Z,
\end{equation}
where \( f \) is any function on \( \mathcal V \). We call \( \jo \) a generalized Jacobi matrix since in most of the literature it is common to define \( \jo \) with  $\sigma\equiv0$. We, however,  allow a more general setup, which, as we explain later, is more natural in the case of Jacobi matrices generated by multiple orthogonality. Keeping this distinctions in mind, throughout the paper we call \( \jo \) from \eqref{sad1} simply a Jacobi matrix on $\mathcal{G}$.

As we already mentioned, we are interested in the connection between Jacobi matrices on graphs and orthogonal polynomials. In full generality of definition \eqref{sad1} such connection no longer exists. However, there are large classes of Jacobi operators on trees that can be defined via multiple orthogonal polynomials.  Spectral theory of Jacobi matrices and Schr\"odinger operators on trees is a vibrant topic of modern mathematical physics, see, e.g., \cite{aizen,sb1,t5,BrDenEl18,t3,t2,t4,t6}. It is conceivable that the powerful tools developed for the analysis of multiple orthogonality, already known to have applications in number theory, statistics, and random matrices, can find new  applications in the analysis of quantum systems.

\section*{Multiple orthogonal polynomials}

The system of polynomials orthogonal on the real line can be generalized to the case of orthogonality with respect to several measures. This multiple orthogonality, being a classical area of approximation theory, has connections to number theory, numerical analysis, etc., see \cite{AI,ismail,VAA} for the introduction to this topic. To define it, consider
\[
\vec{\mu}\ddd (\mu_1,\mu_2), \quad \supp \mu_k \subseteq \R, \quad \text{and} \quad \n\ddd (n_1,n_2)\in\Z_+^2,\quad |\n|\ddd n_1+n_2,
\]
where we assume that all the moments of the measures \( \mu_1,\mu_2 \) are finite.

\begin{Def}
Polynomials $A_\n^{(1)}(x)$ and $A_\n^{(2)}(x)$, \( \deg A_\n^{(k)}\leqslant n_k-1 \), \( k\in\{1,2\} \), that satisfy
\begin{eqnarray}
\label{sad15}
\int_{\mathbb{R}}x^m\big(A_\n^{(1)}(x)d\mu_1(x)+A_\n^{(2)}(x)d\mu_2(x)\big)=0, \quad m\in \{0,\ldots, |\n|-2\}\,,
\end{eqnarray}
are called type I multiple orthogonal polynomials (type I MOPs). We assume that $A_\n^{(k)}(x)\not\equiv0$ unless $n_k-1<0$. Furthermore, non-identically zero polynomial $P_\n(x)$ is called type II multiple orthogonal polynomial (type II MOP) if it  satisfies
\begin{eqnarray}
 \label{1.2}
 \deg P_\n\leqslant|\n|, \quad \int_{\mathbb{R}} P_{\n}(x)x^md\mu_k(x)=0 \quad  {\rm\it for~all~~} m\in \{0,\ldots, n_k-1\} {\rm\it~~and~~} k\in \{1,2\}.
\end{eqnarray}
\end{Def}

Polynomials of the first and second types always exist. The question of uniqueness is more involved. If every $ P_{\n}(x)$ has degree exactly $|\n|$,  then the multi-index $\n$ is called \emph{normal} and we choose the following normalization
\[
 P_{\n}(x)=x^{|\n|}+\cdots\,,
\]
i.e., the polynomial $ P_{\n}(x)$ is monic. It turns out that $\n$ is normal  if and only if the following linear form
\begin{equation}
\label{Qn}
Q_{\n}(x)\ddd A_\n^{(1)}(x)d\mu_1(x)+A_\n^{(1)}(x)d\mu_2(x)
\end{equation}
is defined uniquely up to multiplication by a constant. In this case \( \deg A_\n^{(k)}= n_k-1 \) and we will normalize the polynomials of the first type by
\begin{equation}
\label{n_2}
\int_{\mathbb{R}} x^{|\n|-1}Q_{\n}(x)=1\,.
\end{equation}

\begin{Def}
The vector $\vec{\mu}$  is called perfect if all the multi-indices $\n\in \Z_+^2$ are normal.
\end{Def}

Besides the orthogonal polynomials, we will need  the functions of the second kind. 

\begin{Def}
\label{def:1.3}
The functions
\begin{equation}
\label{Ln}
L_{\n}(z) \ddd  \int_{\mathbb{R}}\frac{Q_{\n}(x)}{z-x} \quad \text{ and }\quad R_\n^{(k)}(z)\ddd \int_\R\frac{P_{\n}(x)d\mu_k(x)}{z-x}, \quad k\in \{1,2\},
\end{equation}
are called functions of the second kind associated to the linear forms \( Q_{\vec n}(x) \) and to polynomials \(P_{\n}(x) \), respectively.
\end{Def}

If $d=1$, type II polynomials $P_{\n}(x)$ are the standard monic polynomials orthogonal on the real line with respect to the measure $\mu_1$ and the polynomials $A_\n^{(1)}(x)$ are proportional to $p_{n-1}(x,\mu_1)$ with the coefficient of proportionality that can be computed explicitly.

In the literature on orthogonal polynomials, the following Cauchy-type integral
\begin{equation}
\label{markov}
\widehat{\mu}(z)\ddd \int_{\mathbb{R}}\frac{d\mu(x)}{z-x}\,, \quad z\not\in\supp\,\mu\,,\quad  \mu\in\mathfrak M,
\end{equation}
is often referred to as a Markov function. If \( \mu_1,\mu_2\in\mathfrak M \), we can rewrite \(  L_\n(z) \) as 
\begin{equation}
\label{Ln-An}
L_\n(z) = A_\n^{(1)}(z)\widehat\mu_1(z) + A_\n^{(2)}(z)\widehat\mu_2(z) - A_\n^{(0)}(z),
\end{equation}
where \( A_\n^{(0)}(z) \) is a polynomial given by
\begin{equation}
\label{An0}
A_\n^{(0)}(z) \ddd \int_\R \frac{ A_\n^{(1)}(z)- A_\n^{(1)}(x)}{z-x}d\mu_1(x) + \int_\R \frac{ A_\n^{(2)}(z)- A_\n^{(2)}(x)}{z-x}d\mu_2(x).
\end{equation}

Similarly to classical orthogonal polynomials on the real line, the above MOPs also satisfy nearest-neighbor lattice recurrence relations. Denote by  $\vec{e}_1\ddd (1,0)$ and $\vec{e}_2\ddd (0,1)$ the standard basis vectors in \( \R^2 \). Assume that
\begin{equation}
\label{as-1}
\vec{\mu}=(\mu_1,\mu_2)\;\;\mbox{is perfect}\,.
\end{equation}
This is an assumption we carry through out the paper. In this case, see, e.g., \cite{ismail,VA11}, there exist real constants \( \{a_{\n,1},a_{\n,2},b_{\n,1},b_{\n,2}\}_{\n\in\Z_+^2}\), which we call the \emph{recurrence coefficients} corresponding to the system \( \vec\mu \), such that linear forms $Q_{\n}(x)$ satisfy
\begin{equation}
\label{1.5}
xQ_{\n}(x)=Q_{\n-\vec{e}_i}(x) + b_{\n-\vec{e}_i,i}Q_{\n}(x) + a_{\n,1}Q_{\n+\vec{e}_1}(x) + a_{\n,2}Q_{\n+\vec{e}_2}(x)\,, \quad \n\in \mathbb{N}^2,
\end{equation}
for each \( i\in\{1,2\} \), while it holds for type II polynomials that
\begin{equation}\label{1.7}
xP_{\n}(x)=P_{\n+\vec{e}_i}(x)+b_{\n,i}P_{\n}(x) + a_{\n,1}P_{\n-\vec{e}_1}(x)+ a_{\n,2}P_{\n-\vec{e}_2}(x)\,, \quad \n\in \Z_+^2,
\end{equation}
again, for each \( i\in\{1,2\} \), where we let  $P_{\n-\vec{e}_l}(x)\equiv0$ when the \( l \)-th components of $\n-\vec{e}_l$ is negative. It is known that
\begin{equation}
\label{coef-a}
a_{\n,i} \neq0,~~ \n\in\mathbb N^2, ~~i\in\{1,2\}, \quad\text{and}\quad
\left\{
\begin{array}{rl}
a_{(n,0),1},~a_{(0,n),2}>0, & n\in\mathbb N, \smallskip \\
a_{(0,n),1}=a_{(n,0),2}\ddd0, & n\in\Z_+,
\end{array}
\right.
\end{equation}
where the first conclusion follows from perfectness and an explicit integral representation for \( a_{\n,i} \), see~\cite[Equation~(1.8)]{VA11}, and the second one is part definition and part a consequence of positivity of parameters $\{a_n\}$ in \eqref{1.19}.

\begin{Rem}
For perfect systems $\vec\mu$, one can show that \eqref{1.5} implies the recursion for the type I polynomials themselves:
\begin{equation}
\label{3sd1}
xA_{\n}^{(j)}(x)=A_{\n-\vec{e}_i}^{(j)}(x) + b_{\n-\vec{e}_i,i}A_{\n}^{(j)}(x) + a_{\n,1}A_{\n+\vec{e}_1}^{(j)}(x) + a_{\n,2}A_{\n+\vec{e}_2}^{(j)}(x)\,, \, \n\in \mathbb{N}^2,\, i,j\in \{1,2\}\,.
\end{equation}
\end{Rem}

The recurrence coefficients  $\{a_{\n,i},b_{\n,i}\}$ are uniquely determined by $\vec{\mu}$. However, when $d>1$, unlike in the one-dimensional case, we can not prescribe them arbitrarily. In fact, coefficients in \eqref{1.5} and \eqref{1.7}  satisfy the so-called ``consistency conditions'', see, e.g., \cite[Theorem~3.2]{VA11} and \cite{ADVA}, which is a system of nonlinear difference equations:
 \begin{eqnarray*}
 b_{\n+\vec{e}_i,j}- b_{\n,j}=b_{\n+\vec{e}_j,i}- b_{\n,i},\\ 
 \sum_{k=1}^2 a_{\n+\vec{e}_j,k} -\sum_{k=1}^2 a_{\n+\vec{e}_i,k}=b_{\n+\vec{e}_j,i}b_{\n,j}-b_{\n+\vec{e}_i,j}b_{\n,i},\\
 a_{\n,i}(b_{\n,j}-b_{\n,i})=a_{\n+\vec{e}_j,i}(b_{\n-\vec{e}_i,j}-b_{\n-\vec{e}_i,i}),
 \end{eqnarray*}
 where $\n\in \mathbb{N}^2$ and $i,j\in \{1,2\}$. Conversely, see \cite[Theorem~3.1]{FilHanVA15}, solution to this nonlinear system is unique and uniquely defines $\vec{\mu}$ (\( \mu_k \)'s are the spectral measures of the Jacobi operators corresponding to the boundary values) provided the boundary values are properly defined. 
 

\part{Jacobi matrices on finite rooted trees}
\label{part1}

The goal of this part of the paper is to prove analogs of \eqref{1.21a} and \eqref{sar_1} for Jacobi matrices \eqref{sad1} on finite trees in the case when these Jacobi matrices are generated by multiple orthogonality.

\section{Definitions and basic properties}
\label{s:ft-def}

\subsection{Finite trees}
\label{ss:TN} Fix $\vec N=(N_1,N_2)\in \mathbb{N}^2$. Truncate $\Z_+^2$ to a discrete rectangle
\[
\mathcal{R}_{\vec N}=\{\n: n_1\leq N_1, n_2\leq N_2\}
\]
and denote by $\mathcal{P}_{\vec N}$ the family of all paths of length $|\vec N|=N_1+N_2$ connecting the points $\vec N=(N_1,N_2)$ and $(0,0)$  (within a path exactly one of the coordinates is decreasing by $1$ at each step).  
\begin{figure}[h!]
\includegraphics[scale=.6]{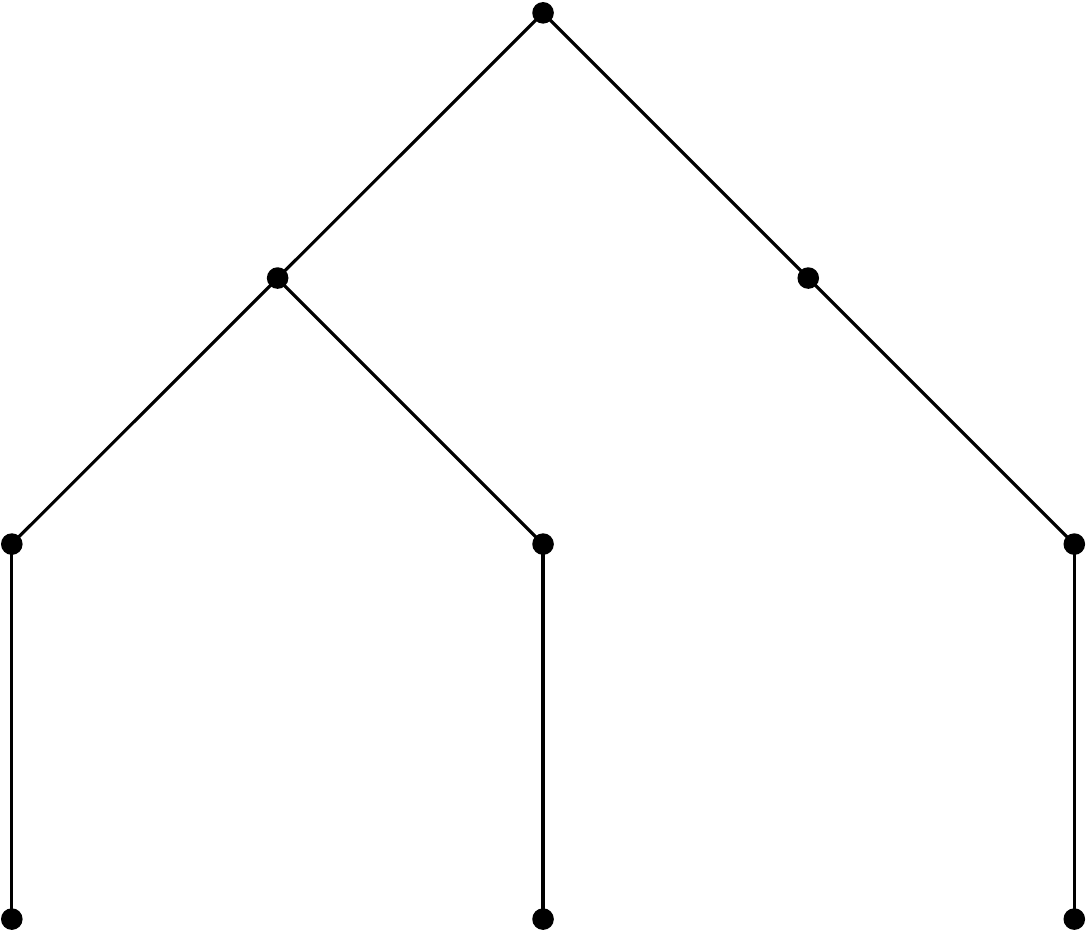}
\begin{picture}(0,0)
\put(-90,157){$(2,1)\sim O$}
\put(-137,110){$(1,1)\sim X_{(p)}=Y_{(p)}$}
\put(-45,110){$(2,0)\sim Z_{(p)}$}
\put(-185,63){$(0,1)\sim X=A_{(p)}$}
\put(-95,63){$(1,0)\sim Y=B_{(p)}$}
\put(-2,63){$(1,0)\sim Z=C_{(p)}$}
\put(-185,0){$(0,0)\sim A=X_{(ch),2}$}
\put(-95,0){$(0,0)\sim B=Y_{(ch),1}$}
\put(-2,0){$(0,0)\sim C=Z_{(ch),1}$}
\end{picture}
\caption{Tree for \( \vec N=(2,1) \).}
\label{fig:finite-tree}
\end{figure}
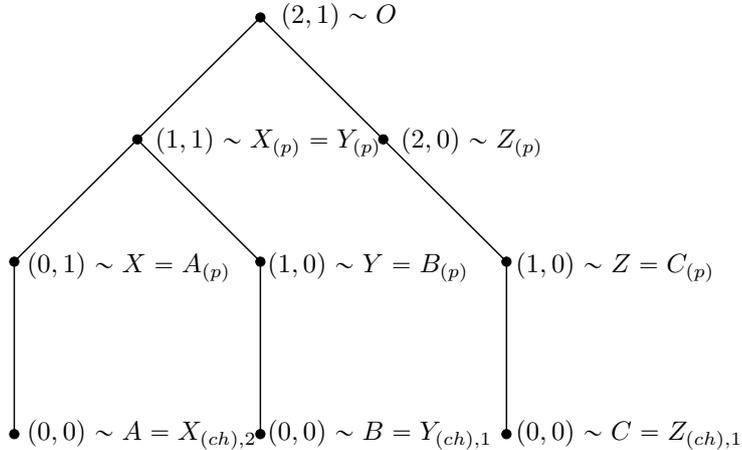
The tree \( \mathcal T_{\vec N} \) is obtained by untwining \( \cal P_{\vec N} \) in such a way that \( \cal P_{\vec N} \) is in one-to-one correspondence with the paths in \( \cal T_{\vec N} \) originating at the root, say \( O \), which corresponds to \( \vec N \), see Figure~\ref{fig:finite-tree} for  $\vec N=(2,1)$. 

We denote by $\cal{V}_{\vec N}$ the set of the vertices of \( \mathcal T_{\vec N} \). The above construction defines a projection $\Pi:\mathcal{V}_{\vec N}\to\mathcal{R}_{\vec N}$ as follows: given $Y\in \mathcal{V}_{\vec N}$ we consider the path from $O$ to $Y$, take the corresponding path on \( \mathcal{R}_{\vec N} \), and let $\Pi(Y)$ to be its the endpoint (the one which is not \( \vec N \)). We denote by \( \ell^2(\mathcal V_{\vec N}) \) the set of all functions on \( \mathcal V_{\vec N} \) with the norm coming from the standard inner product \( \langle\cdot,\cdot\rangle \).

As agreed before, we denote by \( Y_{(p)} \) the ``parent'' of \( Y \). To distinguish the ``children'' of a vertex \( Y \) we introduce an index function \( \iota  \) by
\begin{equation}
\label{iota}
\iota: \mathcal V_{\vec N} \to \{ 1,2\} \quad Z\mapsto \iota_Z ~~\text{ such that }~~ \Pi(Z_{(p)}) = \Pi(Z) + \vec e_{\iota_Z}.
\end{equation}
Then, if \( Y=Z_{(p)} \), we write \( Z=Y_{(ch),\iota_Z} \), see Figure~\ref{fig:finite-tree}. 
We further let
\[
ch(Y)\ddd \big\{i:~n_i>0,~\Pi(Y)=(n_1,n_2)\big\}
\]
to be the index set of the children of \( Y \). It will be convenient to introduce an artificial vertex \( O_{(p)} \), a formal parent of the root \( O \). We do not include \( O_{(p)} \) into \( \mathcal V_{\vec N} \), but we do extend every function \( f \) on \( \mathcal V_{\vec N} \) to \( O_{(p)} \) by setting \( f_{O_{(p)}} = 0 \) (recall that we denote the value of a function \( f \) at \( Y\in\mathcal V_{\vec N} \) by \( f_Y \)).

\subsection{Jacobi matrices generated by multiple orthogonality} \label{sub_1}

Let \( \vec \mu \) be a perfect system and \( \{ a_{\n,i},b_{\n,i} \} \) be its recurrence coefficients, see \eqref{1.5} and \eqref{1.7}. In this subsection, we specialize definition \eqref{sad1} to the case of finite trees \( \mathcal T_{\vec N} \) and Jacobi matrices whose potentials \( V,W \), and the signature $\sigma$ come from \( \vec\mu \).

Fix \( \vec\kappa\in\R^2 \) such that \( |\vec\kappa|=\kappa_1+\kappa_2=1 \). We define the potentials \( V=V^{\vec\mu},W=W^{\vec\mu}:\mathcal V_{\vec N} \to \R \) (as with most quantities depending on \( \vec \mu \), we drop the dependence on \( \vec\mu \) from notation) by
\begin{equation}
\label{potentials}
V_O \ddd \kappa_1 b_{\vec N,1} + \kappa_2 b_{\vec N,2}, \quad W_O\ddd1, \quad \text{and} \quad V_Y \ddd b_{\Pi(Y),\iota_Y},  \quad W_Y \ddd \big|a_{\Pi(Y_{(p)}),\iota_Y}\big|, ~~Y\neq O.
\end{equation}
This definition is consistent with \eqref{sad1} if we let \( W_{Y_{(p)},Y} = W_{Y,Y_{(p)}} = W_Y \) (for trees, neighboring vertices always form child/parent pairs). We further choose function \( \sigma:\mathcal V_{\vec N} \to \{0,1\} \) to recover the signs of the recurrence coefficients \( a_{\n,i} \). Namely, we set \( \sigma_Y \) to be such that
\begin{equation}
\label{nu-fun}
(-1)^{\sigma_Y}W_Y = a_{\Pi(Y_{(p)}),\iota_Y},~~Y\neq O,~~\text{and}~~\sigma_O\ddd0
\end{equation}
(observe that \( W_Y>0 \) since \( a_{\n+\vec e_i,i}\neq 0 \) by \eqref{coef-a}). To relate back to the definition given in \eqref{sad1}, we set \( \sigma_{[Y,Y_{(p)}]} = 0 \) and \( \sigma_{[Y_{(p)},Y]} = \sigma_Y \). With these definitions, \eqref{sad1} specializes to
\begin{equation}
\label{jackn}
(\jackn f)_Y \ddd V_Yf_Y + W_Y^{1/2}f_{Y_{(p)}} +  \sum_{l\in ch(Y)} (-1)^{\sigma_{Y_{(ch),l}}} W_{Y_{(ch),l}}^{1/2}f_{Y_{(ch),l}},
\end{equation}
which we call a \emph{Jacobi matrix on a finite tree} $\cal{T}_{\vec N}$.

Let $P_{\n}(z)$ be the type II MOPs corresponding to the multi-index \( \n \) with respect to \( \vec\mu \), see \eqref{1.2}. We consider \( z\in\C_+\) as a parameter and put
\begin{equation}
\label{gevn1}
p_Y(z) \ddd m_Y^{-1}P_Y(z), \quad P_Y(z) \ddd P_{\Pi(Y)}(z), \quad \text{and} \quad m_Y\ddd   \prod_{Z\in {\rm path}(Y,O)} W_Z^{-1/2},
\end{equation}
where  ${\rm path} (Y,O)$ is the non-self-intersecting path connecting $Y$ and $O$ that includes both $Y$ and \( O \). Obviously, all three functions \( p \), \( P \), and \( m \) depend on \( \vec \mu \). To uniformize the notation,  let us formally set
\begin{equation}
\label{POp}
P_{\Pi(O_{(p)})}(z) \ddd \kappa_1 P_{\vec N+\vec{e}_1}(z) + \kappa_2 P_{\vec N+\vec{e}_2}(z).
\end{equation}
Given \( X\in\mathcal V_{\vec N} \), denote by \( \mathcal T_{\vec N [X]} \) the subtree of \( \mathcal T_{\vec N} \) with root at \( X \) and by \( \mathcal V_{\vec N [X]}\) the set of its vertices. Let \( \mathcal{J}_{[X]} \) and \( p^{[X]} \) be the restriction of \( \jackn \) and \( p \) to \( \mathcal T_{\vec N [X]} \) and \( \mathcal V_{\vec N [X]} \), respectively. Then, it follows from \eqref{1.7} that
\begin{equation}
\label{gevn2}
\mathcal{J}_{[X]} p^{[X]}(z) = zp^{[X]}(z) - \big(m_X^{-1}P_{\Pi(X_{(p)})}(z)\big)\delta^{(X)},
\end{equation}
which is an  identity reminiscent of \eqref{1.21a}.

\subsection{Conditions on \( \vec\mu \)}

Recall that \( \vec\mu \) is a perfect system since, otherwise, its recurrence coefficients might not exist for all \( \n\in\mathcal R_{\vec N} \) and \( \jackn \) is undefined. Besides that, we place one more set of conditions on  \( \vec \mu \). Denote by \( E_{\Pi(Y)} \) the set of zeroes of \( P_{\Pi(Y)}(x) \), \( Y\in\mathcal V_{\vec N}\cup \{O_{(p)}\} \) (recall \eqref{POp}). Notice that \( E_{\Pi(O_{(p)})} =E_{\vec N+\vec e_i}\) when \( \vec\kappa=\vec e_i \), \( i\in\{1,2\} \). Our additional assumptions on $\vec \mu$ are 
\begin{equation}
\label{as-4}
\left\{
\begin{array}{rl}
E_{\Pi(Y)} \subset \R,~~ \#E_{\Pi(Y)} = |\Pi(Y)|, & Y\in\mathcal V_{\vec N}\cup \{O_{(p)}\}, \medskip \\
E_{\Pi(Y)} \cap E_{\Pi(Y_{(p)})} = \varnothing, & Y\in\mathcal V_{\vec N},
\end{array}
\right.
\end{equation}
where we put \( |\Pi(O_{(p)})| \ddd |\vec N| +1\) and \( \#S \) denotes the cardinality of \( S \). That is, we assume that all  zeroes of the polynomials \( P_{\Pi(Y)}(x) \) are real and simple, and that \( P_{\Pi(Y)}(x) \) and \( P_{\Pi(Y_{(p)})}(x) \) do not have common zeroes.

All the classical examples of type II MOPs satisfy \eqref{as-4}. Indeed, for  Angelesco systems, see Part~\ref{part3} further below, multiple Hermite polynomials \cite[Section~5.1]{VA11}, multiple Laguerre polynomials of the second kind \cite[Section~5.4]{VA11}, multiple Charlier polynomials \cite[Section~5.2]{VA11}, and multiple Meixner polynomials of the first kind \cite[Section~3.3]{as1}, it holds that
\begin{equation}
\label{as-2}
a_{\n,i}>0 \quad  \n \in \mathbb N^2,~~i\in\{1,2\}.
\end{equation}
This, together with perfectness (all the above examples form perfect systems) implies, see \cite[Theorem~2.2]{as1}, that
\begin{equation}
\label{interlacing}
x_{\n+\vec e_i,1}<x_{\n,1}<x_{\n+\vec e_i,2}<x_{\n,2}<\ldots<x_{\n,|\n|}<x_{\n+\vec e_i,|\n|+1}
\end{equation}
for any \( i\in\{1,2\} \), where we write \( E_\n=\{x_{\n,1},\ldots,x_{\n,|\n|}\} \). That is, the zeroes of  \( P_\n(x) \) and  \( P_{\n+\vec e_i}(x) \) interlace. Hence, the only conditions that remain to be checked in \eqref{as-4} are those that involve \( O_{(p)} \) and they, of course, depend on \( \vec\kappa \). The positivity of \( a_{\n,i} \), i.e., the condition \eqref{as-2},  is not satisfied by other classical systems such as Nikishin systems, see Section~\ref{s:ap1} further below, 
multiple Laguerre polynomials of the first kind \cite[Section~5.3]{VA11}, Jacobi-Pi\~neiro polynomials \cite[Section~5.5]{VA11}, and multiple Meixner polynomials of the second kind \cite[Section~3.7]{as1}.  However, it is known that type II MOPs form the so-called AT-systems and  their zeroes again satisfy  \eqref{interlacing}   for all just listed examples, see~\cite{as1}. Hence,  all  conditions in \eqref{as-4},  except for the ones involving \( O_{(p)} \), are satisfied automatically.

\section{Spectral analysis}

\subsection{Spectrum and eigenvalues}
One can readily see from \eqref{gevn2} that every \( E\in E_{\Pi(O_{(p)})} \) is an eigenvalue and
\begin{equation}
\label{bEO}
\jackn b(E,O_{(p)}) = Eb(E,O_{(p)}), \quad b(E,O_{(p)}) \ddd p(E).
\end{equation}
We call \( b(E,O_{(p)}) \) the \emph{trivial canonical eigenvector}. To identify the remaining eigenvalues and eigenvectors, we set
\begin{equation}
\label{EkN}
\mathcal E_{\vec\kappa,\vec N} \ddd E_{\Pi(O_{(p)})}\cup \bigcup_{Y\in\mathcal V_{\vec N}:~\#ch(Y)=2} E_{\Pi(Y)}.
\end{equation}
The condition \( \#ch(Y)=2 \) is equivalent to \( \Pi(Y)\in\mathbb N^2 \). Hence, the set \( \mathcal E_{\vec \kappa,\vec N} \) consists of $E_{\Pi(O_{(p)})}$ and the zeroes of type II MOPs that are ``truly'' multiple orthogonal, i.e., they satisfy orthogonality conditions on both intervals. Given $E\in \mathcal E_{\vec \kappa,\vec N}$, let \( \jnt \) be the set of {\it joints} corresponding to \( E \) defined by
\begin{equation}
\label{jnt}
\jnt\ddd \big\{Y\in \mathcal V_{\vec N}: \, P_{Y}(E)=0 \text{~~and~~}  \#{  ch}(Y)=2\big\}.
\end{equation}
If \( E\in E_{\Pi(O_{(p)})} \) and 
\( E\notin \bigcup_{Y\in\mathcal V_{\vec N}:~\#ch(Y)=2} E_{\Pi(Y)}\), then \( \jnt=\varnothing \); otherwise,  \( \jnt\neq\varnothing \).  To each \( X\in\jnt \), we associate a special vector. To define it, recall that \( W_Y>0 \) for all \( Y \), see the remark after formula \eqref{nu-fun}, and that \( p_{X_{(ch),l}}(E)\neq 0 \) by \eqref{as-4} when \( X\in\jnt \). We will need a standard  notation: if $\cal{B}$ is a subset of a graph $\cal{G}$, the symbol $ \chi_{\cal{B}}$ denotes its characteristic function. Given $E\in \mathcal E_{\vec\kappa,\vec N}$, \( X\in\jnt \), let
\begin{equation}
\label{bEX}
b(E,X) \ddd p(E) \left( \frac{(-1)^{\sigma_{X_{(ch),2}}}\chi_{\cal{T}_{\vec N[X_{(ch),2}]}}}{W_{X_{(ch),2}}^{1/2}p_{X_{(ch),2}}(E)} - \frac{(-1)^{\sigma_{X_{(ch),1}}}\chi_{\cal{T}_{\vec N[X_{(ch),1}]}}}{W_{X_{(ch),1}}^{1/2}p_{X_{(ch),1}}(E)} \right) ,
\end{equation}
where, as before, \( \mathcal T_{\vec N[Z]} \) denotes the  subtree of \( \mathcal T_{\vec N} \) with root at \( Z \). Anticipating the forthcoming theorem, we call each \( b(E,X) \) a \emph{canonical eigenvector} (it follows right away from \eqref{jackn} that \( \jackn b(E,X) \) is also supported on \( \cal{T}_{\vec N[X_{(ch),1}]}\cup\cal{T}_{\vec N[X_{(ch),2}]} \)). Finally, we set
\begin{equation}
\label{jnts}
\jnts \ddd \left\{ \begin{array}{ll}\jnt, & E\not\in E_{\Pi(O_{(p)})}, \smallskip \\  \jnt \cup\{O_{(p)}\}, & E\in E_{\Pi(O_{(p)})}. \end{array}\right.
\end{equation}
Definitions \eqref{bEO}, \eqref{EkN}, \eqref{bEX}, and \eqref{jnts} are needed for the following theorem, which is the main result of this part.

\begin{Thm}
\label{sde_1} 
Let \( \vec\mu \) be a perfect system of measures on the real line for which \eqref{as-4} holds and \( \jackn \) be the corresponding Jacobi matrix defined in \eqref{jackn}. Then
\[
\sigma(\jackn) = \mathcal E_{\vec \kappa,\vec N}.
\]
Given \( E \in \sigma(\jackn)\), a particular basis of the eigenspace corresponding to \( E \) is given by
\[
\big\{b(E,X):~X\in\jnts\big\}
\] 
and the geometric multiplicity of $E$, we call it $g_E$, is given by 
\[ 
g_E=\#\jnts.
\]
Moreover, the system 
\[
\big\{b(E,X):~X\in\jnts,~~E\in\sigma(\jackn)\big\}
\]
is a basis for \( \ell^2(\mathcal V_{\vec N}) \).
\end{Thm}

We illustrate the construction of the canonical eigenvectors for a simple case of $\jo_{\vec e_2,(2,1)}$, see Figure~\ref{fig:finite-tree}. There are 9 vertices and 9 eigenvalues: $E_{(2,2)}$ has 4 roots, $E_{(2,1)}$ provides 3 roots, and $E_{(1,1)}$ gives 2 more. Assume that {\it all 9 zeroes are distinct}. Then, $\jnts=\{O_{(p)}\}$ for the roots \( E \) in $E_{(2,2)}$ and each such root defines a trivial canonical eigenvector $p(E)$. Every $E$ in $E_{(2,1)}$ is a simple eigenvalue and we have that \( \jnts = \{O\} \). The corresponding canonical eigenvector $b(E,O)$ is equal to zero at $O$, has the same values as \( p(E)/(W_{X_{(p)}}^{1/2}p_{X_{(p)}}(E)) \) and $p(E)/(W_{Z_{(p)}}^{1/2}p_{Z_{(p)}}(E))$ at the vertices $X_{(p)},X,Y,A,B$ and vertices $Z_{(p)},Z,C$, respectively.  Finally, consider $E\in E_{(1,1)}$ for which \( \jnts = \{X_{(p)}\} \). The canonical eigenvector $b(E,X_{(p)})$ is zero at points $O,X_{(p)},Z_{(p)},Z,C$. Its values at $X,A$ are equal to the ones of \( p(E)/(W_X^{1/2}p_X(E)) \) and its values at $Y,B$ coincide with the values of \(p(E)/(W_Y^{1/2}p_Y(E)) \) there.

\subsection{\(\mathfrak S\)-self-adjointness}

When \( \sigma\equiv0 \) in \eqref{nu-fun}, or equivalently, \eqref{as-2} holds, the corresponding Jacobi matrix is self-adjoint and thus has an orthogonal basis of eigenvectors. When \( \sigma\not\equiv0 \) this is no longer the case. However, there exists an indefinite inner product given by a diagonal matrix \(\mathfrak S \) with diagonal entries equal  \( \pm1 \) such that Jacobi matrices are \( \mathfrak S \)-self-adjoint. The general theory of \( \mathfrak S \)-self-adjoint operators (see, e.g., \cite{gb}) does not guarantee that their eigenvectors span \(  \ell^2(\mathcal V_{\vec N}) \) (that is, that \( \jackn \) has no Jordan blocks, i.e., that it has a \emph{simple structure}). Yet, this is indeed the case for Jacobi matrices.

Let, as before, ${\rm path} (Y,O)$ be the non-self-intersecting path connecting $Y$ and $O$ that includes both $Y$ and \( O \).  Define a diagonal matrix \( \mathfrak S \) on \( \mathcal T_{\vec N} \) by
\begin{equation}
\label{matS}
\mathfrak{S}\delta^{(O)}\ddd \delta^{(O)} \quad \text{and} \quad \mathfrak{S}\delta^{(Y)}\ddd (-1)^{\sum_{Z\in{\rm path}(Y,O)}\sigma_Z}\delta^{(Y)}, ~~~ Y\neq O.
\end{equation}
The diagonal matrix $\mathfrak{S}$ defined this way assigns either $+1$ or $-1$ to a vertex $Y$  depending on whether the number of ``negative'' edges connecting $O$ to $Y$ is even or odd. We define an indefinite inner product \( [\cdot,\cdot] \) by
\begin{equation}
\label{iip}
[f,g]\ddd \big\langle\mathfrak{S}f,g\big\rangle, \quad f,g\in\ell^2(\mathcal V_{\vec N}).
\end{equation}
Denote  the number of vertices $Y\in\mathcal V_{\vec N}$ such that $[\delta^{(Y)},\delta^{(Y)}]=\pm1$ by \( i_\pm \). If \( \sigma\equiv0 \), the matrix \( \mathfrak S \) is the identity matrix and \( [\cdot,\cdot]=\langle\cdot,\cdot\rangle\), \( i_+=\#\mathcal V_{\vec N} \) while \( i_-=0 \). We let $\ell^2_\mathfrak{S}(\mathcal V_{\vec N})$ denote the corresponding indefinite inner product vector space, which is sometimes called a finite-dimensional Krein space. 

A matrix $\mathcal A$ is called $\mathfrak{S}$-self-adjoint if
\begin{equation}
\label{Ssa}
[\mathcal Af,g]=[f,\mathcal Ag] 
\end{equation}
for all vectors $f$ and $g$.
 Notice that \eqref{Ssa} is equivalent to $\mathfrak{S}\mathcal A=\mathcal A^*\mathfrak{S}$, where $\mathcal A^*$ is the adjoint of $\mathcal A$ in the original inner product $\langle\cdot,\cdot\rangle$. 
Since $\mathfrak{S}^2$ is the identity matrix, multiplying identity $\mathfrak{S}\mathcal A=\mathcal A^*\mathfrak{S}$ from the left and from the right by $\mathfrak{S}$ gives us \( \mathcal A\mathfrak S=\mathfrak S\mathcal A^*\). Thus, $\mathcal A$ is $\mathfrak{S}$-self-adjoint if and only if  $\mathcal A^*$ is $\mathfrak{S}$-self-adjoint. Clearly, when $\mathfrak{S}$ is the identity matrix, i.e., when \eqref{as-2} holds, condition \eqref{Ssa} is equivalent to $\mathcal A$ being self-adjoint in the standard inner product.

\begin{Prop} 
\label{sadj}
Jacobi matrices $\jackn$ and  $\jackn^*$ are $\mathfrak{S}$-self-adjoint.
\end{Prop}

\subsection{\(\mathfrak S\)-orthogonalization}
\label{ss:123}

In this subsection, we show that the basis of canonical eigenvectors, which is yielded by Theorem \ref{sde_1}, can be used to construct $\mathfrak{S}$-orthogonal basis of eigenvectors. To this end, we notice that eigenspaces that correspond to two different real eigenvalues are already $\mathfrak{S}$--orthogonal. Indeed, this is due to the following identity
 \[
E_1[\Psi_1,\Psi_2]=[\jackn\Psi_1,\Psi_2]=[\Psi_1,\jackn\Psi_2]=[\Psi_1,E_2\Psi_2]=E_2[\Psi_1,\Psi_2],
\]
where \( E_1,E_2 \) are eigenvalues of \( \jackn \) and \( \Psi_1,\Psi_2 \) are corresponding eigenvectors. Thus, we only need to focus on each individual eigenspace.

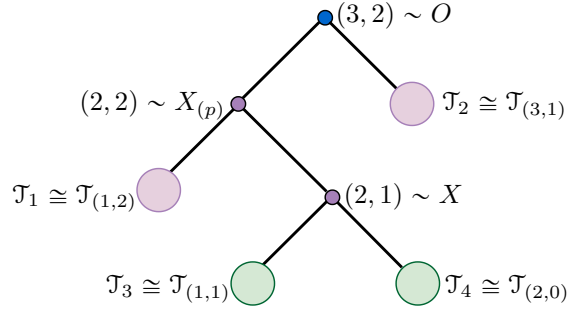
\begin{figure}[h!]
\includegraphics[scale=.75]{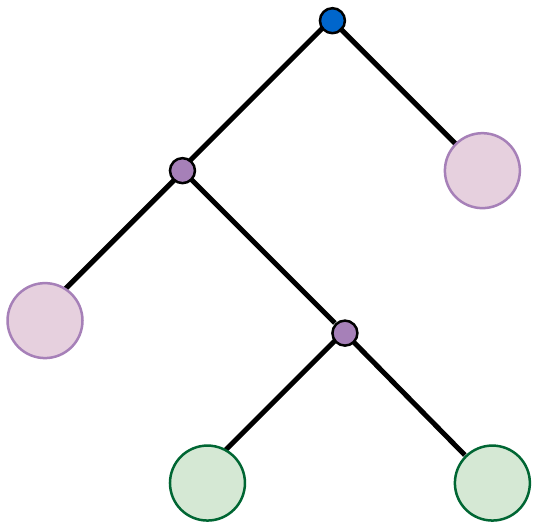}
\begin{picture}(0,0)
\put(-44,108){$(3,2)\sim O$}
\put(-140,75){$(2,2)\sim X_{(p)}$}
\put(-165,40){$\mathcal T_1 \cong \mathcal T_{(1,2)}$}
\put(-41,40){$(2,1)\sim X$}
\put(-4,75){$\mathcal T_2 \cong \mathcal T_{(3,1)}$}
\put(-130,6){$\mathcal T_3 \cong \mathcal T_{(1,1)}$}
\put(-3,6){$\mathcal T_4 \cong \mathcal T_{(2,0)}$}
\end{picture}
\caption{Partition of \( \mathcal V_{(3,2)} \) into waves \( \mathcal W_1(E) \) (blue), \( \mathcal W_2(E) \) (purple), and \( \mathcal W_3(E) \) (green) when \( \jnt=\{O,X\} \).}
\label{fig:partition}
\end{figure}

Suppose $E$ is an eigenvalue and $\jnt\neq \emptyset$. That guarantees that $g(E)>1$ and $\{b(E,X)\}$ is the basis in the eigenspace. We start with some geometric constructions on the tree and a few definitions. 
Let us first partition \( \mathcal V_{\vec N} \) into a collection of disjoint ``waves''. Define the canopy of \( \mathcal T_{\vec N}\) by \( \mathcal C\ddd\Pi^{-1}(0,0) \). If \( O\in\jnt \), we set the first \emph{wave} and its \emph{front} simply to be \( \{O\} \), that is, $\mathcal{W}_1(E)=\mathcal{F}_1(E)=\{O\}$. Otherwise, we define $\mathcal{F}_1(E)$ to be the set of vertices from $\mathcal{C}\cup\jnt$ that can be connected to $O$ by a path which does not contain elements of $\jnt$ in its interior. We then let the wave $\mathcal{W}_1(E)$ to be the union of all the vertices on these paths, including the endpoints.  To define $\mathcal{F}_2(E)$, consider all the vertices in  $(\mathcal C \cup \jnt)\setminus \mathcal W_1(E)$ that can be connected to a vertex in $\mathcal{F}_1(E)$ by a path which does not contain vertices of $\jnt$ in its interior. The second wave $\mathcal{W}_2(E)$ is then defined as the set of all the vertices on these paths,  including the ones from $\mathcal{F}_2(E)$, but excluding the ones from $\mathcal{F}_1(E)$ (so, \( \mathcal W_1(E) \cap \mathcal W_2(E) = \varnothing \)). We continue this process until all of $\mathcal V_{\vec N}$ is exhausted.

\begin{Exa}
Consider \( \mathcal T_{(3,2)} \) and assume that \( \jnt=\{O,X\} \), where \( \Pi(X) = (2,1) \), see Figure~\ref{fig:partition}. Then, 
\[
\mathcal W_1(E) = \{O\}, \quad \mathcal W_2(E)=\{X_{(p)},X\}\cup \mathcal V(\mathcal T_1)\cup \mathcal V(\mathcal T_2), \quad \text{and} \quad \mathcal W_3(E) = \mathcal V(\mathcal T_3)\cup \mathcal V(\mathcal T_4),
\]
where \( \mathcal T_1 \), \( \mathcal T_2 \), \( \mathcal T_3 \), and \( \mathcal T_4 \) are the subtrees with the roots at the sibling of \( X \), the sibling of \( X_{(p)} \),  \( X_{(ch),1} \), and \( X_{(ch),2} \), respectively, and \( \mathcal V(\mathcal T) \) is the set of vertices of a subtree \( \mathcal T \). Moreover, it holds that
\[
\mathcal F_1(E) = \{ O \}, \quad \mathcal F_2(E) = \{ X \} \cup \big( \mathcal C \cap (\mathcal V(\mathcal T_1)\cup \mathcal V(\mathcal T_2)) \big), \quad \text{and} \quad \mathcal F_2(E) = \mathcal C \cap (\mathcal V(\mathcal T_3)\cup \mathcal V(\mathcal T_4)).
\]  
\end{Exa}

Suppose all constructed fronts and waves are enumerated by $\{\cal{F}_1,\ldots,\cal{F}_p\}$ and $\{\cal{W}_1,\ldots,\cal{W}_p\}$. To produce $\mathfrak{S}$-orthogonal basis out of $\{b(E,X)\}$, we start at the canopy and go up the tree. Consider the canonical eigenvectors corresponding to \( E \) that are supported inside the last wave $\mathcal W_p(E)$. Each of these eigenvectors has support on a subtree sitting inside $\mathcal{W}_p(E)$ and having the root at a vertex of the previous front $\mathcal{F}_{p-1}(E)$. As their supports are disjoint, they are $\mathfrak{S}$-orthogonal. Call their span $\mathcal S_p(E)$.   Next, take all the canonical eigenvectors that have support inside $\mathcal{W}_{p-1}(E)\cup\mathcal{W}_{p}(E)$ and that were not chosen before. For each of them, take its $\mathfrak{S}$-perpendicular to $\mathcal S_p(E)$. By construction, it is nonzero. These new vectors are still eigenvectors and they are $\mathfrak{S}$-orthogonal to each other because they are supported on different subtrees as well as $\mathfrak{S}$-orthogonal to the previously considered eigenvectors  by constructions. Denote by \( \mathcal S_{p-1}(E) \) the span of these $\mathfrak{S}$-perpendiculars and previously considered eigenvectors spanning $\mathcal S_p(E)$. If we continue going up the tree in this fashion, we will produce an  $\mathfrak{S}$-orthogonal basis of the $E$-eigenspace. Since all eigenspaces are $\mathfrak{S}$-orthogonal, we have constructed a $\mathfrak{S}$-orthogonal set of eigenvectors. By scaling, we can make sure that this basis is $\mathfrak{S}$-orthonormal.

We want to finish by explaining how our result fits into the general spectral theory of $\mathfrak{S}$-self-adjoint operators. We say that a vector $\psi$ is $\mathfrak{S}$-positive if $[\psi,\psi]>0$ and $\mathfrak{S}$-negative if $[\psi,\psi]<0$. It is  $\mathfrak{S}$-neutral if $[\psi,\psi]=0$. Suppose $\{\psi_1,\ldots,\psi_n\}$ is a $\mathfrak{S}$-orthogonal basis of \( \ell^2(\mathcal V_{\vec N}) \). It known, see \cite[Proposition~2.2.3]{gb} and Lemma~\ref{lem:7} further below, that
\[ 
\#\big\{j:~\psi_j~\text{is $\mathfrak{S}$-negative}\big\}=i_- \quad \text{and} \quad \#\big\{j:~\psi_j~\text{is $\mathfrak{S}$-positive}\big\}=i_+\,,
\]
where the numbers \( i_\pm \) were defined right after \eqref{iip}. Enumerate  the $\mathfrak{S}$-positive and  $\mathfrak{S}$-negative vectors in the basis $\{\psi_1,\ldots,\psi_n\}$ by \( \{\psi_1^+,\ldots,\psi_{i_+}^+\} \) and by 
\( \{\psi_1^-,\ldots,\psi_{i_-}^-\} \), respectively. We clearly have a $\mathfrak{S}$-orthogonal sum decomposition
\[
\ell^2_\mathfrak{S}(\mathcal V_{\vec N})=H_+\oplus_{\mathfrak{S}} H_-, \quad H_\pm=\text{span}\big\{\psi^\pm_1,\ldots,\psi^\pm_{i_\pm}\big\},
\]
where $H_+$ and $H_-$ are positive and negative subspaces. In the  case of our $\mathfrak{S}$-self-adjoint Jacobi matrices, we just illustrated that  such a basis $\{\psi_1,\ldots,\psi_n\}$ can be built out of canonical eigenvectors.  That provides the concrete realization of the Spectral Theorem for $\mathfrak{S}$-self-adjoint matrices, see, e.g.,  \cite[Theorem~5.1.1]{gb}.

\section{Proofs of the main results}

\begin{proof}[Proof of Proposition~\ref{sadj}]
By formula \eqref{Ssa} and the remark that comes after it, we need to check that  \( \mathfrak{S}\jackn = \jackn^*\mathfrak{S} \), which is the same as checking 
\[
\big\langle f, \mathfrak{S}\jackn g\big\rangle = \big\langle \jackn f, \mathfrak{S} g\big\rangle
\]
for all vectors $f,g\in\ell^2(\mathcal V_{\vec N})$. Since $\jackn$ only contains self-interaction and interaction between neighbors, it is enough to consider cases $f=\delta^{(Z)}$ and $g=\delta^{(X)}$ where either \( Z=X \) or $Z\sim X$. It follows from \eqref{jackn} that
\[
\jackn\delta^{(X)} = (-1)^{\sigma_X}W_X^{1/2}\delta^{(X_{(p)})} + V_X\delta^{(X)} + \sum_{l\in ch(X)} W_{X_{(ch),l}}^{1/2}\delta^{(X_{(ch),l})},
\]
where we agree that \( \delta^{(O_{(p)})}\equiv0 \). It further follows from \eqref{matS} that
\[
\mathfrak S\jackn\delta^{(X)} = \big[\delta^{(X)},\delta^{(X)}\big]\bigg( W_X^{1/2}\delta^{(X_{(p)})} + V_X\delta^{(X)} + \sum_{l\in ch(X)} (-1)^{\sigma_{X_{(ch),l}}}W_{X_{(ch),l}}^{1/2}\delta^{(X_{(ch),l})}\bigg).
\]
Now, it is a simple matter of examining  three cases: when \( Z=X \), \( Z=X_{(p)} \), and \( Z=X_{(ch),l} \).
\end{proof}

It will be convenient for us to split the proof Theorem~\ref{sde_1} into several lemmas.  Let $X_{(g)}$ denote the parent of $X_{(p)}$. Recall that we extend all  functions on \( \ell^2(\mathcal V_{\vec N}) \) to \( O_{(p)} \) by zero.

\begin{Lem}
\label{lem:1}
Let $E\in\sigma(\jackn)$ and $\Psi$ be a corresponding eigenvector. If \( \Psi_X\neq 0 \) and \( \Psi_{X_{(p)}} = 0 \), then \( E \in E_{\Pi(X_{(p)})} \). Moreover, if we also have \( \Psi_{X_{(g)}}=0 \), then \( X_{(p)}\in\jnt \). Finally, we have an inclusion
\[
 \sigma(\jackn)\subseteq\mathcal E_{\vec \kappa,\vec N}.
 \]
\end{Lem}
\begin{proof}
Denote by \( \mathcal T_{[X]} \) the subtree of \( \mathcal T_{\vec N} \) with root at \( X \) and by \( \jo_{[X]} \) the restriction of \( \jackn \)  to  \( \mathcal T_{[X]} \). By the conditions of the lemma, \( E \) is also an eigenvalue of \( \jo_{[X]} \) with an eigenvector \( \chi_{\mathcal T_{[X]}}\Psi  \).  We can restrict the indefinite inner product to \( \mathcal T_{[X]} \) as well keeping the same notation $[\cdot,\cdot]$. Notice that 
$\jo_{[X]}=\chi_{\mathcal T_{[X]}}\jackn\chi_{\mathcal T_{[X]}}$ is $\mathfrak S$-self-adjoint with respect to this restriction.

The function
\begin{equation}\label{ad1}
F(z)\ddd \left[(\jo_{[X]}-z)^{-1}\chi_{\mathcal T_{[X]}}\Psi,\delta^{(X)}\right] = \frac{(\mathfrak S\Psi)_X}{E-z}, \quad (\mathfrak S\Psi)_X \neq 0,
\end{equation}
is well-defined in small punctured neighborhood of $E$ because operator $\jo_{[X]}-z$ is invertible there. Since $\jo_{[X]}$ is $\mathfrak S$-self-adjoint, we can write
\begin{equation}
\label{ad2}
F(z) = \left[\chi_{\mathcal T_{[X]}}\Psi,(\jo_{[X]}-\bar{z})^{-1}\delta^{(X)}\right] = -\frac{[\chi_{\mathcal T_{[X]}}\Psi,m_X\chi_{\mathcal T}p(\bar z)]}{P_{\Pi(X_{(p)})}(z)},
\end{equation}
where we also used \eqref{gevn2} and the fact polynomials \( P_\n(x) \) have real coefficients.  Since $E$ is a pole of $F(z)$ by \eqref{ad1}, the denominator in the right hand side of \eqref{ad2} vanishes at $E$ and we have \( E \in E_{\Pi(X_{(p)})}\) as claimed.

To prove the second statement of the lemma, we only need to show that \( X \) has a sibling, see \eqref{jnt}. That is true since otherwise
\[
0 = E\Psi_{X_{(p)}} = (\jackn \Psi)_{X_{(p)}} = V_{X_{(p)}}\Psi_{X_{(p)}} + W_{X_{(p)}}^{1/2}\Psi_{X_{(g)}} +  W_X^{1/2}\Psi_X = W_X^{1/2}\Psi_X
\]
by \eqref{jackn}, which is clearly impossible as \( W_X>0 \) and $\Psi_X\neq 0$. 

Consider the last claim. Let $E$ be an eigenvalue and $\Psi$ be its eigenfunction. If $\Psi_O\neq 0$, we have $E\in E_{\Pi(O_{(p)})}\subseteq \mathcal E_{\vec \kappa,\vec N}$ by the definition. If $\Psi_O=0$, let \( Z \) be a vertex with the shortest path to \( O \) among all vertices \( X \) for which \( \Psi_X\neq0 \) and \( \Psi_Y=0 \) for all \( Y\in\mathrm{path}(X,O) \), \(Y\neq X \). Since $Z\neq O$, \( Z_{(p)} \in\jnt \) by the second claim and therefore \( E\in E_{\Pi(Z_{(p)})}\subseteq \mathcal E_{\vec\kappa,\vec N} \).
\end{proof}

\begin{Rem}
 Notice that assumption \eqref{as-4} was not used in the proof.
\end{Rem}

\begin{Lem}
\label{lem:2}
Let $E\in \mathcal E_{\vec \kappa,\vec N}$ and \( X\in\jnts \). Then, \( E\in\sigma(\jackn) \) and \( b(E,X) \) is a corresponding eigenvector.
\end{Lem}
\begin{proof}
Let \( E \) be a zero \( P_{\Pi(O_{(p)})}(x) \).  In this case \eqref{gevn2} states that \( \jackn p(E) =Ep(E) \) and therefore \( E \) is indeed an eigenvalue with an eigenvector \( b(E,O_{(p)}) \). Now, let \( E \in \mathcal E_{\vec \kappa,\vec N} \) and \( X \in \jnt \). We need to show that \( b(E,X) \) is an eigenvector with eigenvalue \( E \). Recall that $\mathcal T_{\vec N [X]}$ denotes the subtree of $\mathcal T_{\vec N}$ which has $X$ as its root and observe that
\[
\big(\jackn b(E,X) \big)_Y =0 = Eb_Y(E,X), \quad Y\not\in \mathcal T_{\vec N [X]},
\]
by the definition of \( b(E,X) \), see \eqref{jackn}. Moreover, let \( \upsilon_i \ddd (-1)^{i+\sigma_{X_{(ch),i}}} W_{X_{(ch),i}}^{-1/2}p_{X_{(ch),i}}^{-1}(E) \). Then
\[
\big(\jackn b(E,X) \big)_X = \sum_{i=1}^2(-1)^{\sigma_{X_{(ch),i}}} W_{X_{(ch),i}}^{1/2}b_{X_{(ch),i}}(E,X) = 0 = Eb_X(E,X)
\]
by  \eqref{jackn} and the choice of \( \upsilon_i \). Furthermore,
\begin{multline*}
\big(\jackn b(E,X) \big)_{X_{(ch),l}} = \big(\jackn \upsilon_l p(E) \big)_{X_{(ch),l}} - (-1)^{\sigma_{X_{(ch),l}}} W_{X_{(ch),l}}^{1/2}  \upsilon_l p_X(E) \\ = E \upsilon_l p_{X_{(ch),l}}(E)  = E b_{X_{(ch),l}}(E,X)
\end{multline*}
by \eqref{gevn2}, definition of \( b(E,X) \), and since \( p_X(E) =0 \). Similarly,
\[
\big(\jackn b(E,X) \big)_Y =  \big(\jackn \upsilon_lp(E) \big)_Y = E\upsilon_lp_Y(E) = Eb_Y(E,X), \quad Y\in\mathcal T_{\vec N[X_{(ch),l}]},
\]
which finishes the proof of the lemma.
\end{proof}

\begin{Lem}
\label{lem:3}
Given \( E \in \mathcal E_{\vec \kappa,\vec N}\), the vectors in the system \( \big\{b(E,X):~X\in\jnts\big\} \) are linearly independent.
\end{Lem}
\begin{proof}
Assume that $E\in E_{\Pi(O_{(p)})}$, the proof for other cases is similar. Let \( \beta(Z) \), \( Z\in\jnts \), be constants such that
\[
\beta(O_{(p)}) b_Y(E,O_{(p)})+\sum_{Z\in\jnt} \beta(Z) b_Y(E,Z)=0
\]
is true for all $Y$. Due to assumption \eqref{as-4} with \( Y=O \) and the very construction of \( b(E,Z) \), it holds that
\[
b_O(E,O_{(p)}) = p_O(E)\neq 0 \quad \text{and} \quad b_O(E,Z) =0,~~Z\in\jnt.
\]
Thus, it must hold that $\beta(O_{(p)})=0$. Next, let \( X\in\jnt \) be any vertex such that the path from $X$ to $O$ contains no other elements in $\jnt$. That and assumption \eqref{as-4} then yield that 
\[
b_{X_{(ch),1}}(E,X) = p_{X_{(ch),1}}(E)\neq 0 \quad \text{and} \quad b_{X_{(ch),1}}(E,Z) = 0, ~~ Z\in\jnt\setminus\{X\}.
\]
Hence, \( \beta(X)=0 \). Going down the tree \( \mathcal T_{\vec N} \) in this fashion, we can inductively show that \( \beta(Z)=0 \) for every \( Z\in\jnts \), thus, proving linear independence. 
\end{proof}

\begin{Lem}
\label{lem:4}
Suppose $\Psi$ is an eigenvector of \( \jackn \) with eigenvalue $E$. If $\Psi_O=0$, then $\Psi_Y=0$ for all $Y\in \mathcal W_1(E)$, where the waves \( \mathcal W_k(E) \) were defined in Section~\ref{ss:123}. 
\end{Lem}
\begin{proof}If $O\in \jnt$, then \( \mathcal W_1(E)=\{O\} \) by definition and the claim is obvious. Otherwise, take \( O_{(ch),l} \in \mathcal W_1(E) \). If \( \Psi_{O_{(ch),l}}\neq 0 \) were true, then it would hold that \( O\in \jnt \) by Lemma~\ref{lem:1} which is a contradiction. Furthermore, if the desired claim were false at another vertex of \(\mathcal W_1(E) \), there would exist \( X\in\mathcal W_1(E) \) such that \( \Psi_X \neq 0 \) and \( \Psi_{X_{(p)}} = \Psi_{X_{(g)}} =0 \), where \( X_{(g)} \) is the parent of \( X_{(p)} \). Then, \( X_{(p)} \in \jnt \) by Lemma~\ref{lem:1}, which contradicts the very definition of \( \mathcal W_1(E) \).
\end{proof}

\begin{Lem}
\label{lem:5}
Given \( E \in \mathcal E_{\vec \kappa,\vec N}\), the system \( \big\{b(E,X):~X\in\jnts\big\} \) spans the subspace of eigenvectors corresponding to \( E \).
\end{Lem}
\begin{proof}
Let $\Psi$  be an eigenvector that corresponds to $E$. First, consider the values of \( \Psi \) on \( \mathcal W_1(E) \). If \( \Psi_O=0 \), then $\Psi_Y=0$ for all $Y\in \mathcal{W}_1(E)$ by Lemma \ref{lem:4} and we set \( \Psi^{(1)} \ddd \Psi \). Otherwise, $\Psi_O\neq 0$ and $E$ is a zero of $P_{\Pi(O_{(p)})}$ according to Lemma~\ref{lem:1}. In particular, \( P_{\vec N}(E) \neq 0 \) due to assumption \eqref{as-4} with \( Y=O \) and so \( p_O(E)\neq 0 \). Then, we set
 \[
\Psi^{(1)} \ddd \Psi - \big(\Psi_O/p_O(E)\big) b(E,O_{(p)}).
 \]
Since \( P_{\Pi(O_{(p)})}(E) = 0 \), it follows from \eqref{gevn2} and the definition of \( b(E,O_{(p)}) \) that $\Psi^{(1)}$ is also an eigenvector corresponding to $E$. Since \( \Psi^{(1)}_O = 0 \), we have \( \Psi^{(1)}_Y = 0 \) for every \( Y\in\mathcal W_1(E) \) by Lemma~\ref{lem:4} as desired.

Second, we consider the values of \( \Psi^{(1)} \) on \( \mathcal W_2(E)\cup\mathcal W_1(E) \). Fix \( X \in\mathcal F_1(E)\setminus\mathcal C \). By the very definition of the first front we have that \( X\in\jnt \). Choose \( \beta(X) \) so that
\[
\Phi_{X_{(ch),1}}=0, \quad \Phi \dd \Psi^{(1)} - \beta(X) b(E,X).
\]
Since \( \Phi \) is an eigenvector corresponding to \( E \) that vanishes at \( X_{(ch),1} \), \( X \), and \( X_{(p)} \), it follows from \eqref{jackn} that
\begin{multline*}
0 = E\Phi_X = (\jackn \Phi)_X = V_X \Phi_X + W_X^{1/2}\Phi_{X_{(p)}} + (-1)^{\sigma_{X_{(ch),1}}} W_{X_{(ch),1}}^{1/2}\Phi_{X_{(ch),1}} \\ + (-1)^{\sigma_{X_{(ch),2}}} W_{X_{(ch),2}}^{1/2}\Phi_{X_{(ch),2}} = W_{X_{(ch),2}}^{1/2}\Phi_{X_{(ch),2}}.
\end{multline*}
Thus, \( \Phi \) vanishes at \( X_{(ch),2} \) as well. Now, as in the proof of Lemma \ref{lem:4}, we apply the second claim of Lemma~\ref{lem:1} to conclude that \( \Phi \) vanishes at all \( Y\in\mathcal T_{\vec N[X]}\cap \mathcal W_2(E) \). Therefore, we can set
\[
\Psi^{(2)} \ddd \Psi^{(1)} - \sum_{X\in \mathcal{F}_1(E)}\beta(X) b_Y(E,X),
\]
which is an eigenvector corresponding to \( E \) that vanishes at all \( Y\in \mathcal W_2(E)\cup\mathcal W_1(E) \).
Continuing in the same way, we decompose $\Psi$ into the sum of canonical eigenvectors. 
\end{proof}

\begin{Lem}
\label{lem:6}
It holds that
\[
\# \,\mathcal{V}_{\vec{N}}=\sum_{E\in \mathcal{E}_{\vec\kappa,\vec N}}\#\, \jnts\,.
\]
\end{Lem}
\begin{proof}
Recall that according to our assumption \eqref{as-4} all zeroes of any polynomial $P_\n(x)$ are simple and there are exactly \( |\n| \) of them since \( \vec\mu \) is perfect. Given an eigenvalue \( E \), each polynomial \( P_\n(x) \), \( \n\in\mathbb N^2 \), such that \( P_\n(E)=0 \), generates as many canonical eigenvectors as the number of vertices \( X \) for which \( \Pi(X) = \n \) (the number of paths from \( \n \) to \(\vec N \) in \( \mathcal R_{\vec N} \)). Hence, the number of the canonical eigenvectors that each polynomial \( P_\n(x) \), \( \n\in\mathbb N^2 \), generates is equal to \( |\n|\cdot \#\Pi^{-1}(\n) \). Therefore, the total number of eigenvectors is equal to
\[
\sum_{E\in \mathcal{E}_{\vec\kappa,\vec N}}\#\, \jnts\, = |\vec N| + 1 + \sum_{\n\in\mathcal R_{\vec N}\cap \mathbb N^2}|\n|\binom{|\vec N|-|\n|}{N_1-n_1},
\]
where \( |\vec N|+1 \) is the number of the trivial canonical eigenvectors (\eqref{as-4} is used here too as well as equality \( \kappa_1+\kappa_2=1\)). The above formula is true for every Jacobi matrix on \( \mathcal T_{\vec N} \), including the self-adjoint ones (that do exist). For the self-adjoint matrices the desired claim is a standard fact of linear algebra (the number of linearly independent eigenvectors of a self-adjoint matrix is equal to the dimension of the space). Hence, it holds for all Jacobi matrices.
\end{proof}

\begin{Rem}
There is an alternative proof of this lemma using an inductive argument.
\end{Rem}

\begin{proof}[Proof of Theorem~\ref{sde_1}]
The first claim follows from Lemmas~\ref{lem:1} and~\ref{lem:2}. The validity of the second one is due to Lemmas~\ref{lem:3} and~\ref{lem:5}. The formula for $g_E$ is a trivial consequence of the second claim. Since all eigenspaces of a linear operator are mutually linearly independent, the last claim follows from Lemma~\ref{lem:6}.
\end{proof}

For reader's convenience, we include the proof of the following standard result.
\begin{Lem}
\label{lem:7}
Suppose $\{\psi_1,\ldots,\psi_n\}$ is a $\mathfrak{S}$-orthogonal basis of \( \ell^2(\mathcal V_{\vec N}) \). Then
\[ 
\#\big\{j:~\psi_j~\text{is $\mathfrak{S}$-negative}\big\}=i_- \quad \text{and} \quad \#\big\{j:~\psi_j~\text{is $\mathfrak{S}$-positive}\big\}=i_+\,.
\]
\end{Lem}
\begin{proof}
Notice first that none of $\{\psi_j\}$ is $\mathfrak{S}$-neutral since otherwise, we would have $[\psi_k,f]=0$ for all \( f\in\ell^2(\mathcal V_{\vec N}) \) and some $k$. In particular, this would yield that
\[
0 = \big|\big[\psi_k,\delta^{(Y)}\big]\big| = |\psi_{k|Y}|\,\big|\big[\delta^{(Y)},\delta^{(Y)}\big]\big| = |\psi_{k|Y}|
\]
for every \( Y\in\mathcal V_{\vec N} \), which is clearly impossible as \( \psi_k\not\equiv0 \). Thus, we can assume that  $[\psi_j,\psi_j]=\pm 1$ for all $j$. Let $k_-$ and \( k_+ \) be the numbers of $\mathfrak{S}$-negative  and $\mathfrak{S}$-positive vectors in $\{\psi_j\}$, respectively. Assume without loss of generality that $\{\psi_1,\ldots, \psi_{k_+}\}$ are $\mathfrak{S}$-positive. Since $\{\psi_j\}$ is a basis, we can write
\[
f =\sum_Y  f_Y \delta^{(Y)} = \sum_j x_j \psi_j,  \quad f\in\ell^2(\mathcal V_{\vec N}),
\]
for some constants \( \{x_j\} \). Let \( \mathcal V_+ \) and \( \mathcal V_- \) be the subsets of \( \mathcal V_{\vec N} \) for which \( \delta^{(Y)} \) is \( \mathfrak S \)-positive and \( \mathfrak S \)-negative, respectively. Clearly, \( \# \mathcal V_\pm=i_\pm \) by definition. Then
\[
\sum_{Y\in \mathcal V_+}|f_Y|^2 - \sum_{Y\in \mathcal V_-} |f_Y|^2 = \langle \mathfrak{S}f,f\rangle=[f,f]=\sum_{j=1}^{k_+}|x_j|^2 - \sum_{j=k_++1}^{n} |x_j|^2 \,.
\]
The desired claim now follows from Sylvester's law of inertia for Hermitian matrices,  \cite[Theorem~X.18]{Gantmacher} (the numbers of positive and negative squares do not depend on the choice of a representation of a Hermitian form).
\end{proof}

\section{Appendix to Part~\ref{part1}}
\label{s:ap1}

In the end of Subsection~\ref{sub_1},  we have listed a number of systems of MOPs whose recurrence coefficients do not satisfy condition \eqref{as-2}. Most of them come from special orthogonality measures and their  recurrence coefficients are known explicitly.  The only exception in that list are Nikishin systems. A vector \( \vec\mu = (\mu_1,\mu_2) \) defines a \emph{Nikishin system} if there exists a measure \( \tau \) such that
\begin{equation}
\label{nik-sys}
d\mu_2(x) = \widehat\tau(x)d\mu_1(x) \quad \text{and} \quad \Delta_1\cap\Delta_\tau = \varnothing,
\end{equation}
where $\widehat\tau(z)$ is the Markov function of \(\tau \), see \eqref{markov}, \( \Delta_1\ddd\mathrm{ch}(\supp\,\mu_1) \), and \( \Delta_\tau\ddd\mathrm{ch}(\supp\,\tau) \) (here, \( \mathrm{ch}(\cdot) \) stands for the convex hull).  
Given two sets $E_1$ and $E_2$, we write $E_1<E_2$ if \( \sup E_1 <\inf  E_2 \). In what follows,  we assume that 
\begin{equation}
\label{dtd1}
\Delta_\tau<\Delta_1.
\end{equation}
The case when $\Delta_\tau>\Delta_1$ can be handled similarly.

It is known that Nikishin systems are perfect \cite{bl,ds,g_b}.
The goal of this appendix is to show that the recurrence coefficients \( \{ a_{\n,1},a_{\n,2}\}_{\n\in\mathbb N^2} \), see \eqref{1.5}--\eqref{1.7}, of Nikishin systems have a definite sign pattern. That explains how the indefinite inner product $\mathfrak{S}$ should be defined to make the associated Jacobi matrix $\mathfrak{S}$-self-adjoint. Recall \eqref{coef-a}.

\begin{Thm}
\label{thm:sign-N}
For all \( \n\in\mathbb N^2 \) and \( j\in\{1,2\} \) it holds that
\[
\sgn\, a_{\n,j} = (-1)^{j-1}, ~~ n_2\leq n_1, \quad \text{and} \quad \sgn\, a_{\n,j} = (-1)^j, ~~ n_2\geq n_1+1.
\]
\end{Thm}

To prove this theorem, let us make the following observation. It holds that
\begin{equation}\label{1sd1}
\frac1{\widehat\tau(z)} - \frac z{m_0(\tau)} + \frac{m_1(\tau)}{m_0^2(\tau)} = \mathcal O\left(\frac1z\right)
\end{equation}
as \( z\to\infty \), where \( m_l(\tau) \dd \int x^ld\tau(x) \). Next, we will use some basic facts from the theory of Herglotz-Nevalinna functions, see Section~\ref{s:Pi} further below. As the left-hand side of \eqref{1sd1} has positive imaginary part in \( \C_+ \) and is holomorphic and vanishing at infinity, there exists a positive measure \( \tau_d \) supported on \( \Delta_\tau \), which we call the \emph{dual measure} of \( \tau \), such that
\begin{equation}
\label{dual-measure}
\frac1{\widehat\tau(z)} - \frac z{m_0(\tau)} + \frac{m_1(\tau)}{m_0^2(\tau)} = -\widehat\tau_d(z).
\end{equation}
The bulk of the proof of Theorem~\ref{thm:sign-N} is contained in Lemmas~\ref{lem:n1} and~\ref{lem:n2}. These lemmas and ideas behind their proofs are not new, see, for example, \cite{grs,bl,ds}, but we decided to include them as their proofs are  short, they are formulated exactly in the way we need, and their inclusion makes the paper as self-contained as possible.

Let \(\{ P_\n(x) \}\) be monic type II MOPs for  Nikishin system \eqref{nik-sys}--\eqref{dtd1}. Define
\begin{equation}
\label{hnj}
h_{\n,j} \ddd \int P_\n^2(x)d\mu_j(x).
\end{equation}
Recall the functions of the second kind \( R_{\n,j}(z) \) defined in \eqref{Ln}. It follows from orthogonality relations \eqref{1.2} that
\begin{equation}
\label{n1}
R_{\n,j}(z)=\frac{1}{p(z)}\int_{\Delta_1}\frac{p(x)P_{\n}(x)}{z-x}d\mu_j(x)
\end{equation}
for any polynomial $p(z)$ such that $\deg p\le n_j$. Moreover, the Taylor expansion of $(z-x)^{-1}$ at infinity gives
\begin{equation}
\label{n9}
R_{\n,j}(z) = \frac{h_{\n,j}}{z^{n_j+1}}\Bigl(1+\mathcal O(z^{-1})\Bigr) \quad \text{as} \quad z\to \infty\,.
\end{equation}
Then, the following lemma takes place.

\begin{Lem}
\label{lem:n1}
Let functions \( R_{\n,j}(z) \) be given by \eqref{Ln} for a Nikishin system \eqref{nik-sys}--\eqref{dtd1} and \( \tau_d \) be the dual measure of \( \tau \). The functions \( R_{\n,j}(z) \) satisfy
\[
\int x^kR_{\n,1}(x)d\tau(x) = 0 \quad \text{and} \quad \int x^kR_{\n,2}(x)d\tau_d(x) = 0
\]
for \( k\leq \min\{n_1,n_2-1\} \) and \( k\leq \min\{n_1-1,n_2-2\} \), respectively. It further holds that
\[
\int x^{n_2}R_{\n,1}(x)d\tau(x) = -h_{\n,2} \quad \text{and} \quad \int x^{n_1}R_{\n,2}(x)d\tau_d(x) = h_{\n,1}
\]
when \( n_2\leq n_1 \) and \( n_2\geq n_1+2\), respectively. Finally, it holds that
\[
\|\tau\|h_{\n,1} - h_{\n,2} = \int x^{n_2}R_{\n,1}(x)d\tau(x) =  \|\tau\| \int x^{n_1}R_{\n,2}(x)d\tau_d(x)
\]
when \( n_2 = n_1+1 \), where \( \|\tau\|=m_0(\tau) \) is the total mass of \( \tau \).
\end{Lem}
\begin{proof}
We only consider the case \( j=2 \), the argument for $j=1$ is similar. Assume that \( k\leq n_1-1 \). Then
\[
0 = \int P_\n(x)x^kd\mu_1(x) = \int P_\n(x)x^k\widehat\tau^{-1}(x)d\mu_2(x).
\]
If we further assume that \( k\leq n_2-2 \), then we get from \eqref{dual-measure} and orthogonality conditions that
\[
0 = -\int P_\n(x)x^k\widehat\tau_d(x)d\mu_2(x).
\]
Thus, we can deduce from Fubini-Tonelli Theorem that
\[
0 = -\int \left( \int \frac{x^kP_\n(x)}{x-y}d\mu_2(x)\right) d\tau_d(y) = \int y^kR_{\n,2}(y)d\tau_d(y)
\]
as claimed, where we used \eqref{n1} with \( p(x)=x^k \). Similarly, we have that
\[
h_{\n,1} =  \int P_\n(x)x^{n_1}d\mu_1(x) =  -\int P_\n(x)x^{n_1}\widehat\tau_d(x)d\mu_2(x) = \int y^{n_1}R_{\n,2}(y)d\tau_d(y)
\]
when \( n_1\leq n_2-2 \). Furthermore, if \( n_1 = n_2-1 \), we get from \eqref{dual-measure} that
\[
h_{\n,1} = -\int P_\n(x)x^{n_1}\widehat\tau_d(x)d\mu_2(x) + \|\tau\|^{-1}\int P_\n(x)x^{n_2}d\mu_2(x) =  \int y^{n_1}R_{\n,2}(y)d\tau_d(y) + \|\tau\|^{-1}h_{\n,2}. \qedhere
\]
\end{proof}

Let \( r_{\n,j}(x) \) be a monic polynomial with zeroes on \( \Delta_\tau \) such that \( R_{\n,j}(x)/r_{\n,j}(x) \) is analytic and non-vanishing on \( \Delta_\tau \). It follows from the previous lemma that \( r_{\n,1}(x) \) has at least  \( \min\{n_1,n_2-1\} + 1 \) different zeroes while \( r_{\n,2}(x) \) has at least  \( \min\{n_1-1,n_2-2\} + 1 \) different zeroes.

\begin{Lem}
\label{lem:n2}
If \( n_2\leq n_1+1 \), \( r_{\n,1}(x) \) has degree exactly \( n_2 \) (in particular, all its zeroes are simple) and \( R_{\n,1}(z)/r_{\n,1}(z) \) is non-vanishing in \( \C\setminus\Delta_1 \). Moreover,
\[
\int x^kP_\n(x)\frac{d\mu_1(x)}{r_{\n,1}(x)} = 0 \quad \text{and} \quad \int P_\n^2(x)\frac{d\mu_1(x)}{r_{\n,1}(x)} = h_{\n,1},
\]
where the first relation holds for any \( k<|\n| \).

Similarly, if \( n_2 \geq n_1+1 \), \( r_{\n,2}(x) \) has degree exactly \( n_1 \) (in particular, all its zeroes are simple) and \( R_{\n,2}(z)/r_{\n,2}(z) \) is non-vanishing in \( \C\setminus\Delta_1 \). Furthermore,
\[
\int x^kP_\n(x)\frac{d\mu_2(x)}{r_{\n,2}(x)} = 0 \quad \text{and} \quad \int P_\n^2(x)\frac{d\mu_2(x)}{r_{\n,2}(x)} = h_{\n,2},
\]
where again the first relation holds for any \( k<|\n| \).
\end{Lem}
\begin{proof}
It follows from the remark before the lemma that \( \deg r_{\n,j} = n_{3-j} + m_j\), \( m_j\geq0 \), in the considered cases. Therefore, it follows from \eqref{n9} that
\[
R_{\n,j}(z)/r_{\n,j}(z) = h_{\n,j} z^{-|\n|-m_j-1} + \mathcal O\left( z^{-|\n|-m_j-2} \right)
\]
as \(z\to\infty \) and the ratio is a holomorphic function in \( \overline\C\setminus\Delta_1 \). Let \( \Gamma \) be a smooth Jordan curve that encircles \( \Delta_1 \) but not \( \Delta_\tau \). Then, by integrating over $\Gamma$ in positive direction we get
\[
0 = \frac1{2\pi\ic}\int_\Gamma s^kR_{\n,j}(s)\frac{ds}{r_{\n,j}(s)} = \int P_\n(x)\left(\frac1{2\pi\ic}\int_\Gamma \frac{s^k}{s-x}\frac{ds}{r_{\n,j}(s)}\right) d\mu_j(x) =  \int x^k P_\n(x)\frac{d\mu_j(x)}{r_{\n,j}(x)} 
\]
for \( k<|\n|+m_j \), by Cauchy theorem, Fubini-Tonelli theorem, and Cauchy integral formula. Since \( d\mu_j(x)/r_{\n,j}(x) \) is a measure of constant sign on \( \Delta_1 \), \( P_\n(x) \) cannot be orthogonal to itself. Thus, \( m_j=0 \). Now, if there existed another real zero \( x_0\not\in\Delta_1\cup\Delta_\tau \) of \( R_{\n,j}(z) \), then the above argument can be applied with \( r_{n,j}(z) \) replaced by \( (z-x_0)r_{\n,j}(z) \) and \( \Gamma \) not containing \( x_0 \) in its interior to arrive at a contradiction, namely, that \( P_\n(x) \) is orthogonal to itself with respect to a measure of constant sign. If \( R_{\n,j}(z_0)=0 \) for some \( z_0\not\in \R \), then \( R_{\n,j}(\bar z_0)=0 \) by conjugate-symmetry, and therefore the above argument can be used with \( (z-z_0)(z-\bar z_0)r_{\n,j}(z) \). Using \eqref{n9} one more time, we get that
\[
h_{\n,j} = \frac1{2\pi\ic}\int_\Gamma s^{|\n|}R_{\n,j}(s)\frac{ds}{r_{\n,j}(s)} = \int x^{|\n|} P_\n(x)\frac{d\mu_j(x)}{r_{\n,j}(x)} = \int P_\n^2(x)\frac{d\mu_j(x)}{r_{\n,j}(x)} 
\]
by orthogonality and since \( P_\n(x) \) is monic.
\end{proof}

\begin{Cor}
\label{cor:n1}
It holds that
\[
\sgn\,h_{\n,1} = 1\quad \text{and} \quad \sgn\,h_{\n,2} = 1 
\]
when \( n_2\leq n_1+1 \) and \( n_2\geq n_1+1\), respectively.
\end{Cor}
\begin{proof}
The claim follows from Lemma~\ref{lem:n2} since  \( \Delta_\tau<\Delta_1 \) while each \( r_{\n,j}(z) \) is a monic polynomial.
\end{proof}

\begin{Cor}
\label{cor:n2}
It holds that
\[
R_{\n,1}(z) = \frac{r_{\n,1}(z)}{P_\n(z)}\int\frac{P_\n^2(x)}{z-x}\frac{d\mu_1(x)}{r_{\n,1}(x)} \quad \text{and} \quad R_{\n,2}(z) = \frac{r_{\n,2}(z)}{P_\n(z)}\int\frac{P_\n^2(x)}{z-x}\frac{d\mu_2(x)}{r_{\n,2}(x)}
\]
when \( n_2\leq n_1+1 \) and \( n_2\geq n_1+1\), respectively.
\end{Cor}
\begin{proof}
We have that
\[
R_{\n,j}(z) = \int \frac{P_\n(x)r_{\n,j}(x)}{z-x}\frac{d\mu_j(x)}{r_{\n,j}(x)}  =  \int \frac{r_{\n,j}(x)-r_{\n,j}(z)}{z-x}P_\n(x) \frac{d\mu_j(x)}{r_{\n,j}(x)} + r_{\n,j}(z)\int\frac{P_\n(x)}{z-x}\frac{d\mu_j(x)}{r_{\n,j}(x)}.
\]
Since \( n_{3-j}-1<|\n| \) is the degree of \( (r_{\n,j}(\cdot)-r_{\n,j}(z))/(z-\cdot)\), it holds that
\[
R_{\n,j}(z) = r_{\n,j}(z)\int\frac{P_\n(x)}{z-x}\frac{d\mu_j(x)}{r_{\n,j}(x)}.
\]
Using the same argument one more time yields the desired claim.
\end{proof}

\begin{Cor}
\label{cor:n3}
It holds that
\[
\sgn\,h_{\n,1} = (-1)^{|\n|+1} \quad \text{and} \quad \sgn\,h_{\n,2} = (-1)^{|\n|} 
\]
when \( n_2 \geq n_1+2 \) and \( n_2 \leq n_1 \), respectively.
\end{Cor}
\begin{proof}
It follows from the previous corollary that 
\begin{equation}
\label{Rrj}
\sgn\left(R_{\n,j}(x)/r_{\n,j}(x) \right) = (-1)^{|\n|+1}, \quad x\in\Delta_\tau,
\end{equation}
when \( n_2 \geq n_1+1 \) for \( j=2 \) and \( n_2 \leq n_1+1 \) for \( j=1 \). The claim now follows from Lemma~\ref{lem:n1} since
\[
h_{\n,1} = \int x^{n_1}R_{\n,2}(x)d\tau_d(x) = \int r_{\n,2}(x)R_{\n,2}(x)d\tau_d(x)
\]
when \( n_2\geq n_1+2 \) and
\[
h_{\n,2} = -\int x^{n_2}R_{\n,1}(x)d\tau(x) = -\int r_{\n,1}(x)R_{\n,1}(x)d\tau(x)
\]
when \( n_2\leq n_1 \).
\end{proof}

\begin{proof}[Proof of Theorem~\ref{thm:sign-N}]
As is well known \cite[Theorem~23.1.11]{ismail}, if we multiply equation \eqref{1.7} by \( x^{n_j-1} \) and integrate agains the measure \( \mu_j \), we will get

\begin{equation}
\label{an-hn}
a_{\n,j} = \frac{\int P_\n(x)x^{n_j}d\mu_j(x)}{\int P_{\n-\vec e_j}(x)x^{n_j-1}d\mu_j(x)} = \frac{h_{\n,j}}{h_{\n-\vec{e}_j,j}},
\end{equation}
where we used \eqref{hnj} and orthogonality relations \eqref{1.2} to get the second equality. The claim of the theorem now follows from Corollaries~\ref{cor:n1} and~\ref{cor:n3}.
\end{proof}

\part{Jacobi matrices on infinite rooted Cayley trees}
\label{part2}

Below we introduce a notion of a Jacobi matrix on an infinite 2-homogenous rooted tree whose coefficients generated by MOPs. 

\section{Definitions}

Let \( \vec \mu \) be a perfect system of measures on the real line with recurrence coefficients \( \{ a_{\n,i},b_{\n,i} \} \), see \eqref{1.5} and \eqref{1.7}. Assume that
\begin{equation}
\label{as-5}
\sup \limits_{\n\in \Z_+^2,~i\in \{1,2\}} |a_{\n,i}|<\infty \quad \text{and} \quad \sup \limits_{\n \in \Z_+^2,~i\in \{1,2\}}|b_{\n,i}|<\infty  \,. 
\end{equation}
Conditions \eqref{as-5} used along the marginal directions imply that the classical Jacobi matrices corresponding to $\mu_1$ and $\mu_2$ have bounded coefficients and therefore $\mu_1,\mu_2\in \mathfrak{M}$.

\subsection{Rooted Cayley tree}
\label{ss:rct}
\begin{figure}[h!]
\includegraphics[scale=.6]{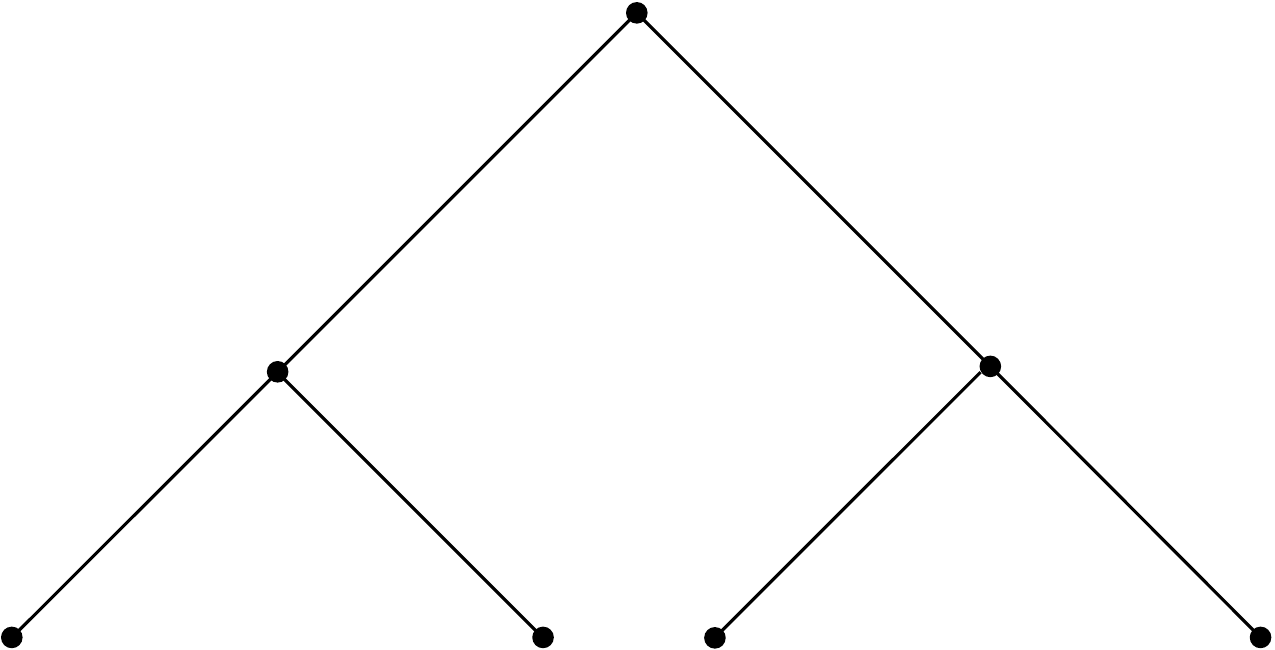}
\begin{picture}(0,0)
\put(-110,110){$(1,1)\sim O=Y_{(p)}$}
\put(-170,47){$(2,1)$}
\put(-48,47){$(1,2)\sim Y=O_{(ch),2}$}
\put(-215,0){$(3,1)$}
\put(-128,0){$(2,2)$}
\put(-95,0){$(2,2)\sim Y_{(ch),1}$}
\put(-2,0){$(1,3)\sim Y_{(ch),2}$}
\end{picture}
\caption{Three generations of \( \mathcal T \).}
\label{fig:inf-tree}
\end{figure}
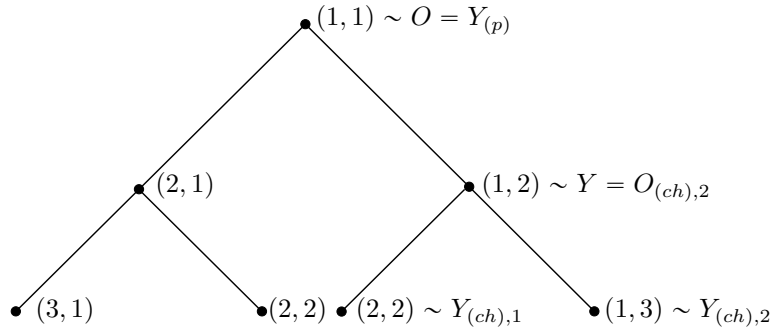

Hereafter, we let \( \mathcal T \) stand for an infinite $2$-homogeneous rooted tree (rooted Cayley tree) and \( \mathcal V \) for the set of its vertices with \( O \) being the root. On the lattice \( \mathbb{N}^2 \), consider an infinite path
\[
\big\{\n^{(1)},\n^{(2)}, \ldots\big\}, \quad \n^{(1)}=\vec 1\ddd(1,1)\quad  \text{and} \quad \n^{(l+1)}=\n^{(l)}+\vec{e}_{k_l}, ~~k_l\in \{1,2\},~~l\in\mathbb N.
\]
Clearly, these are paths for which, as we move from $\vec 1$ to infinity, the multi-index of each next vertex is increasing by $1$ at exactly one position.  Each such path can be mapped bijectively to a non-self-intersecting path on \( \mathcal T \) that starts at $O$, see Figure~\ref{fig:inf-tree}. This construction defines a projection $\Pi:\mathcal{V}\to\mathbb{N}^2$ as follows: given $Y\in \mathcal{V}$ we consider a path from $O$ to $Y$, map it to a path on \( \mathbb{N}^2 \) and let $\Pi(Y)$ be the endpoint of the mapped path. Every vertex $Y\in \mathcal{V}$, which is different from $O$, has a unique parent, which we denote by \( Y_{(p)} \). That allows us to define the following index function:
\begin{equation}
\label{imath}
\imath:\mathcal V \to \{1,2\}, \quad Y\mapsto \imath_Y ~\text{ such that }~ \Pi(Y) = \Pi(Y_{(p)}) + \vec e_{\imath_Y}.
\end{equation}
This way, if \( Z=Y_{(p)} \), then we write that \( Y=Z_{(ch),\iota_Y} \), see Figure~\ref{fig:inf-tree}. Recall that for a function \( f \) on \( \mathcal V \), we  denote its value at a vertex \( Y\in\mathcal V \)  by \( f_Y \). As before, we introduce an artificial vertex \( O_{(p)} \), a formal parent of the root \( O \). We do not include \( O_{(p)} \) into \( \mathcal V \), but we do extend every function \( f \) on \( \mathcal V \) to \( O_{(p)} \) by setting \( f_{O_{(p)}} = 0 \). We denote  the space of square-summable functions on \( \mathcal V \) by \( \ell^2(\mathcal V) \) and  the standard inner product generating \( \ell^2(\mathcal V) \) by \( \langle\cdot,\cdot\rangle \).

The above construction differs from the one in Section~\ref{s:ft-def} in the following ways: the projection \( \Pi \) maps onto the lattice \( \mathbb N^2 \), not \( \Z_+^2 \); the values \( |\Pi(Y)| \) increase rather than decrease as we go down the tree;  the index function \( \iota_Y \) now tells which coordinate of \( \Pi(Y_{(p)}) \) needs to increase rather than decrease to get \( \Pi(Y) \).

\subsection{Jacobi matrices}

In this subsection we specialize definition \eqref{sad1} to the case of Jacobi matrices on \( \mathcal T \) whose potentials \( V,W \) come from \( \vec\mu \). As in the previous part, we fix \( \vec\kappa\in\R^2 \) such that \( |\vec\kappa| = 1 \). We define the potentials \( V=V^{\vec\mu},W=W^{\vec\mu}:\mathcal V \to \R \) (again, as with most quantities depending on \( \vec \mu \), we drop the dependence on \( \vec\mu \) from notation) by
\begin{equation}
\label{jack1}
V_O \ddd \kappa_1 b_{(0,1),1} + \kappa_2 b_{(1,0),2}, \quad W_O\ddd1, \quad \text{and} \quad V_Y \ddd b_{\Pi(Y_{(p)}),\iota_Y},  \quad W_Y \ddd \big|a_{\Pi(Y_{(p)}),\iota_Y}\big|, ~~Y\neq O.
\end{equation}
Notice the difference in definition of \( V \) as compared to \eqref{potentials}. As before, this definition is consistent with \eqref{sad1} if we let \( W_{Y_{(p)},Y} = W_{Y,Y_{(p)}} = W_Y \). We further choose function \( \sigma:\mathcal V \to \{0,1\} \) to recover the signs of the recurrence coefficients \( a_{\n,i} \) exactly as in \eqref{nu-fun}. With these definitions, \eqref{sad1} specializes to the following Jacobi matrix \( \jack=\jack^{\vec\mu} \) on \( \mathcal T \):
\begin{equation}
\label{jack2}
(\jack f)_Y \ddd  V_Yf_Y + W_Y^{1/2}f_{Y_{(p)}} + (-1)^{\sigma_{Y_{(ch),1}}}W_{Y_{(ch),1}}^{1/2}f_{Y_{(ch),1}} + (-1)^{\sigma_{Y_{(ch),2}}}W_{Y_{(ch),2}}^{1/2}f_{Y_{(ch),2}}.
\end{equation}
Due to its local nature, \( \jack \) is defined on the set of all functions on \( \mathcal V \). Moreover, assumption \eqref{as-5} also shows that it is a bounded operator on \( \ell^2(\mathcal V) \) and therefore we can talk about its spectrum \( \sigma(\jack) \). Notice also that if \( a_{\n,i}>0 \) when \( n_i>0 \), \( \n=(n_1,n_2) \), i.e., if they satisfy \eqref{as-2} (recall also \eqref{coef-a}), the operator \( \jack \) is self-adjoint. Otherwise, let \( \mathfrak S \) be a diagonal matrix on \( \mathcal T \) defined by \eqref{matS} and \( [\cdot,\cdot] \) be the corresponding indefinite inner product on \( \ell^2(\mathcal V) \) given by \eqref{iip}. We define \( \mathfrak S\)-self-adjointess exactly as in \eqref{Ssa}.

\begin{Prop} 
Jacobi matrices $\jack$ and  $\jack^*$ are $\mathfrak{S}$-self-adjoint.
\end{Prop}
\begin{proof}
The operator $\jack$ is bounded in $\ell^2(\cal{V})$ and checking its $\mathfrak{S}$--self-adjointness is identical to the proof of Proposition~\ref{sadj} from the previous part.
\end{proof}

\section{Basic Properties}

Recall the functions \( L_\n(z) \) introduced in \eqref{Ln}. We consider \( z\in\C_+\) as a parameter and define
\begin{equation}
\label{gev1}
l_Y(z) \ddd m_Y^{-1}L_Y(z), \quad L_Y(z) \ddd L_{\Pi(Y)}(z), \quad \text{and} \quad m_Y\ddd   \prod_{Z\in {\rm path}(Y,O)} W_Z^{-1/2},
\end{equation}
where  ${\rm path} (Y,O)$ is the non-self-intersecting path connecting $Y$ and $O$ that includes both $Y$ and \( O \). More generally, we agree that any function \( f=\{f_\n\} \) on the lattice \( \mathbb N^2 \) is also a function on \( \mathcal V \) whose values are \( f_Y \ddd f_{\Pi(Y)} \). 

Recall definition \eqref{markov} of a Markov function. It will be convenient to formally set \( \Pi(O_{(p)}) \ddd \vec\varkappa \ddd (\kappa_2,\kappa_1)\) and
\begin{equation}
\label{LPO}
L_{\vec\varkappa}(z) \ddd \kappa_2 L_{\vec e_1}(z) + \kappa_1 L_{\vec e_2}(z) = \big(\varkappa_1\|\mu_1\|^{-1}\big)\widehat\mu_1(z) + \big(\varkappa_2\|\mu_2\|^{-1}\big)\widehat\mu_2(z),
\end{equation}
where the second equality follows straight from the definition \eqref{sad15} and the normalization \eqref{n_2}. The reason we introduced \( \vec\varkappa \) is that this way the meaning of \( L_{\vec e_i}(z) \) is still the same.

Given \( X\in\mathcal V \), we shall denote by \( \mathcal T_{[X]} \) the subtree of \( \mathcal T \) with the root at \( X \) (in this case \( \mathcal T_{[O]} = \mathcal T \)). We also let \( \mathcal V_{[X]} \)  be the set of vertices of \( \mathcal T_{[X]} \) and denote  the restriction of the inner product in \( \ell^2(\mathcal V) \) to \( \mathcal V_{[X]} \) by the same symbol \( \langle\cdot,\cdot\rangle \).  The notation \( \jox \) and \( l_{[X]} \) stand for the restrictions of \( \jack \) and \( l \) to \( \mathcal T_{[X]} \) and \( \mathcal V_{[X]} \), respectively. In general,  \( f_{[X]} \) will be used to denote the restriction of any function \( f \), defined on \( \mathcal V \) initially, to  the subset \( \mathcal V_{[X]} \).

\begin{Prop}
\label{prop:gev}
It holds that 
\begin{equation}
\label{gev2}
\jack l(z) = zl(z) - L_{\vec\varkappa}(z)\delta^{(O)}.
\end{equation}
Given \( X\in\mathcal V \), \( X\neq O \), we also have
\begin{equation}
\label{gev3}
\jox l_{[X]}(z) = zl_{[X]}(z) - m_X^{-1}L_{X_{(p)}}(z)\delta^{(X)}.
\end{equation}
\end{Prop}
\begin{proof}
By integrating \eqref{1.5} against \( (z-x)^{-1} \) and noticing that \( |\n|\geq2 \) for any \( \n\in\mathbb N^2 \), we get that
\[
zL_\n(z) = L_{\n-\vec e_j}(z) + b_{\n-\vec e_j,j} L_\n(z) + a_{\n,1} L_{\n+\vec e_1}(z) + a_{\n,2} L_{\n+\vec e_2}(z), \quad j\in\{1,2\}.
\]
Fix \( j \) and let \( Y\neq O \), \( \Pi(Y) = \n \), be such that \( \Pi(Y_{(p)}) = \n-\vec e_j \). Then, the above relation can be rewritten as
\[
zL_Y(z) = L_{Y_{(p)}}(z) + V_Y L_Y(z) + (-1)^{\sigma_{Y_{(ch),1}}} W_{Y_{(ch),1}} L_{Y_{(ch),1}}(z) + (-1)^{\sigma_{Y_{(ch),2}}}W_{Y_{(ch),2}} L_{Y_{(ch),2}}(z).
\]
It follows immediately from \eqref{gev1} and the above formula that \( zl_Y(z) = \big(\jack l(z)\big)_Y \), \( Y\neq O \). Similarly, it holds that
\begin{multline*}
zl_O(z) = zL_O(z) =\\ \kappa_1 \left(L_{\vec 1-\vec e_1}(z) + b_{\vec 1-\vec e_1,1}L_O(z) + (-1)^{\sigma_{O_{(ch),1}}} W_{O_{(ch),1}} L_{O_{(ch),1}}(z) + (-1)^{\sigma_{O_{(ch),2}}} W_{O_{(ch),2}} L_{O_{(ch),2}}(z)\right) \\
+  \kappa_2 \left(L_{\vec 1-\vec e_2}(z) + b_{\vec 1-\vec e_2,2}L_O(z) + (-1)^{\sigma_{O_{(ch),1}}} W_{O_{(ch),1}} L_{O_{(ch),1}}(z) + (-1)^{\sigma_{O_{(ch),2}}} W_{O_{(ch),2}} L_{O_{(ch),2}}(z)\right),
\end{multline*}
which finishes the proof \eqref{gev2} (recall that \( \kappa_1+\kappa_2=1\)). 

 Consider the second claim of the lemma. Given any $f$ defined on $\mathcal{V}$, we can use  \eqref{jack2} to get 
\begin{equation}
\label{jo-rest}
\jox f_{[X]} = (\jack f)_{[X]}- W_X^{1/2}f_{X_{(p)}}\delta^{(X)}.
\end{equation}
Since \( W_X^{1/2}m_{X_{(p)}}^{-1} = m_X^{-1} \), \eqref{gev3} follows from \eqref{jo-rest} applied to \( f=l \).
\end{proof}

We need to introduce  \emph{Green's functions} of $\jox$. They are defined by
\[
G(Y,X;z) \ddd \left\langle(\jox-z)^{-1}\delta^{(X)},\delta^{(Y)} \right\rangle,
\]
where \( X\in\mathcal V \) and \( Y\in\mathcal V_{[X]} \). Using Proposition~\ref{prop:gev} we can obtain the following conditional result.  Since \( \jox \) is a bounded operator, there exists \( R_{[X]}>0 \) such that \( \sigma(\jox)\subset\{|z|\leq R_{[X]}\} \).  Let $C_{[X]}$ denote the unbounded component of the complement of \( \sigma(\jox)\cup\supp \mu_1 \cup\supp \mu_2 \cup \{z: L_{\Pi(X_{(p)})}(z)=0\}\).

\begin{Prop}
\label{prop:Green}
If there exists \( R>0 \) such that \( l(z)\in\ell^2(\mathcal V) \) for \( |z|>R \), then \( l(z)\in\ell^2(\mathcal V) \) for all \( z \in C_{[X]}\),  and for all such \( z \) we have that
\begin{equation}
\label{Rz}
(\jox-z)^{-1}\delta^{(X)} = -m_Xl_{[X]}(z)/L_{\Pi(X_{(p)})}(z).
\end{equation}
In particular, \(G(Y,X;z)\) extends to a holomorphic function in \( C_{[X]} \) by
\begin{equation}
\label{GYX}
G(Y,X;z) = - \frac{m_X}{m_Y} \frac{L_{\Pi(Y)}(z)}{L_{\Pi(X_{(p)})}(z)}.
\end{equation}
\end{Prop}
\begin{proof} 
Let $z\in C_{[X]}$ be such that \( |z|>|R| \). If $X=O$, identity \eqref{Rz}  follows immediately from \eqref{gev2}. If $X\neq O$, it follows from \eqref{gev3}. Formula \eqref{GYX} for such \( z \) now follows from the definition and \eqref{gev1}. Moreover, since \( G(Y,X;z) \) is an analytic function of \( z\not\in\sigma(\jox) \) and \(   {L_{\Pi(Y)}(z)}/{L_{\Pi(X_{(p)})}(z)} \) is analytic in $z\notin \supp \mu_1 \cup\supp \mu_2 \cup \{z: L_{\Pi(X_{(p)})}(z)=0\}$, the full claim follows by analytic continuation.
\end{proof}

There are other functions that satisfy algebraic identities similar to \eqref{gev2}.  To introduce them, we first recall \eqref{sad15} and \eqref{An0}. Set
\begin{equation}
\label{Lambdas}
\Lambda_Y^{(k)}(z) \ddd m_Y^{-1}\,A_{\Pi(Y)}^{(k)}(z), \quad k\in\{0,1,2\}. 
\end{equation}
Observe that \( \Lambda_O^{(0)} = 0 \). For any functions \( f,g \) on \( \mathcal V \) and a fixed vertex \( Z\in\mathcal V \) we introduce  new function on $\mathcal V$ by
\begin{equation}
\label{commut}
[f,g]^{(Z)} \ddd f_Zg - fg_Z\,.
\end{equation}
We call it the commutator of functions \( f,g \) with respect to the vertex \( Z \).

\begin{Prop}
\label{prop:2.1.4}
The following algebraic identities hold
\[
\jack \Lambda^{(k)}(z) = z\Lambda^{(k)}(z) - \kappa_{3-k}A_{\vec e_k}^{(k)}(z)\delta^{(O)},
\]
for each \( k\in\{1,2\} \), as well as
\[
\jack\Lambda^{(0)}(z) = z\Lambda^{(0)}(z).
\]
Furthermore, let \( X\neq O \). Then, for any \( k,l\in\{0,1,2 \} \) it holds that
\begin{equation}\label{1sd4}
\jox \big[\Lambda^{(k)}(z),\Lambda^{(l)}(z)\big]^{(X_{(p)})}_{[X]} = z\big[\Lambda^{(k)}(z),\Lambda^{(l)}(z)\big]^{(X_{(p)})}_{[X]}.
\end{equation}
\end{Prop}
\begin{proof}
We can repeat the proof of Proposition~\ref{prop:gev} with \( L_\n(z) \) replaced by \( A_\n^{(k)}(z) \), \( k\in\{1,2\} \), and using \eqref{3sd1} instead of \eqref{1.5} to get that
\[
z\Lambda^{(k)}(z) = \jack \Lambda^{(k)}(z) + \big(\kappa_1A_{(0,1)}^{(k)}(z) + \kappa_2A_{(1,0)}^{(k)}(z)\big)\delta^{(O)}.
\]
Since \( A_{\vec1-\vec e_k}^{(k)}(z) \equiv 0 \), the first claim follows. We further get from \eqref{An0} and \eqref{gev2} that
\begin{eqnarray*}
\jack \Lambda^{(0)}(z) &=& \widehat\mu_1(z) \jack \Lambda^{(1)}(z) + \widehat\mu_2(z) \jack \Lambda^{(2)}(z) - \jack l(z) \\
& = & z\Lambda^{(0)}(z) - \left(\kappa_2 A_{(1,0)}^{(1)}(z) \, \widehat\mu_1(z) + \kappa_1 A_{(0,1)}^{(2)}(z) \, \widehat\mu_2(z) - L_{\vec\varkappa}(z)\right)\delta^{(O)}. 
\end{eqnarray*}
Since  \( A_{(1,0)}^{(1)}(z) = \|\mu_1\|^{-1} \) and \( A_{(0,1)}^{(1)}(z) = \|\mu_2\|^{-1} \) by \eqref{n_2}, the second claim follows from \eqref{LPO}.  To prove the third claim, observe that
\[
\jox\Lambda_{[X]}^{(k)}(z) = z\Lambda_{[X]}^{(k)}(z) - W_X^{1/2}\Lambda_{X_{(p)}}^{(k)}(z)\delta^{(X)}
\]
by \eqref{jo-rest}. The desired identity \eqref{1sd4} now easily follows from the definition~\eqref{commut}.
\end{proof}

\begin{Rem}
The relations of Proposition~\ref{prop:2.1.4} should be regarded as  algebraic identities and we do not claim that the functions involved belong to the Hilbert space \( \ell^2(\mathcal V) \) for any given \( z \).
\end{Rem}

The spectral theory of Jacobi matrices \eqref{jack2} under the sole condition \eqref{as-5} is currently  beyond our reach.  In Part~\ref{part3}, however, we consider a large class of multiorthogonal  systems, known as Angelesco systems, for  which this analysis is possible. 

\section{Appendix to Part~\ref{part2}}
\label{s:ap2}

In Part~\ref{part3}, we will explain that the so-called  Angelesco systems  generate bounded and self-adjoint Jacobi matrices.   In the current  appendix, we show that Nikishin systems,  see Section~\ref{s:ap1}, do not generate bounded Jacobi matrices, in general. We need some notation first.  Recall that a measure \( \mu \) supported on an interval \( \Delta=[\alpha,\beta] \) is called a \emph{Szeg\H{o} measure} if
\begin{equation}
\label{Gmu}
G(\mu) \ddd \exp\left\{\frac1\pi\int_{\Delta}\frac{\log\mu^\prime(x)dx}{\sqrt{(x-\alpha)(\beta-x)}}\right\}>0,
\end{equation}
where \( d\mu(x)=\mu^\prime(x)dx+d\mu_\mathrm{sing}(x) \) and \( \mu_\mathrm{sing} \) is singular to Lebesgue measure.

\begin{Thm}
\label{thm:nik-b}
Let \( \vec\mu \) be a Nikishin system \eqref{nik-sys}--\eqref{dtd1} and \( \{ b_{\n,1},b_{\n,2},a_{\n,1},a_{\n,2}\}_{\n\in\Z_+^2} \) be the corresponding recurrence coefficients, see \eqref{1.5}--\eqref{1.7}. Then, there exists a constant \( C_{\vec\mu} \) such that
\begin{equation}
\label{nik-bound}
\sup_{\n\in\Z_+^2} |b_{\n,i}| \leq C_{\vec\mu}, \quad\sup_{\n\in\Z_+^2:n_2\leq n_1~\text{or}~n_2\geq n_1+2} |a_{\n,i}| \leq C_{\vec\mu} 
\end{equation}
for any \( i\in\{1,2\} \). Assume further that the measures \( \mu_1\) and \(\tau \) are Szeg\H{o} measures. Then,
\begin{equation}
\label{nik-unbound}
\lim_{n\to\infty} a_{(n,n+1),1} = -\infty \quad \text{and} \quad \lim_{n\to\infty} a_{(n,n+1),2} = \infty.
\end{equation}
\end{Thm}

It is conceivable that  the Szeg\H{o} condition for the measures can be relaxed.  However, we assume it to simplify the proof. Our result shows that even for nice measures \( \mu_1,\tau \) the corresponding Nikishin system does not generate a bounded Jacobi matrix. In the remaining part of this section, we prove Theorem~\ref{thm:nik-b}.

\begin{Lem}
\label{lem:n3}
There exists a constant \( C_{\vec\mu} \) such that
\[
\sup_{\n\in\Z_+^2,~i\in\{1,2\}} |b_{\n,i}| \leq C_{\vec\mu}. 
\]
\end{Lem}
\begin{proof}
We continue to use notation from Section~\ref{s:ap1}. The following argument is taken from \cite{rocha1}. Divide recursion relations \eqref{1.7} by \( xP_\n(x) \) and integrate over a contour \( \Gamma \) that encircles \( \Delta_1 \) in positive direction to get
\[
\frac1{2\pi\ic}\int_\Gamma\left(1-\frac{P_{\n+\vec e_i}(z)}{P_\n(z)}\right)dz = b_{\n,i} + \frac1{2\pi\ic}\int_\Gamma\sum_{k=1}^2 a_{\n,k}\frac{P_{\n-\vec e_k}(z)}{P_\n(z)}dz.
\]
The last integral is zero by Cauchy theorem applied at infinity. Therefore,
\[
|b_{\n,i}| \leq \frac1{2\pi}\int_\Gamma\left| \frac{P_{\n+\vec e_i}(z)}{P_\n(z)}\right||dz|.
\]
It is known that the zeroes of \( P_{\n+\vec e_i}(z) \) and \( P_\n(z) \) interlace. Indeed, this follows from \cite[Theorem~2.1]{as1}, see also \cite[Corollary~1]{u1}, since it was established in \cite{g_b} that \( \{1,\widehat\tau\} \) is an AT system on \( \Delta_1 \) relative to \( \mu_1 \), see also \cite[page 782]{g_s1}. Thus, it holds that
\begin{equation}
\label{rat-int}
\frac{P_{\n+\vec e_i}(z)}{P_\n(z)} = (z-x_{\n+\vec e_i,1})(z-x_{\n+\vec e_i,|\vec n|+1}) \sum_{l=1}^{|\n|}\frac{c_{\n,l}}{z-x_{\n,l}},
\end{equation}
where \( c_{\n,l}>0 \) and \( \sum_{l=1}^{|\n|}c_{\n,l} = 1 \), and \( x_{\vec m,1}<x_{\vec m,2}<\cdots<x_{\vec m,|\vec m|}\) are the zeroes of \( P_{\vec m}(x) \), which all belong to \( \Delta_1 \). That shows boundedness of \( |P_{\n+\vec e_i}(z)/P_\n(z)| \) on \( \Gamma \) independently of \( \n \) and therefore proves the desired claim.
\end{proof}

Let \( r_{\n,1}(z),r_{\n,2}(z) \) be polynomials from Lemma~\ref{lem:n2}.

\begin{Lem}
\label{lem:n4}
If multi-indices \( \n \) and \( \n+\vec e_i \) both belong to the region \( \{(n_1,n_2):~n_2\leq n_1+1\} \), then the zeroes of \( r_{\n,1}(z) \) and \( r_{\n+\vec e_i,1}(z) \) interlace. Similarly, if multi-indices \( \n \) and \( \n+\vec e_i \) both belong to the region \( \{(n_1,n_2):~n_2\geq n_1+1\} \), then the zeroes of \( r_{\n,2}(z) \) and \( r_{\n+\vec e_i,2}(z) \) interlace.
\end{Lem}
\begin{proof}
The first claim was shown in \cite[Theorem~2.1]{agr}. The proof of the second claim is identical provided  one knows that  the functions \( AR_{\n,2}(z) + BR_{\n+\vec e_i,2}(z) \) have no more than \( n_1+1 \) zeroes in $\R\backslash \Delta_1$ all of which are simple (\(A,B\) are arbitrary real constants). The last property can be established exactly as in Lemma~\ref{lem:n2}, where  the cases \(A=0,B=1\) and \(A=1,B=0\) were considered.
\end{proof}

Recall that if \( q_n(x;\mu) \) is the \( n \)-th monic orthogonal polynomial with respect to measure \( \mu \) on the real line, then $q_n(x;\mu)$ is the unique minimizer of the following variational problem:
\begin{equation}
\label{var-prop}
\int q_n^2(x;\mu)d\mu(x) = \min\left\{\int q^2(x)d\mu(x):~q(x)=x^n+q_{n-1}x^{n-1}+\cdots+q_1x+q_0, \, \{q_j\}\in \R\right\}.
\end{equation}
\begin{Lem}
\label{lem:n5}
We have a bound
\[
\sup_{i\in\{1,2\},~\n\in\Z_+^2:n_2\leq n_1~\text{or}~n_2\geq n_1+2} |a_{\n,i}| \leq C_{\vec\mu}.
\]
\end{Lem}
\begin{proof}
As shown in Lemma~\ref{lem:n2}, it holds that
\[
h_{\n,1} = \int P_\n^2(x)\frac{d\mu_1(x)}{r_{\n,1}(x)}
\]
when \( n_2\leq n_1+1 \). Recall that the monic polynomials \( r_{\n,1}(x) \) and \( r_{\n-\vec e_1,1}(x) \) have degree \( n_2 \) and all their zeroes belong to \( \Delta_\tau \) when \( n_2\leq n_1 \) by Lemma~\ref{lem:n2}. Let \( 0<l<L \) be constants given by
\begin{equation}
\label{lL}
l^{-1} \ddd \max\{|x-y|:~x\in\Delta_1,~~y\in\Delta_\tau\} \quad \text{and} \quad L^{-1} \ddd \min\{|x-y|:~x\in\Delta_1,~~y\in\Delta_\tau\}.
\end{equation}
Then, when \( n_2\leq n_1 \), it follows from Lemma~\ref{lem:n4} that \( lr_{\n-\vec e_1,1}(x) \leq L r_{\n,1}(x) \) for any \( x\in\Delta_1>\Delta_\tau \). Let \( \gamma \) be the midpoint of \( \Delta_1 \) and \( |\Delta_1| \) be its length. Using \eqref{var-prop}, we get that 
\begin{eqnarray*}
h_{\n-\vec e_1,1} &\geq& \frac4{|\Delta_1|^2}\int (x-\gamma)^2P_{\n-\vec e_1}^2(x)\frac{d\mu_1(x)}{r_{\n-\vec e_1,1}(x)} \geq \frac4{|\Delta_1|^2}\frac lL \int (x-\gamma)^2P_{\n-\vec e_1}^2(x)\frac{d\mu_1(x)}{r_{\n,1}(x)} \\
& \geq &\frac4{|\Delta_1|^2}\frac lL \min\left\{ \int q^2(x)\frac{d\mu_1(x)}{r_{\n,1}(x)}:~q(x) = x^{|\n|} + \cdots \right\} = \frac{4}{|\Delta_1|^2} \frac lL h_{\n,1}.
\end{eqnarray*}
Therefore, it follows from \eqref{an-hn} that
\begin{equation}\label{1sd8}
|a_{\n,1}| = |h_{\n,1}/h_{\n-\vec e_1,1}| \leq (|\Delta_1|^2L)/(4l), \quad n_2\leq n_1.
\end{equation}
Furthermore, we get from recursion relations \eqref{1.7} that
\begin{equation}
\label{rat-Pns}
x - b_{\n,i}= \frac{P_{\n+\vec{e}_i}(x)}{P_{\n}(x)} + a_{\n,1}\frac{P_{\n-\vec{e}_1}(x)}{P_{\n}(x)}+ a_{\n,2}\frac{P_{\n-\vec{e}_2}(x)}{P_{\n}(x)}.
\end{equation}
Take \( x=\beta_1+1 \), where \( \Delta_1=[\alpha_1,\beta_1] \). The interlacing property used in Lemma~\ref{lem:n3}, see \eqref{rat-int}, implies that the ratios of polynomials in the above formula are positive and bounded above and away from zero independently of \( \n \). Thus, it follows from Lemma~\ref{lem:n3} and \eqref{1sd8} that 
\[
|a_{\n,2}| \leq C_{\vec\mu}, \quad n_2\leq n_1,
\]
for some constant \( C_{\vec\mu} \) independent of \( \n \), which is not necessarily the same as in Lemma~\ref{lem:n3}. The proof in the case \( n_2\geq n_1+2 \) is absolutely analogous: we first use Lemmas~\ref{lem:n2} and~\ref{lem:n4} to show boundedness of \( a_{\n,2} \) and then use recurrence relations \eqref{1.7} and Lemma~\ref{lem:n3} to deduce boundedness of \( a_{\n,1} \).
\end{proof}

We are left with proving \eqref{nik-unbound}. To proceed, let us recall some results from  \cite{St00}. Consider the function
\[ 
\psi(z) = z + \sqrt{z^2-1},
\] 
which  is the conformal map of \( \overline\C\setminus [-1,1] \) onto \( \overline\C\setminus\{|z|\leq1\} \) such that \( \psi(\infty)=\infty \) and \( \psi^\prime(\infty)>0 \). Let \( \mu \) be a Szeg\H{o} measure on \( [-1,1] \) and \( \{a_{2n,i}\}_{i=1}^{2n}\subset\overline\C\setminus [-1,1] \) define  a sequence  of multi-sets of complex numbers that are conjugate-symmetric and satisfy
\begin{equation}
\label{2sd1}
\lim_{n\to\infty} \sum_{i=1}^{2n}\left(1-|\psi(a_{2n,i})|^{-1}\right) = \infty\,.
\end{equation}
We emphasize that the elements in each multi-set $\{a_{2n,i}\}_{i=1}^{2n}$ can be equal to each other and  some of them can be equal to $\infty$. Let \( m_n \) be the number of finite elements in \( \{a_{2n,i}\}_{i=1}^{2n} \). Set 
\[ 
w_{2n}(z) \ddd \prod_{i=1}^{2n}(1-z/a_{2n,i}) \quad \text{and} \quad \widetilde w_{2n}(z) \ddd \prod_{|a_{2n,i}|<\infty}(z-a_{2n,i}),
\] 
which are polynomials of degree \( m_n\leq2n \) (\(\widetilde w_{2n}(z) \) is the monic renormalization of \( w_{2n}(z) \)). Conjugate-symmetry of $\{a_{2n,i}\}_{i=1}^{2n}$ guarantees that $w_{2n}(z)$ is real on the real line. Notice that \(w_{2n}(z)\equiv1\) when \( a_{2n,i}=\infty \) for all \( i\in\{1,\ldots,2n\} \). If \(\gamma_n \) is the leading coefficient of the \( n \)-th polynomial orthonormal with respect to the measure \( |w_{2n}(x)|^{-1}d\mu(x) \), then
\begin{equation}
\label{gamman2}
\gamma_n^{-2} \ddd \min\left\{\int q^2(x)\frac{d\mu(x)}{|w_{2n}(x)|}:q(x)=x^n+\cdots\right\}\,,
\end{equation}
see \eqref{var-prop}. It was shown in \cite[Corollary 1]{St00} that
\[
\lim_{n\to\infty} \gamma_n^{-2}2^{2n}\prod_{|a_{2n,i}|<\infty}\left(\frac{\psi(a_{2n,i})}{2a_{2n,i}} \right)= 2G(\mu),
\]
where \( G(\mu) \) was introduced in \eqref{Gmu}. Furthermore, if \( \widetilde\gamma_n \) is defined via \eqref{gamman2} with \( |w_{2n}(x)|^{-1}d\mu(x) \) replaced by \( |\widetilde w_{2n}(x)|^{-1}d\mu(x) \). Then it clearly holds that
\begin{equation}
\label{stahl0}
\lim_{n\to\infty}  \widetilde\gamma_n^{-2} 2^{2n-m_n}\prod_{i=1}^{m_n}{\psi(a_{2n,i})} = 2G(\mu)\,.
\end{equation}

More generally, let $\nu$ be a Szeg\H{o} measure on an interval $\Delta=[\alpha,\beta]$ and $d_n(z)$ be a monic polynomial of degree $m_n\le 2n$ with all its zeroes belonging to an interval $\Delta^*$ such that $\Delta^*\cap \Delta=\emptyset$. Define 
\begin{equation}
\label{Deltanur1}
\Delta_n(\nu,d_n) \ddd  \min\left\{\int q^2(x)\frac{d\nu(x)}{|d_n(x)|}:q(x)=x^n+\cdots\right\}\,.
\end{equation}
By rescaling the variables  as \( x(s) = |\Delta|(s+1)/2 + \alpha \), we get from \eqref{stahl0} that
\begin{equation}
\label{stahl}
\lim_{n\to\infty} \Delta_n(\nu,d_n)  (4/|\Delta|)^{2n-m_n}\prod_{i=1}^{m_n}{\psi(s_{2n,i})} = |\Delta| G(\nu),
\end{equation}
where \( \{x(s_{n,i})\}_{i=1}^{m_n} \) are the zeroes of \( d_n(x) \). 

We will need the following auxiliary statement.
\begin{Lem}If $G(\tau)>-\infty$, then $G(\tau_{d})>-\infty$. That is, the dual measure \( \tau_d \) is a Szeg\H{o} measure when \( \tau \) is a Szeg\H{o} measure. \label{2sd3}
\end{Lem}
\begin{proof}
 It follows from Proposition~\ref{prop:Markov} further below that \( \tau_d^\prime(x) \) exists almost everywhere on \( \Delta_\tau \) and
\[
\tau_d^\prime(x) = |\widehat\tau_+(x)|^{-2}\tau^\prime(x) \geq \big(\pi^2\tau^\prime(x)\big)^{-1}
\]
for a.e. \( x\in\Delta_\tau \). Thus, if we let \( w_\tau(x) \ddd \sqrt{(x-\alpha_\tau)(\beta_\tau-x)} \), it holds that
\begin{eqnarray*}
G(\tau_d) &\geq& \exp\left\{\frac1\pi\int\log\left(\frac1{w_\tau(x)\tau^\prime(x)}\right)\frac{dx}{w_\tau(x)} + \frac1\pi\int\log\left(\frac{w_\tau(x)}{\pi^2}\right)\frac{dx}{w_\tau(x)} \right\} \\ 
&=& C_\tau \exp\left\{\frac1\pi\int\log\left(\frac1{w_\tau(x)\tau^\prime(x)}\right)\frac{dx}{w_\tau(x)}\right\} \geq \frac{ C_\tau}{\int \tau^\prime(x)dx}>0,
\end{eqnarray*}
where we used Jensen's inequality at the last step. 
\end{proof}

\begin{Lem}
\label{lem:n6}
Assume that \( \tau \) is a Szeg\H{o} measure. Then, there exists a constant \( C_{\vec{\mu}} \) such that
\[
C_{\vec{\mu}}^{-1} \leq h_{(n,n+1),2}/h_{(n,n+1),1} \leq C_{\vec{\mu}}
\]
for all $n\in \mathbb{N}$.
\end{Lem}
\begin{proof}
Let \( \n=(n,n+1) \) and \( l,L \) be as in \eqref{lL}. It follows from Corollary~\ref{cor:n2} and Lemma~\ref{lem:n2} that
\begin{equation}
\label{lhL}
lh_{\n,j} \leq |P_\n(x)R_{\n,j}(x)/r_{\n,j}(x)| \leq Lh_{\n,j}, \quad x\in\Delta_\tau,
\end{equation}
(recall that \( h_{\n,j}>0 \) for such \( \n \) by Corollary~\ref{cor:n1}). We further get from Lemma~\ref{lem:n1} that
\[
\int x^k r_{\n,2}(x) \frac{R_{\n,2}(x)d\tau_d(x)}{r_{\n,2}(x)} = 0
\]
for \( k\leq n-1 \) and \( \deg\,r_{\n,2} = n \), where the measure \( R_{\n,2}(x)d\tau_d(x)/r_{\n,2}(x) \) is non-negative on \( \Delta_\tau \), see \eqref{Rrj}. Therefore,
\begin{equation}
\label{Rn2-min}
\int x^n R_{\n,2}(x)d\tau_d(x) = \int r_{\n,2}^2(x) \frac{R_{\n,2}(x)d\tau_d(x)}{r_{\n,2}(x)} = \min_{q(x)=x^n+\cdots} \int q^2(x) \frac{R_{\n,2}(x)d\tau_d(x)}{r_{\n,2}(x)},
\end{equation}
where we used \eqref{var-prop} for the last equality. One can readily check that
\[
\min_{q(x)=x^n+\cdots}  \int q^2(x)d\tau_1(x) \leq \min_{q(x)=x^n+\cdots}  \int q^2(x)d\tau_2(x)
\]
if \( \tau_1(B)\leq\tau_2(B) \) for all Borel sets \( B \). Hence, it follows from \eqref{lhL} that
\[
l^2h_{\n,2}\Delta_n(\tau_d,P_\n^*) \leq \min_{q(x)=x^n+\cdots} \int q^2(x) \frac{R_{\n,2}(x)d\tau_d(x)}{r_{\n,2}(x)} \leq L^2h_{\n,2}\Delta_n(\tau_d,P_\n^*),
\]
where \( \Delta_n(\tau_d,P_\n^*) \) is defined via \eqref{Deltanur1}, \( P_\n^*(x) = P_\n(x)/(x-x_{\n,2n+1}) \), and we denote the zeroes of \( P_\n(x) \) by \( x_{\n,1}<\ldots <x_{\n,2n+1} \) (we stripped one zero from \( P_\n(x) \) since \( \deg\,P_\n=2n+1>2n=2\deg r_{\n,2} \)). Similarly, we get that
\[
lh_{\n,1}\Delta_{n+1}(\tau,P_\n) \leq \min_{q(x)=x^{n+1}+\cdots} \int q^2(x) \frac{R_{\n,1}(x)d\tau(x)}{r_{\n,1}(x)} \leq Lh_{\n,1}\Delta_{n+1}(\tau,P_\n)
\]
(here, we do not need to strip zeroes from \( P_\n(x) \) since \( \deg\,P_\n=2n+1<2(n+1)=2\deg r_{\n,1} \)). Then, it follows from the last claim of Lemma~\ref{lem:n1} (the equality of the integrals), \eqref{Rn2-min}, and a similar formula for \( R_{\n,1}(x) \) that
\begin{equation}
\label{rathn1}
\frac1{\|\tau\|}\frac l{L^2} \frac{\Delta_{n+1}(\tau,P_\n)}{\Delta_n(\tau_d,P_\n^*)} \leq \frac{h_{\n,2}}{h_{\n,1}} \leq \frac1{\|\tau\|}\frac L{l^2} \frac{\Delta_{n+1}(\tau,P_\n)}{\Delta_n(\tau_d,P_\n^*)}.
\end{equation}
By Lemma \ref{2sd3}, we know that $\tau_d$ is a Szeg\H{o} measure and we can apply formula \eqref{stahl} to control the ratios in the left-hand and right-hand sides of \eqref{rathn1}. We get
\begin{equation}
\label{rathn2}
C_{\Delta_\tau,\Delta_1}\frac{G(\tau)}{G(\tau_d)} \leq \liminf_{n\to\infty} \frac{\Delta_{n+1}(\tau,P_\n)}{\Delta_n(\tau_d,P_\n^*)} \leq \limsup_{n\to\infty} \frac{\Delta_{n+1}(\tau,P_\n)}{\Delta_n(\tau_d,P_\n^*)} \leq C'_{\Delta_\tau,\Delta_1}\frac{G(\tau)}{G(\tau_d)},
\end{equation}
where  $C_{\Delta_\tau,\Delta_1}$ and $C'_{\Delta_\tau,\Delta_1}$  depend only on the intervals $\Delta_\tau$ and $\Delta_1$. The desired claim now follows from \eqref{rathn1} and \eqref{rathn2}.
\end{proof}

\begin{Lem}
\label{lem:n7}Assume that \( \tau \) is a Szeg\H{o} measure. Then,
\[
\lim_{n\to\infty} h_{(n,n),1}/h_{(n,n),2} = \infty.
\]
\end{Lem}
\begin{proof}
Let \( \n=(n,n) \). Similarly to \eqref{Rn2-min}, it follows from Lemma~\ref{lem:n1}, \eqref{Rrj}, and \eqref{var-prop} that
\[
h_{\n,2} = - \int x^n R_{\n,1}(x)d\tau(x) = \min_{q(x)=x^n+\cdots} \int q^2(x) \frac{|R_{\n,1}(x)|d\tau(x)}{|r_{\n,1}(x)|} .
\]
As in the previous lemma, we get from \eqref{lhL} that
\[
0\leq \frac{h_{\n,2}}{h_{\n,1}} \leq L \min_{q=x^n+\cdots} \int q^2(x) \frac{d\tau(x)}{|P_\n(x)|} = L \Delta_n(\tau,P_\n). 
\]
Again, as in the previous lemma, let \( x(s)=|\Delta_\tau|(s+1)+\alpha_\tau \), \( \Delta_\tau=[\alpha_\tau,\beta_\tau] \). Then, we get from \eqref{stahl} that
\[
\lim_{n\to\infty} \Delta_n(\tau,P_\n) \prod_{i=1}^{2n}\psi(s_{\n,i}) = |\Delta_\tau|G(\tau),
\]
where \( x_{\n,i} = x(s_{\n,i}) \), \( i\in\{1,\ldots,2n\} \), are the zeroes of \( P_\n(x) \). Since \( \psi(x(s_{\n,i}))\geq\psi(x(\alpha_1))>1 \), it holds that \( \lim_{n\to\infty} \Delta_n(\tau,P_\n)=0 \), which finishes the proof of the lemma.
\end{proof}

\begin{Lem}
\label{lem:n8}
Assume that \( \mu_1 \) is a Szeg\H{o} measure. Then, there exists a constant \( C_{\vec\mu} \) such that
\[
C_{\vec\mu}^{-1} \leq h_{(n,n+1),1}/h_{(n,n),1} \leq C_{\vec\mu}
\]
holds for all $n\in \mathbb{N}$.
\end{Lem}
\begin{proof}
As shown in Lemma~\ref{lem:n2}, it holds that \( \deg r_{\n,1} = n_2 \) and
\[
h_{\n,1} = \int P_\n^2(x)\frac{d\mu_1(x)}{r_{\n,1}(x)} = \min_{q(x)=x^{|\n|}+\cdots}\int q^2(x)\frac{d\mu_1(x)}{r_{\n,1}(x)} = \Delta_{|\n|}(\mu_1,r_{\n,1})
\]
when \( n_2\leq n_1+1 \), where we also used property \eqref{var-prop} and definition \eqref{Deltanur1}. Let \( x(s)=|\Delta_1|(s+1)/2+\alpha_1 \), where \( \Delta_1=[\alpha_1,\beta_1] \). Then, it follows from \eqref{stahl} that
\[
\lim_{n\to\infty}\Delta_{2n}(\mu_1,r_{(n,n),1}) \big(4/|\Delta_1|\big)^{3n}\prod_{i=1}^n|\psi(s_{(n,n),i})| = |\Delta_1| G(\mu_1),
\]
where \( x_{\n,i}=x(s_{\n,i}) \), \( i\in\{1,\ldots,n\} \), are the zeroes of \( r_{\n,1}(x) \), and
\[
\lim_{n\to\infty}\Delta_{2n+1}(\mu_1,r_{(n,n+1),1})(4/|\Delta_1|\big)^{3n+1}\prod_{i=1}^{n+1}|\psi(s_{(n,n+1),i})| = |\Delta_1|G(\mu_1).
\]
Recall that according to Lemma~\ref{lem:n4}, the zeroes of $r_{\n,1}(x)$ and $r_{\n+\vec{e}_l,1}(x)$ interlace as long as both $\n$ and $\n+\vec{e}_l$ belong to the set $\{n_2\le n_1+1\}$. Thus,
\[
|\psi(x(\alpha_\tau))|^{-1} \prod_{i=1}^{n+1}|\psi(s_{(n,n+1),i})| \leq \prod_{i=1}^n|\psi(s_{(n,n),i})| \leq |\psi(x(\beta_\tau))|^{-1} \prod_{i=1}^{n+1}|\psi(s_{(n,n+1),i})|,
\]
where, as before, we write \( \Delta_\tau=[\alpha_\tau,\beta_\tau]<\Delta_1 \). Therefore, by combining  the previous estimates, we get that
\[
C^{-1}_{\vec\mu}\leq \liminf_{n\to\infty} \frac{h_{(n,n+1),1}}{h_{(n,n),1}} \leq \limsup_{n\to\infty} \frac{h_{(n,n+1),1}}{h_{(n,n),1}} \leq C_{\vec\mu}
\]
which yields the desired claim. 
\end{proof}

\begin{proof}[Proof of Theorem~\ref{thm:nik-b}]
The proofs of the claims in \eqref{nik-bound} are contained in Lemmas~\ref{lem:n3} and~\ref{lem:n5}. It further follows from \eqref{an-hn} that
\[
a_{(n,n+1),2} = \frac{h_{(n,n+1),2}}{h_{(n,n),2}} =  \left(\frac{h_{(n,n+1),2}}{h_{(n,n+1),1}} \right)\cdot \left(\frac{h_{(n,n),1}}{h_{(n,n),2}}\right) \cdot \left(\frac{h_{(n,n+1),1}}{h_{(n,n),1}}\right).
\]
Thus, the second claim of \eqref{nik-unbound} is a consequence of Lemmas~\ref{lem:n6}--\ref{lem:n8}. The first claim of \eqref{nik-unbound} now follows from the considerations adduced right after \eqref{rat-Pns}.
\end{proof}

\part{Jacobi matrices of Angelesco systems}
\label{part3}

In this part, we consider Angelesco systems \cite{Ang19}. These are systems \( \vec\mu=(\mu_1,\mu_2) \) that satisfy
\begin{equation}
\label{1.10}
\Delta_1\cap\Delta_2=\varnothing, \quad \Delta_i\ddd  \mathrm{ch}(\supp\, \mu_i)=[\alpha_i,\beta_i],
\end{equation}
where, as before, \( \mathrm{ch}(\cdot) \) stands for the convex hull of a set. Without loss of generality, we assume that $\Delta_1<\Delta_2$ (recall that  we write $E_1<E_2$ if two sets $E_1$ and $E_2$ satisfy \( \sup E_1 <\inf  E_2 \)). Note that $\Delta_1,\,\Delta_2$ is a system of two closed intervals separated by an open one. It will be convenient to use notation 
\begin{equation}
\label{mustar}
d\mu^\star(x) \ddd \chi_{\Delta_1}(x)d\mu_1(x) + \chi_{\Delta_2}(x)d\mu_2(x)
\end{equation}
to define the ``concatenation'' of measures \( \mu_1 \) and \( \mu_2 \).

It is known, see \cite[Appendix~A]{ady1}, that Angelesco systems satisfy not only \eqref{as-1} and \eqref{as-5}, but also \eqref{as-4}. In particular, Jacobi matrices \( \jack \) of such systems are bounded and self-adjoint. It is also known that \( l(z)\in\ell^2(\mathcal V) \) for \( |z|>R \) and some \( R>0 \), see \cite[Proposition~4.2]{ady1}. The function $L_{\n}(z)$  has no zeroes outside $\Delta_1\cup\Delta_2$ (see, e.g., Lemma \ref{lem:L}) for any $\n$. Therefore, \eqref{Rz} holds everywhere in \( \sigma(\jack)\cup\supp\,\mu_1\cup\supp\,\mu_2 \) as \( \sigma(\jack)\subset\R \).

\section{Poisson integrals}
\label{s:Pi}

Our primary working tool in studying spectral properties of \( \jack \) is the Green's function \( G(Y,X;z) \), whose boundary behavior we investigate via formula \eqref{GYX}. To ease referencing while doing so, we gather some well-known properties of functions harmonic in \( \C_+ \), the upper half-plane, in this section.

\begin{Prop}
\label{prop:Garnett}
Let \( v(z) \) be a function harmonic in \( \C_+ \) and such that
\begin{equation}
\label{sup-u}
\sup_{y>0} \int_\R |v(x+\ic y)|^pdx < \infty
\end{equation}
for some \( p\geq 1 \). Then, there exists a finite (generally signed) measure \( \mu \) on \( \R \) such that
\begin{equation}
\label{v-recover}
v(x+\ic y) = \int_\R P_z(t)d\mu(t), \quad P_z(t) \ddd \frac1\pi\Im\left(\frac1{t-z}\right),\quad z=x+iy\,,
\end{equation}
where \( P_z(t) \) is known as the Poisson kernel. The measure \( \mu \) is constructed as
\begin{equation}
\label{trace}
v(x+\ic y)dt  \cws d\mu(x) \quad \text{as} \quad y\to 0^+,
\end{equation}
where \( \cws \) denotes the weak$^*$ convergence of measures.  The limit 
\begin{equation}
\label{harm-lim}
\mu^\prime(x) = \lim_{y\to0^+} v(x+\ic y)
\end{equation}
exists for Lebesgue almost all \( x \) on the real line (the limit in \eqref{harm-lim} can be taken in non-tangential sense) and
\begin{equation}
\label{mu-decomp}
d\mu(x)=\mu^\prime(x)dx + d\mu_\mathrm{sing}(x),
\end{equation}
where \( \mu_\mathrm{sing} \) is singular to Lebesgue measure. For each \( p>1 \), \eqref{sup-u} holds if and only if \( \mu_\mathrm{sing}\equiv0 \) and \( \mu^\prime\in L^p(dx) \).
\end{Prop}
\begin{proof}
This proposition is a combination of Theorem~I.3.1,~I.3.5, and I.5.3 of \cite{Garnett}.
\end{proof}

Hereafter we use the following convention: for a closed interval \( \Delta \), we let \( \Delta^\circ \) be the corresponding open interval. We denote by \( \mathrm{DC}(I) \) the set of Dini-continuous functions on \( I\in\{\Delta,\Delta^\circ\} \) (see, e.g., p.105 in \cite{Garnett}).

\begin{Prop}
\label{prop:an}
Let   \( v(z) = \Im f(z) \) for some function $f(z)$ analytic in \( \C_+ \) which satisfies  \begin{equation}\label{2sd10} \lim_{y\to +\infty}f(iy)= 0.\end{equation}
\begin{itemize}
\item[(1)] If \( v(z) \) satisfies \eqref{sup-u} for some \( p>1 \), then so does \( f(z) \). \smallskip
\item[(2)]  Suppose \( v(z) \) satisfies \eqref{sup-u} with $p=1$, the measure \( \mu \), defined in \eqref{trace}, is absolutely continuous on some open, possibly unbounded, interval \( I \), and  \( \mu^\prime \in \mathrm{DC}(I) \), then \( f(z) \) extends continuously to \( I \) from \( \C_+ \).
\end{itemize}
\end{Prop}
\begin{proof} Given condition \eqref{2sd10}, we can write $f=v+i\widetilde v$, where $\widetilde v$ is the harmonic conjugate of $v$.
Now, the proof  follows by applying  a combination of Theorem~III.2.3 and Corollary~III.1.4 in \cite{Garnett}.
\end{proof}

The following result provides an integral representation of functions that are harmonic and  positive in $\C_+$.

\begin{Prop}
\label{prop:Her}
A function \( u(z) \) is positive harmonic in \( \C_+ \) if and only if
\begin{equation}
\label{u-pos}
u(x+\ic y) = by + \int_\R P_z(t)d\mu(t),
\end{equation}
where \( b\geq0 \) and \( \mu \) is a positive measure satisfying \( \int_{\R}(1+x^2)^{-1}d\mu(x)<\infty \). Given such \( u(z) \), the measure \( \mu \) can be obtained via \eqref{trace}.
\end{Prop}
\begin{proof}
These claims are contained in \cite[Theorem I.3.5]{Garnett}.
\end{proof}

The function \( m(z) \) belongs to \( \HN \), the Herglotz-Nevanlinna class, if it is holomorphic in \( \C_+ \) and has non-negative imaginary part there. 
Such functions allow the following unique integral representation \cite{Garnett}
\begin{equation}\label{df1}
m(z) = \frac{1}{\pi}\int_{\R}\left(\frac{1}{x-z} - \frac{x}{x^2+1}\right)\,d\mu(x)+ bz+\tilde a, \quad z\in \mathbb{C}_+,
\end{equation}
where $\tilde a \in\R$, and $b,\mu$ are as in \eqref{u-pos}. If \( m(z) \) has a holomorphic continuation to a punctured neighborhood of infinity (where its has a simple pole), the measure \( \mu \) is compactly supported and the above representation becomes
\begin{equation}
\label{hg}
m(z) = -\pi^{-1}\widehat\mu(z) + bz+a, \quad z\in \mathbb{C}_+,
\end{equation}
where $b \ge 0$, $a \in\R$, and \( \widehat\mu(z) \) is the Markov function of \( \mu \), see \eqref{markov}. Notice that 
\[\Im m(z) = by + \int_\R P_z(t)d\mu(t).\]
Motivated by \eqref{trace}, we shall set
\begin{equation}
\label{imp}
\Im m^+ \ddd\mu\,.
\end{equation}

We will be particularly interested in reciprocals $\widehat\mu^{-1}$ of Markov functions $\widehat \mu$. It follows straight from the definition  that \( \widehat\mu^{-1}\in \HN \). Since $\mu$ is positive and has compact support,  there exist a compactly supported positive measure \( \upsilon \) and a real number \( a \) such that
\begin{equation}
\label{recip-mark}
\widehat\mu^{-1}(z) = a + \|\mu\|^{-1} z - \widehat\upsilon(z).
\end{equation}
We called the measure $\upsilon$ dual to $\mu$, see \eqref{dual-measure}.  Let \( \mathrm{DC}_0(\Delta) \subset \mathrm{DC}(\Delta) \) be the subset of functions that vanish at the endpoints of a closed interval \( \Delta \).

\begin{Prop}
\label{prop:Markov}
Let \( \mu \)  be compactly supported non-negative measure and $\mu_{\rm sing}$ denote its singular part. It holds that
\begin{itemize}
\item[(1)] The traces \( \widehat\mu_{\pm}(x) \ddd \lim_{y\to 0^\pm} \widehat\mu(x+\ic y) \) exist and are finite almost everywhere on the real line. \medskip
\item[(2)] \( \mu^\prime(x) = -\pi^{-1}\Im\big(\widehat\mu_+(x)\big) \) almost everywhere on the real line. \medskip
\item[(3)]  \( \mu(\{E\}) = \lim_{y\to0^+} \ic y\widehat\mu(E+\ic y) \) and \( \supp\,\mu_\mathrm{sing} \subseteq \big\{x:~-\lim_{y\to0^+} \Im\big(\widehat\mu(x+\ic y) \big)=\infty \big\}. \) \medskip
\item[(4)] If \( \supp\,\mu=\Delta \), \( \mu \) is absolutely continuous, and \( \mu^\prime \in \mathrm{DC}_0(\Delta)\), then \( \widehat\mu(z) \) extends continuously to \( \R \) from \( \C_+ \) and from \( \C_- \). Moreover, $\widehat\mu_+(x)=
\overline{\widehat\mu_-(x)}$.
 \medskip
\item[(5)] If, in addition to assumptions in (4), we have \( \mu^\prime(x)>0 \) for  \( x\in\Delta^\circ \), then \( \widehat\mu_\pm(x) \neq 0 \), \( x\in \R \).
\end{itemize}
\end{Prop}
\begin{proof} (1) This claim follows from \cite[Theorem~I.5.3, Lemma~III.1.1, and Theorem~III.2.1]{Garnett}. \smallskip
\\
(2) The claim is a restatement of \eqref{harm-lim}. \smallskip
\\
(3) These statements can be found in \cite[Proposition~1]{pearson} and \cite[Proposition~2.3.12]{Simon2}.  \smallskip
\\
(4) This claim follows from Proposition~\ref{prop:an}(2) since \( \mu^\prime\in\mathrm{DC}(I) \) for any open interval \( I \) containing \( \Delta \). \smallskip
\\
(5) Since \( \Im\big(\widehat\mu_+(x)\big) = -\pi\mu^\prime(x) \) by claim (2), it is non-vanishing on \( \Delta^\circ \). Moreover, \( \Re\big(\widehat\mu_+(x)\big) = \widehat\mu(x) \) for \( x\not\in\Delta^\circ \) and therefore is monotonically decreasing there while also equal to zero at infinity. Thus, it is necessarily non-vanishing for \( x\not\in\Delta^\circ \).
\end{proof}

\section{Reference measures}

As we mentioned before, formula \eqref{GYX} is central to our analysis and therefore we need to study  the functions \( L_\n(z) \). Below, we shall often refer to the auxiliary lemmas proven in Section~\ref{s:ap3}. 

\begin{Lem}
\label{lem:Ln}
Assume that the measure \( \mu_k \) is supported  on \( \Delta_k \) and is absolutely continuous with \( \mu_k^\prime \in \mathrm{DC}_0(\Delta_k) \) and \( \mu_k^\prime(x)>0 \) for \( x\in\Delta_k^\circ \), \( k\in\{1,2\} \). Then, given \( \n\in \mathbb N^2 \), the function \( L_\n(z) \) extends continuously to the real line from \( \C_+ \) and, in particular, the function \( |L_\n(x)| \) is well-defined, continuous, and non-vanishing on the whole real line. 
\end{Lem}
\begin{proof}
It follows from \eqref{Ln-An} and Proposition~\ref{prop:an}(2) that \( L_\n(z) \) extends continuously to the real line from the upper and lower half-planes.  Actually, as \( L_{\n+}(x) \) and \( L_{\n-}(x) \) are complex-conjugates of each other, \( |L_\n(x)| \) is well-defied and continuous on all of $\C$. It follows from Lemma~\ref{lem:L}(3) that it is non-vanishing outside of \( \Delta_1^\circ\cup\Delta_2^\circ \). We further get from Proposition~\ref{prop:Markov}(2-4) that \( \Im L_{\n+}(x) = -\pi A_\n^{(k)}(x)\mu_k^\prime(x) \) on \( \Delta_k \). Thus,  \( |L_\n(x)| \) is non-vanishing outside of zeroes of \(  A_\n^{(k)}(x) \). However, we show in Lemma~\ref{lem:L}(4) that \( \Re L_{\n+}(E)\neq 0 \) for each such zero \( E \). 
\end{proof}

In the case of systems \( \vec  \mu \) satisfying conditions of Lemma~\ref{lem:Ln} we can introduce  \emph{``reference measures''} as
\begin{equation}
\label{omega-n}
|L_\n(x)|^{-2}d\mu^\star(x),
\end{equation}
where \( \mu^\star \) is the concatenated measure from \eqref{mustar}. When \( \vec  \mu \) is no longer smooth, we use the general theory of Herglotz-Nevalinna functions to introduce them. We start with a few definitions. Given $\xi\in(\beta_2,\alpha_1)$, define $D_{\n,\xi}(z)$ by \begin{equation}
\label{ur3}
D_{\n,\xi}(z) \ddd (-1)^{n_2}(z-\xi)A_\n^{(1)}(z)A_\n^{(2)}(z)
\end{equation}
and non-negative function $S_{\n,\xi}(x)$ by
\[
S_{\n,\xi}(x) \ddd |x-\xi|^{-1}\left(\chi_{\Delta_1}(x)\,|A_\n^{(2)}(x)|^{-1} + \chi_{\Delta_2}(x)\,|A_\n^{(1)}(x)|^{-1}\right).
\]
Let \( E_\n \) be the set of zeroes of \( A_\n^{(1)}(z)A_\n^{(2)}(z) \). For each \( E\in E_\n \), we define an auxiliary measure \( \nu_{\n,E} \) by
\begin{equation}
\label{nunE}
d\nu_{\n,E}(x) \ddd \frac{D_{\n,\xi}(x)A_\n^{(1)}(x)}{(x-E)^2} d\mu_1(x) +  \frac{D_{\n,\xi}(x)A_\n^{(2)}(x)}{(x-E)^2} d\mu_2(x)\,.
\end{equation}
This is a well-defined measure on \( \Delta_1\cup\Delta_2 \)  since each \( E\in\Delta_k \) is a double zero of the respective numerator. In  Lemma~\ref{lem:L}(2), we prove that $\nu_{\n,E}(x)$ is in fact positive provided that $\xi\in (\beta_2,\alpha_1)$. Recall that $\HN$ stands for the Herglotz-Nevanlinna class.
\begin{Prop}
\label{prop:RefM}
Given \( \n\in\mathbb N^2 \), it holds that \( (D_{\n,\xi}L_\n)^{-1}\in\HN \) for any \( \xi\in(\beta_2,\alpha_1) \). There exists a non-negative measure \( \omega_\n \) (the reference measure) supported on \( \Delta_1\cup\Delta_2 \) such that 
\begin{equation}
\label{DnLnOmn}
\frac1{D_{\n,\xi}(z)L_\n(z)}  = \int_\R\frac{S_{\n,\xi}(x)d\omega_\n(x)}{x-z} +\sum_{E:~D_{\n,\xi}(E)=0} \frac{\zeta_{\n,\xi}(E)}{E-z} + a_{\n,\xi} + b_{\n,\xi} z,\quad z\in \C_+\,,
\end{equation}
where \( a_{\n,\xi}\in\R \), \( b_{\n,\xi}>0 \) and the numbers \( \zeta_{\n,\xi}(E) \ddd -(D_{\n,\xi}^\prime(E)L_{\n+}(E))^{-1} \) are well-defined and positive for every zero \( E \) of \( D_{\n,\xi}(x) \) (in fact, \(\zeta_{\n,\xi}(E) = \|\nu_{\n,E}\| - \nu_{\n,E}(\{E\})\) for each \( E\in E_\n \)). Measure $\omega_\n$ has no atoms at the zeroes of $ D_{\n,\xi}(z)$. Moreover, if \( \vec \mu \) satisfies the conditions of Lemma~\ref{lem:Ln}, then \( d\omega_\n(x) \) is equal to \eqref{omega-n}.
\end{Prop}
\begin{proof}
It is shown in Lemma~\ref{lem:L}(2) that the linear form \( D_{\n,\xi}(x)Q_\n(x) \) is, in fact, a non-negative measure on \( \Delta_1\cup\Delta_2 \) for any \( \xi\in(\beta_1,\alpha_2) \), and, according to Lemma~\ref{lem:L}(1), the Markov function of this measure is equal to \( D_{\n,\xi}(z)L_\n(z) \). Therefore, \( (D_{\n,\xi}L_\n)^{-1}\in\HN \) and we get from \eqref{recip-mark} that there exist constants \( b_{\n,\xi}>0 \), \( a_{\n,\xi}\in \R \), and a non-negative measure \( \upsilon_{\n,\xi} \) such that
\begin{equation}
\label{up-nxi}
\big(D_{\n,\xi}(z)L_\n(z)\big)^{-1}  - a_{\n,\xi} - b_{\n,\xi} z  = -\pi^{-1}\widehat\upsilon_{\n,\xi}(z).
\end{equation}
The measure \( \upsilon_{\n,\xi} \) has a point mass at \( \xi \) since \( D_{\n,\xi}(z)L_\n(z) \) is holomorphic around \( \xi \) and has a simple zero there. The mass at \( \xi \) can be computed via Proposition~\ref{prop:Markov}(3), where one needs to observe that  \( D_{\n,\xi}^\prime(\xi)L_\n(\xi)<0 \) because Markov functions have negative derivatives on the real line away from the support of the defining measure.  If $E\in E_{\n}$, it follows from Proposition~\ref{prop:Markov}(3) and Lemma~\ref{lem:L}(4) that
\begin{eqnarray}
\nonumber
\upsilon_{\n,\xi}(\{E\}) &=& -\pi \lim_{y\to0^+}\Bigl(iy(D_{\n,\xi}(E+\ic y)L_\n(E+\ic y))^{-1}\Bigr)  \\
\label{up-E}
&=& -\pi(D_{\n,\xi}^\prime(E) L_{\n+}(E))^{-1} = \pi(\|\nu_{\n,E}\|-\nu_{\n,E}(\{E \}))^{-1}>0.
\end{eqnarray}
Hence, the reference measure \( \omega_\n \) introduced in the proposition is equal to
\[
d\omega_\n(x) = \pi^{-1}S_{\n,\xi}^{-1}(x)d\upsilon_{\n,\xi}(x) - \sum_{E:~D_\n(\xi;E)=0} S_{\n,\xi}^{-1}(E)\zeta_{\n,\xi}(E)d\delta_E(x)
\] 
and it has no atoms at zeroes of $D_\n(\xi;x)$.
To show that \( \omega_\n \) is indeed independent of \( \xi \), let us derive an explicit expression for it when \( \vec \mu \) satisfies the condition of Lemma~\ref{lem:Ln}. We know from Lemma~\ref{lem:L}(1,2) and Proposition~\ref{prop:Markov}(2-4) that
\begin{equation}
\label{up-1}
\Im(L_{\n+}(x)) = -\pi \left( \chi_{\Delta_1}(x) A_\n^{(1)}(x)\mu_1^\prime(x) + \chi_{\Delta_2}(x) A_\n^{(2)}(x)\mu_2^\prime(x) \right).
\end{equation}
It further follows from Lemma~\ref{lem:Ln} that \( |L_\n(x)| \) is continuous and non-vanishing on the real line. Therefore, for any \( x\not\in E_\n \) we get that
\begin{equation}
\label{up-2}
-\pi^{-1}\Im(\widehat\upsilon_{\n,\xi+}(x)) = - |L_\n(x)|^{-2} D_{\n,\xi}^{-1}(x) \, \Im(L_{\n+}(x)) < \infty.
\end{equation}
Thus, Proposition~\ref{prop:Markov}(3) yields that the support of the singular part of \( \upsilon_{\n,\xi} \) is a subset of the zeroes of \( D_{\n,\xi}(z) \) (actually is equal to it by what precedes). Hence, in this case \( \omega_\n \) is an absolutely continuous measure and it follows from Proposition~\ref{prop:Markov}(2) that
\[
d\omega_\n(x) = \pi^{-1}S_{\n,\xi}^{-1}(x)\upsilon_{\n,\xi}^\prime(x)dx = -\pi^{-2}S_{\n,\xi}^{-1}(x)\Im\big(\widehat\upsilon_{\n,\xi+}(x)\big)dx = |L_\n(x)|^{-2}d\mu^\star(x)
\]
as claimed, where we used \eqref{ur3}, \eqref{up-1}, \eqref{up-2}, and Lemma~\ref{lem:L}(2) to get the last equality.

Let \( \vec\mu \) be any Angelesco system and \( \{\vec\mu_m\} \) be a sequence of Angelesco systems satisfying conditions of Lemma~\ref{lem:Ln} and such that \( \mu_{m,l} \cws \mu_l \) as \( m\to\infty \), \( l\in\{1,2\} \).  Since the moments of \( \mu_{m,l} \) converge to the corresponding moments of \( \mu_l \),  MOPs with respect to \( \vec\mu_m \) converge uniformly on compact subsets of \( \C \) to the corresponding MOP with respect to \( \vec\mu \).  Thus, linear forms \eqref{Qn} with respect to \( \vec\mu_m \) converge in the weak$^*$ topology to the corresponding linear form with respect to \( \vec\mu \).   Therefore, their functions of the second kind \eqref{Ln} converge uniformly on closed subsets of \( \overline\C\setminus(\Delta_1 \cup \Delta_2) \) to the respective function of the second kind with respect to \( \vec\mu \). Since compactly supported measures on the real line are uniquely determined by their moments and those moments are the Laurent coefficients at infinity of the respective Markov function, it also holds that  the measures \eqref{nunE} and \eqref{up-nxi} defined with respect to \( \vec\mu_m \) converge in the weak$^*$ topology to \( \nu_{\n,E} \) and \( \upsilon_{\n,\xi} \), respectively. Notice that if $E\in E_n$ and  \( \mu^\star \) has no atom at \( E \), it holds that \( \upsilon_{\n,\xi}(\{E\}) = \pi\|\nu_{\n,E}\|^{-1} \) by \eqref{up-E}.  In particular, this is the case for each \( \vec\mu_m \). Thus, the weak$^*$ limit of the reference measures corresponding to \( \vec \mu_m \), which is obviously independent of \( \xi \), is equal to
\begin{equation}
\label{weak-ref}
\pi^{-1}S_{\n,\xi}^{-1}(x)d\upsilon_{\n,\xi}(x) -  \sum_{E:~D_\n(\xi;E)=0} \frac{S_{\n,\xi}^{-1}(E)}{\|\nu_{\n,E}\| }d\delta_E(x) =  d\omega_\n(x) + \sum_{E\in E_\n} \frac{S_{\n,\xi}^{-1}(E) \, \nu_{\n,E}(\{E\}) }{\|\nu_{\n,E}\|(\|\nu_{\n,E}\| - \nu_{\n,E}(\{E\})) } d\delta_E(x).
\end{equation}
Fix \( E\in E_\n \). Let \( k \in\{1,2\} \) be such that \( E\in \Delta_k \). Recall the definition of \( S_{\n,l,k}(x) \) in \eqref{SX} further below.  We get from the very definition of \( \nu_{\n,E} \) in \eqref{nunE}, \eqref{up-E}, and Lemmas~\ref{lem:1.4.2} and~\ref{lem:L}(5) that
\begin{equation}
\label{value-mass-p}
\frac{S_{\n,\xi}^{-1}(E) \, \nu_{\n,E}(\{E\}) }{\|\nu_{\n,E}\|(\|\nu_{\n,E}\| - \nu_{\n,E}(\{E\})) } = \frac{S_{\n,\xi}^{-1}(E)\,A_{\n+\vec e_l}^{(k)}(E)}{D_{\n,\xi}^\prime(E)S_{\n,l,k}(E)} \frac{D_{\n,\xi}^\prime(E)(A_\n^{(k)})^\prime(E)}{D_{\n,\xi}^\prime(E)L_{\n+}(E)}\mu_k(\{E\}) = \frac{Q_{\n+\vec e_l}(\{E\})}{S_{\n,l,k}(E)L_{\n+}(E)}.
\end{equation}
As the above expression is independent of \( \xi \), so is the measure~\( \omega_\n \).
\end{proof}

\section{Green's functions}
\label{s:gf}

In this section, we study functions \( G(Y,X;z) \) using equation \eqref{GYX}. The Spectral Theorem applied to  self-adjoint operator $\jox$ gives
\[
G(Y,X;z)=\langle(\jox-z)^{-1}\delta^{(X)},\delta^{(Y)}\rangle=\int \frac{d\langle P_{[X],\lambda}\delta^{(X)},\delta^{(Y)}\rangle }{\lambda-z}\,,
\]
where $\{P_{[X],\lambda}\}$ is the family of orthoprojectors associated with $\jox$.
The function $F(\lambda)= \langle P_{[X],\lambda}\delta^{(X)},\delta^{(Y)}\rangle$ has bounded variation and can be written as as difference of two non-decreasing function. Therefore, $G(Y,X;z)$ is a difference of two $\HN$ functions and the nontangential boundary values $G(Y,X;x)_{\pm}$ are defined a.e. on $\R$.  

 Let \( \mathcal T_{[X]} \) be the subtree with the root at \( X \) and \( \rho_{[X]}=\langle P_{[X],\lambda} \delta^{(X)},\delta^{(X)}\rangle \) be the spectral measure of \( \delta^{(X)} \) restricted to \( \mathcal T_{[X]} \), see \eqref{spectral-m}, where we also write \( \rho_O \) for \( \rho_{[O]} \) (we use square brackets to emphasize that \( \rho_{[X]} \) is a spectral measure of \( \delta^{(X)} \) with respect to a subtree and not the whole tree).  Then
\begin{equation}
\label{spectral-Markov}
G(X,X;z) = -\widehat\rho_{[X]}(z) \quad \text{and therefore} \quad \Im G(X,X)^+ = \pi\rho_{[X]}.
\end{equation}

Statements \eqref{spectral-Markov} and \eqref{GYX} provide  a non-trivial application of the operator theory to the theory of orthogonal polynomials. They say that the ratio of Markov functions of two ``consecutive''  linear forms \( Q_{\n+\vec e_l}(x) \) and \( Q_\n(x) \) is also a Markov function! Below, we shall verify it in a different way by providing  ``explicit'' expressions for \( \rho_{[X]} \) and more generally \( \Im G(Y,X;x)_+ \). Again, we often refer to the auxiliary lemmas proven in Section~\ref{s:ap3}.

\subsection{Function \( L_{\vec\varkappa}(z) \)}

By \eqref{GYX},  \( G(O,O;z)=-L_{\vec1}(z)/L_{\vec\varkappa}(z) \). While the behavior of the numerator \( L_{\vec1}(z) \) for smooth measures is described by Lemma~\ref{lem:Ln}, we have not yet addressed the behavior of \( L_{\vec\varkappa}(z) \). Recall that function \( L_{\vec\varkappa}(z) \) was defined in \eqref{LPO} and 
\begin{equation}
\label{LPO1}
L_{\vec\varkappa}(z)= \kappa_2 L_{\vec e_1}(z) + \kappa_1 L_{\vec e_2}(z) = \big(\varkappa_1\|\mu_1\|^{-1}\big)\widehat\mu_1(z) + \big(\varkappa_2\|\mu_2\|^{-1}\big)\widehat\mu_2(z), \,\,  \vec\varkappa = (\kappa_2,\kappa_1)\,.
\end{equation}

\begin{Lem}
\label{lem:Lkappa}
The set \( E_{\vec\varkappa}\ddd \{E:~~L_{\vec\varkappa}(E) = 0,~~E\in\R\setminus(\Delta_1\cup\Delta_2)\} \) is either empty or has exactly one element in it. It is empty when \( \vec\varkappa=\vec e_i \), \( i\in\{1,2\} \). If \( E\in E_{\vec\varkappa} \) exists, it is necessarily a simple zero of \( L_{\vec\varkappa}(x) \). If \( \vec\mu \) satisfies the assumptions of Lemma~\ref{lem:Ln}, then \( L_{\vec\varkappa}(z) \) extends continuously from \( \C_+ \) to \( \R \) and the function \( |L_{\vec\varkappa}(x)| \) is well-defined, continuous and non-vanishing on \( \R \) except for a possible single zero that belongs to \( \R\setminus(\Delta_1^\circ\cup\Delta_2^\circ) \).
\end{Lem}
\begin{proof}[Proof of Lemma~\ref{lem:Lkappa}]

The function \( L_{\vec\varkappa}(z)=\sum_{i=1}^2\varkappa_i \widehat\sigma_i(z)\), \( \sigma_i = \|\mu_i\|^{-1}\mu_i \), is analytic in \( \overline\C\setminus(\Delta_1\cup\Delta_2) \) and we are looking for its zeroes on the real line away from the intervals \( \Delta_1,\Delta_2 \).  Observe that the equation \( L_{\vec\varkappa}(x)=0 \) has no solutions on the set of interest when \( \vec\varkappa = \vec e_i \), \( i\in\{1,2\} \), since in this case it is a Markov function and Markov functions have no zeroes in the finite plane away from the convex hull of the support of the defining measure. When \( \varkappa_i>0 \), \( i\in\{1,2\} \), we have that \( L_{\vec\varkappa}(x)>0 \) for \( x\in (\beta_2,\infty) \) and \( L_{\vec\varkappa}(x)<0 \) for \( x\in(-\infty,\alpha_1) \) as one can see from \eqref{LPO1}. Since both functions \( \varkappa_i\widehat\sigma_i(x) \) are decreasing in the gap \( (\beta_1,\alpha_2) \), but one of them is negative and one is positive, there can be at most one solution there. When \( \varkappa_1\varkappa_2<0 \), there cannot be any solutions in \( (\beta_1,\alpha_2) \). To show that there is at most one solution in  \( (-\infty,\alpha_1)\cup(\beta_2,\infty) \) in this case, notice that the original equation can be rewritten as \( -(\widehat\sigma_1/\widehat\sigma_2)(x) = \varkappa_2/\varkappa_1 \). The ratio \( -(\widehat\sigma_1/\widehat\sigma_2)(z) \) is a Markov function of a measure supported on \( \Delta_1\cup\Delta_2 \). Indeed, it follows from \eqref{trace} that
\begin{equation}
\label{rhoO}
\Im(\widehat\sigma_1/\widehat\sigma_2)(x+\ic y)dx \cws \widehat\sigma_2^{-1}(x)d\big(\Im\widehat\sigma_1\big)^+(x) + \widehat\sigma_1(x)d\big(\Im\widehat\sigma_2^{-1}\big)^+(x),
\end{equation}
which is indeed a positive measure supported on \( \Delta_1\cup\Delta_2 \) since \( \widehat\sigma_2(x)<0 \), \( x\in \Delta_1 \), and \( \widehat\sigma_1(x)>0 \), \( x\in \Delta_2 \).  Markov functions are monotonically decreasing on the real line away from the support and are positive/negative to the right/left of the convex hull of the support of the defining measure. Thus, any equation of the form \( (\widehat\sigma_1/\widehat\sigma_2)(x) = \tau \neq 0 \) can have at most two solution away from \( \Delta_1\cup\Delta_2 \), one in the gap and one outside the gap, which proves the desired conclusion. 

Continuity of \( |L_{\vec\varkappa}(x)| \) when \( \vec\mu \) satisfies condition of Lemma~\ref{lem:Ln} can be shown exactly as in the proof of that lemma. Since \( \Im L_{\vec\varkappa\pm}(x) = \mp\pi \varkappa_k \sigma_k^\prime(x) \) on \( \Delta_k \) by Proposition~\ref{prop:Markov}(2-4), it vanishes at the endpoints of the intervals \( \Delta_1,\Delta_2 \). Hence, the traces \( L_{\vec\varkappa\pm}(x) \) are real at those points and the considerations of the previous paragraph can be extended from open intervals to closed ones. Since \( \Im L_{\vec\varkappa\pm}(x) \) does not vanish on \( \Delta_1^\circ\cup\Delta_2^\circ \), there cannot be any other zeroes.
\end{proof}

Notice that for Dini-continuous measures, \( |L_{\vec\varkappa}(x)| \) can vanish at an endpoint of the intervals \( \Delta_1,\Delta_2\). 

\subsection{Green's functions at \( O \)}
We already know from the Spectral Theory that \( G(O,O;z) \in \HN\).  However, we can see it directly. Recall that  \( \sigma_i = \|\mu_i\|^{-1}\mu_i \) and  define
 \begin{equation} \Xi_{\vec\mu} \ddd \int_\R td \sigma_2(t) - \int_\R td\sigma_1(t)\,.
 \label{3sd6}
 \end{equation} We have $\Xi_{\vec\mu}>0$ since it is a difference of the centers of mass of probability measures supported on disjoint intervals with \( \supp\,\sigma_1<\supp\,\sigma_2 \).
Assuming that \( \varkappa_1\neq 0 \) (the case \( \varkappa_2\neq0 \) can be treated absolutely analogously), we have that
\begin{equation}
\label{mO}
\Xi_{\vec\mu}\, G(O,O;z) = -\Xi_{\vec\mu}\,\frac{L_{\vec 1}(z)}{L_{\vec\varkappa}(z)}= -\frac{\widehat\sigma_2(z)-\widehat\sigma_1(z)}{\varkappa_2\widehat\sigma_2(z)+\varkappa_1\widehat\sigma_1(z)} = \frac1{\varkappa_1} - \frac1{\varkappa_1\varkappa_2+\varkappa_1^2(\widehat\sigma_1/\widehat\sigma_2)(z)},
\end{equation}
where  we used $\varkappa_1+\varkappa_2=1$, \eqref{LPO1}, and Lemma~\ref{lem:1.4.1}. Since \( \varkappa_1^2,\Xi_{\vec\mu}>0 \),  \( G(O,O;\cdot)\in\HN \) if and only if \( (\widehat\sigma_1/\widehat\sigma_2)\in\HN \). The claim \( (\widehat\sigma_1/\widehat\sigma_2)\in\HN \) has been shown in the proof of Lemma~\ref{lem:Lkappa} above, see \eqref{rhoO}.  

Let \( S_O(x) \) be a positive function on \( \Delta_1\cup\Delta_2 \) given by
\begin{equation}
\label{987-1}
S_O(x) \ddd \big(\Xi_{\vec\mu}\|\mu_1\|\|\mu_2\|\big)^{-1}\left(\widehat\mu_1(x)\chi_{\Delta_2}(x) - \widehat\mu_2(x)\chi_{\Delta_1}(x) \right).
\end{equation}
This function will be used to obtain a convenient formula for the generalized eigenfunction $\Psi$, introduced in the following proposition (for the general theory of eigenfunction expansion, check \cite{berezan}).

\begin{Prop}
\label{prop:rhoO}
Let \( E_{\vec\varkappa} \) be as in Lemma~\ref{lem:Lkappa}. We have \( \supp \rho_O \subseteq \Delta_1\cup\Delta_2\cup E_{\vec\varkappa} \) and
\begin{equation}
\label{987-2}
d\Im G(Y,O)^+(x) = \pi\Psi_Y(O;x)d\rho_O(x),
\end{equation}
where \( \Psi(O;E)=l(E)/L_{\vec{1}}(E) \) for \( E\in E_{\vec\varkappa} \), 
\begin{equation}
\label{987-3}
\Psi(O;x) = S_O^{-1}(x)\left( \Lambda^{(0)}(x)\frac{\varkappa_k}{\|\mu_k\|} - (-1)^k\widehat\mu_{3-k}(x)\left(\Lambda^{(2)}(x)\frac{\varkappa_1}{\|\mu_1\|} - \Lambda^{(1)}(x)\frac{\varkappa_2}{\|\mu_2\|}  \right) \right)
\end{equation}
for \( x\in\Delta_k, \, k\in \{1,2\} \), and otherwise \( \Psi(O;x)=0 \). Furthermore, it holds that
\begin{equation}
\label{987-4}
\jack \Psi(O;x) = x\Psi(O;x) \quad \text{and} \quad \delta_Y^{(O)} = \int\Psi_Y(O;x)d\rho_O(x).
\end{equation}
If \( \vec \mu \) satisfies conditions of Lemma~\ref{lem:Ln} and  
 \begin{equation}\label{4sd1}
 \ (\mu_k^\prime(x))^{-1}\in  L^p(\Delta_k) 
 \end{equation}
  for some \( p>1 \) and each \( k\in\{1,2\} \), then
\begin{equation}
\label{987-0}
d\rho_O(x) = S_O(x)|L_{\vec\varkappa}(x)|^{-2} d\mu^\star(x) + \sum_{E\in E_{\vec\varkappa}}(L_{\vec 1}/L_{\vec\varkappa}^\prime)(E) d\delta_E(x).
\end{equation}

\end{Prop}
\noindent {\bf Remark.} Assumption \eqref{4sd1} is a non-essential technical condition which we use solely to simplify the discussion of the behavior of \( \rho_O \) around a zero of \( |L_{\vec\varkappa}(x)| \) when the latter happens to be an endpoint of either \( \Delta_1 \) or \( \Delta_2 \).
\begin{proof}
The first claim follows from \eqref{mO} and the definition of \( E_{\vec\varkappa} \) in Lemma~\ref{lem:Lkappa}. Assume first that \( \vec\mu \) satisfies conditions of Lemma~\ref{lem:Ln} with the additional integrability assumption \eqref{4sd1}. We get from Lemma~\ref{lem:Lkappa} that \( |L_{\vec\varkappa}(x)| \) is continuous on the real line with at most one zero, say \( E \), that belongs to \( \R\setminus(\Delta_1^\circ\cup\Delta_2^\circ) \). Since \( -G(O,O;z) \) is a Markov function by \eqref{mO} and the explanation right after, it follows from Lemma~\ref{lem:Ln} and Proposition~\ref{prop:Markov}(2,3) that \( \rho_O \) is an absolutely continuous measure except for a possible mass point at \( E \). When \( E \) is not an endpoint of \( \Delta_1 \) or \( \Delta_2 \), we get from  Proposition~\ref{prop:Markov}(3) that \( \rho_O \) indeed has a mass point at \( E \) of mass \( (L_{\vec 1}/L_{\vec\varkappa}^\prime)(E) \). If \( E \) is an endpoint of either \( \Delta_1 \) or \( \Delta_2 \), we deduce from Proposition~\ref{prop:Markov}(3) and Lemma~\ref{lem:boot-strap} further below that \( E \) is not a mass point (this is exactly where the \(L^p\)-integrability is used). Hence, it only remains to compute the absolutely continuous part of \( \rho_O \), that is, \( \pi^{-1}\Im\big( G(O,O;x)_+\big)\), see again Proposition~\ref{prop:Markov}(2). To this end, it holds that
\[
G(O,O;x)_+ = (-1)^k \frac{\widehat\sigma_{3-k}(x)}{\Xi_{\vec\mu}} \frac{\varkappa_{3-k}\widehat\sigma_{3-k}(x) + \varkappa_k\widehat\sigma_{k-}(x)}{|L_{\vec\varkappa}(x)|^2} - (-1)^k \frac{\widehat\sigma_{k+}(x)}{\Xi_{\vec\mu}}\frac{ \varkappa_{3-k}\widehat\sigma_{3-k}(x)+\varkappa_k \widehat\sigma_{k-}(x) }{|L_{\vec\varkappa}(x)|^2}
\]
for \( x\in\Delta_k \), \( k\in\{1,2\} \), where again \( \sigma_k=\|\mu_k\|^{-1}\mu_k \) are the normalized measures. By taking the imaginary part of both sides and using $\widehat\sigma_{k-}(x)=\overline{\widehat\sigma_{k+}(x)}$ and  \( \varkappa_1+\varkappa_2=1 \), we get that
\begin{eqnarray*}
\Im\big( G(O,O;x)_+\big) &=& \frac{(-1)^k}{\Xi_{\vec\mu}} \frac{\varkappa_k\widehat\sigma_{3-k}(x)\Im\big(\widehat\sigma_{k-}(x)\big) -\varkappa_{3-k}\widehat\sigma_{3-k}(x)\Im\big(\widehat\sigma_{k+}(x)\big)}{|L_{\vec\varkappa}(x)|^2} \\
&=& \frac{(-1)^k}{\Xi_{\vec\mu}} \frac{\widehat\sigma_{3-k}(x)\Im\big(\widehat\sigma_{k-}(x)\big)}{|L_{\vec\varkappa}(x)|^2} = \pi\frac{(-1)^k \, \widehat\mu_{3-k}(x)}{\Xi_{\vec\mu}\|\mu_1\|\|\mu_2\|}\frac{\mu_k^\prime(x)}{|L_{\vec\varkappa}(x)|^2} = \pi\frac{S_O(x)\mu_k^\prime(x)}{|L_{\vec\varkappa}(x)|^2}
\end{eqnarray*}
for \( x\in\Delta_k \) by the very definition \eqref{987-1}, which finishes the proof of \eqref{987-0}.

Let us still assume that \( \vec\mu \) satisfies condition of Lemma~\ref{lem:Ln} with the additional integrability assumption  \eqref{4sd1}. Set \( g_Y(z) \) to be \( G(Y,O;z) \) when \( E_{\vec\varkappa} = \varnothing \) or \( E\in E_{\vec\varkappa} \) is an endpoint of \( \Delta_1 \) or \( \Delta_2 \) and otherwise set it to be \( G(Y,O;z)-(L_Y/L_{\vec\varkappa})(E)(E-z)^{-1} \). Then, \( -g_O(z) \) is a Markov function of an absolutely continuous measure with an \( L^p \)-density for some \( p>1 \) by Lemma~\ref{lem:boot-strap}. It follows from the last claim of Proposition~\ref{prop:Garnett} and Proposition~\ref{prop:an}(1) that both real and imaginary parts of \( g_O(z) \) satisfy \eqref{sup-u} with this \( p \). Since \( g_Y(z)=m_Y(L_Y/L_O)(z)g_O(z) \) and \( (L_Y/L_O)(z) \) extends continuously to the real line from the upper half-plane by Lemma~\ref{lem:Ln}, the imaginary part of \( g_Y(z) \) satisfies \eqref{sup-u} with the same \( p \) as well. Thus, it follows from the last claim of Proposition~\ref{prop:Garnett} that \( \Im g_Y(z) \) is a Poisson integral of an absolutely continuous measure whose density is equal to \( \Im (g_{Y+}(x) ) \). Now, we get from \eqref{LPO}, \eqref{Ln-An}, and \eqref{GYX} that
\[
G(Y,O;z) = |L_{\vec \varkappa}(z)|^{-2}\left(\Lambda_Y^{(0)}(z) - \Lambda_Y^{(1)}(z)\widehat\mu_1(z) - \Lambda_Y^{(2)}(z) \widehat \mu_2(z)\right)\left(  \big(\varkappa_1\|\mu_1\|^{-1}\big)\overline{\widehat\mu_1(z)} + \big(\varkappa_2\|\mu_2\|^{-1}\big)\overline{\widehat\mu_2(z)} \right).
\]
Since \( \Im(\widehat\mu_{k+}(x)) = -\pi\mu_k^\prime(x) \) by Proposition~\ref{prop:Markov}(2-4), it holds that
\[
\Im\big(g_{Y+}(x)\big) = \pi\mu_k^\prime(x)|L_{\vec \varkappa}(x)|^{-2}\left( \Lambda_Y^{(0)}(x)\frac{\varkappa_k}{\|\mu_k\|} - (-1)^k\widehat\mu_{3-k}(x)\left(\Lambda^{(2)}(x)\frac{\varkappa_1}{\|\mu_1\|} - \Lambda^{(1)}(x)\frac{\varkappa_2}{\|\mu_2\|}  \right) \right)
\]
for \( x\in\Delta_k \). That clearly yields \eqref{987-2} and \eqref{987-3} in the considered case.

 If the system  \( \vec\mu \) does not satisfy the assumptions of Lemma~\ref{lem:Ln} with the additional integrability assumption, approximate \( \vec\mu \) in the weak$^*$ topology by a sequence \( \{\vec\mu_m \}\) of measures that do satisfy them as it was done in the proof of Proposition~\ref{prop:RefM}. The explanation given there shows that the spectral measures and measures generated by Green's functions corresponding to \( \vec\mu_m \) will converge in the weak$^*$ sense to \( \rho_O \) and \( \Im G(Y,O)^+\) corresponding to \( \vec\mu \), respectively. This convergence will clearly preserve \eqref{987-2} and \eqref{987-3}.

The first algebraic identity of \eqref{987-4} is a direct consequence of the first two claims of Proposition~\ref{prop:2.1.4}. To prove the second identity, notice that
\[
G(Y,O;z) = \int_\R\frac{\Psi_Y(O;x)d\rho_O(x)}{x-z}
\]
by \eqref{v-recover} and \eqref{trace}. Now, since \( \Psi_Y(O;x)d\rho_O(x) \) has finite total variation, the above formula, Fubini-Tonelli Theorem, and Cauchy integral formula give that
\[
\int_\R \Psi_Y(O;x)d\rho_O(x) = \frac1{2\pi\ic}\int_\Gamma G(Y,O;z)dz = \frac{1}{2\pi \ic}\oint_{\Gamma}\Bigl((\jack-z)^{-1}\delta^{(O)}\Bigr)_Ydz = \delta_Y^{(O)},
\]
where \( \Gamma \) encircles \( \sigma(\jack)\cup\Delta_1\cup\Delta_2 \cup E_{\vec\varkappa}\) in the positive direction, the second identity is just definition \eqref{GYX}, and the last one is a part of the Spectral theorem for self-adjoint operators.
\end{proof}

\subsection{Green's functions at \( X\neq O \)}

Recall definition \eqref{commut} of the commutator of two functions with respect to a given vertex as well as definitions of functions \( \Lambda^{(k)}(x) \) in \eqref{Lambdas}. Given \( X\neq O \), set
\[
\widetilde \Psi(X;x) \ddd m_{X_{(p)}}m_X\sum_{k=1}^2\left( \big[\Lambda^{(k)}(x),\Lambda^{(0)}(x)\big]^{(X_{(p)})} +  \big[\Lambda^{(3-k)}(x),\Lambda^{(k)}(x)\big]^{(X_{(p)})} \widehat\mu_{3-k}(x)\right)\chi_{\Delta_k}(x)
\]
to be a function on \( \mathcal V \) that depends of a parameter \( x\in\Delta_1\cup\Delta_2 \). Clearly, each \( \widetilde\Psi_Y(X;x) \) extends analytically from each interval \( \Delta_k \).

\begin{Lem}
\label{lem:Commut}
Given \( X\in\mathcal V \), \( X\neq O \), it holds that \( S_X(x) \ddd \widetilde\Psi_X(X;x)>0 \) and it is continuous for \( x\in\Delta_1\cup\Delta_2 \).
\end{Lem}

We prove Lemma~\ref{lem:Commut} further below in Section~\ref{ss:An}. Recall Proposition~\ref{prop:RefM} and our notation that \( \jox \) stands for the restriction of \( \jack \) to \( \mathcal T_{[X]} \). If $E\in \R$ and $X\in \mathcal{V}$, the symbol $Q_X(\{E\})$  stands for the mass of the form $Q_{\Pi(X)}$ at point $E$, i.e., 
\[
Q_X(\{E\})\ddd A_{\Pi(X)}^{(1)}(E)\mu_1(\{E\})+A_{\Pi(X)}^{(2)}(E)\mu_2(\{E\})\,.
\]

In the next result, we explain how the spectral measure $\rho_{[X]}$ from \eqref{spectral-Markov} is related to the reference measure at the point $X_{(p)}$. We also introduce $\Psi(X;x)$, a function on $\mathcal{V}_{[X]}$ which is a  generalized eigenfunction of the operator \( \jox \).

\begin{Prop}
\label{prop:rhoX}
Let \( X\neq O \) and \( E_Y \) be the set of zeroes of the polynomial \( \Lambda_Y^{(1)}(x)\Lambda_Y^{(2)}(x) \). It holds that
\begin{equation}
\label{378-1}
d\rho_{[X]}(x) = S_X(x)d\omega_{X_{(p)}}(x) + \sum_{E\in E_{X_{(p)}}} Q_X(\{E\}) L_{X_{(p)}+}^{-1}(E)\, d\delta_E(x),
\end{equation}
where the numbers \( Q_X(\{E\}) L_{X_{(p)}+}^{-1}(E)\geq0 \) are well-defined and non-negative for each \( E\in E_{X_{(p)}} \). Moreover, it holds that
\begin{equation}
\label{378-2}
d\Im G(Y,X)^+(x) = \pi\Psi_Y(X;x)d\rho_{[X]}(x),
\end{equation}
for every \( Y\in\mathcal V_{[X]} \), where \( \Psi(X;x) = S_X^{-1}(x)\widetilde\Psi(X;x) \). Furthermore, it holds that
\begin{equation}
\label{378-3}
\jox \Psi(X;x) = x\Psi(X;x) \quad \text{and} \quad \delta_Y^{(X)} = \int\Psi_Y(X;x)d\rho_{[X]}(x).
\end{equation}
\end{Prop}
\noindent {\bf Remark.} It follows directly from definition that $\Psi$ satisfies a normalization $\Psi_X(X;x)=1$.

\begin{proof}
Assume first that \( \vec\mu \) satisfies conditions of Lemma~\ref{lem:Ln}. Recall that the traces \( \widehat\mu_{k\pm}(x) \) are continuous on the real line and are complex conjugates of each other, see Proposition~\ref{prop:an}(2). It follows from \eqref{GYX} and \eqref{Ln-An} that
\begin{multline*}
\Im(G(Y,X;z))  =  -m_X m_{X_{(p)}}|L_{X_{(p)}}(z)|^{-2} \Im\left( \left(\Lambda_Y^{(1)}(z)\widehat\mu_1(z) + \Lambda_Y^{(2)}(z)\widehat\mu_2(z)-\Lambda_Y^{(0)}(z)\right)\right. \times \\ \left.\overline{\left(\Lambda_{X_{(p)}}^{(1)}(z)\widehat\mu_1(z) + \Lambda_{X_{(p)}}^{(2)}(z)\widehat\mu_2(z) - \Lambda_{X_{(p)}}^{(0)}(z)\right)} \right).
\end{multline*}
Since the first kind MOPs have real coefficients, a straightforward algebraic computation and Lemma~\ref{lem:Ln} imply that \( \Im(G(Y,X;z)) \) has continuous traces on the real line and
\begin{equation}
\label{GYX+}
\Im(G(Y,X;x)_+)  = -\widetilde\Psi_Y(X;x) |L_{X_{(p)}}(x)|^{-2} \Im\big(\widehat\mu_{k+}(x)\big), \quad x\in \Delta_k, ~~k\in\{1,2\}.
\end{equation}
In particular, we get from Proposition~\ref{prop:Markov}(2-4), Lemmas~\ref{lem:Ln} and~\ref{lem:Commut} that \( \Im(G(X,X;z)) \) extends continuously to the real line where it has a continuous and non-negative trace. Thus, it follows from the maximum principle for harmonic functions that \( \Im(G(X,X;\cdot)) \in \HN \), the fact that we already know from the general Spectral Theory. Since \( -G(X,X;z) \) is holomorphic at infinity, it is indeed a Markov function. Formula \eqref{378-1} now follows from Propositions~\ref{prop:Markov}(2,3) and~\ref{prop:RefM} since \( Q_X(\{E\})=0 \) for any \( E \) by absolute continuity of \( \mu^\star \). Since \( \rho_{[X]} \) is absolutely continuous with continuous density, we get from the last claim of Proposition~\ref{prop:Garnett} and Proposition~\ref{prop:an}(1) that both real and imaginary parts of \( G(X,X;z) \) satisfy \eqref{sup-u} for any \( p>1 \). Since \( G(Y,X;z) = m_Y(L_Y/L_X)(z)G(X,X;z) \) and \( (L_Y/L_X)(z) \) extends continuously to the real line by Lemma~\ref{lem:Ln}, the imaginary part of \( G(Y,X;z) \) also satisfies \eqref{sup-u} for any \( p>1 \). Thus, \( \Im(G(Y,X;z)) \) is a Poisson integral of an absolutely continuous measure with density given by \( \Im(G(Y,X;x)_+) \), which, together with \eqref{GYX+}, proves \eqref{378-2} in the considered case.

If the system  \( \vec\mu \) does not satisfy assumptions of Lemma~\ref{lem:Ln}, approximate \( \vec\mu \) in the weak$^*$ topology by systems \( \vec\mu_m \) that do satisfy these assumptions as it was done in the proof of Proposition~\ref{prop:RefM}. The explanation given there shows that the spectral measures corresponding to \( \vec\mu_m \) converge in the weak$^*$ sense to \( \rho_{[X]} \), the spectral measure corresponding to \( \vec\mu \). On the other hand, the right-hand sides of \eqref{378-1} corresponding to \( \vec\mu_m \) will converge weak$^*$ to \( S_X(x) \) times the measure in \eqref{weak-ref}.  Formula \eqref{378-1} now follows from \eqref{value-mass-p} and from the identity \( S_X(x)=S_{\n,l,k}(x) \) for \( x\in\Delta_k \) which holds by the definition of \( S_{\n,l,k}(x) \) in \eqref{SX}, where \( \n=\Pi(X_{(p)}) \) and \( l=\iota_X \).  As \( \Im G(Y,X)^+\) is the weak$^*$ limit of the corresponding measures with respect to \( \vec\mu_m \), the validity of \eqref{378-2} follows as well.

The first algebraic identity of \eqref{378-3} is a direct consequence of the third claim of Proposition~\ref{prop:2.1.4}. The second one can be justified exactly as in Proposition~\ref{prop:rhoO}.
\end{proof}

\section{Cyclic subspaces}
\label{s:cs}

In this section we derive an orthogonal decomposition of \( \ell^2(\mathcal V) \) into a direct sum of cyclic subspaces.

\subsection{Trivial cyclic subspaces}

Let \( X\in\mathcal V \) and \( \alpha(x) \) be a polynomial. Formulas \eqref{987-4} and \eqref{378-3} immediately allow us to conclude that
\begin{equation}
\label{mly-10}
\alpha(\jox)\delta^{(X)} = \int \alpha(\jox)\Psi(X;x)d\rho_{[X]}(x) = \int \alpha(x)\Psi(X;x) d\rho_{[X]}(x) \ddd \widehat\alpha \in \ell^2(\mathcal V_{[X]}),
\end{equation}
where the last conclusion trivially holds as \( \alpha(\jox)\delta^{(X)} \) is compactly supported in this case. Of course, \eqref{mly-10} can be further extended to continuous functions on \( \Delta_1\cup\Delta_2 \) using the Spectral Theorem. Namely, let $\{P_{[X],\lambda}\}$ be the orthogonal spectral decomposition for $\jox$. Then, it holds that
\[
\alpha(\jox)\delta^{(X)} \ddd \left(\int \alpha(\lambda)dP_{[X],\lambda}\right)\delta^{(X)}\in \ell^2(\mathcal V_{[X]}).
\]
In fact, we can say more. Let \( \mathfrak C_{[X]}^{(X)} \) be the cyclic subspace of \( \ell^2(\mathcal V_{[X]}) \) generated by \( \delta^{(X)} \), that is,
\[
\mathfrak C_{[X]}^{(X)} \ddd \, \overline{ \mathrm{span}\left\{\jox^n\delta^{(X)}:~n\in\Z_+\right\} } = \overline{ \left\{\alpha(\jox)\delta^{(X)}:~\alpha~\text{ is a polynomial}\right\} }.
\]
The next result is an analog of Proposition~\ref{lem1}, where $\Psi$ plays the role of orthogonal polynomials.

\begin{Prop}
\label{prop:triv-c}
Fix \( X\in\mathcal V \). The map
\begin{equation}
\label{864-0}
\alpha(x) \mapsto \widehat \alpha = \big\{\widehat\alpha_Y\big\}_{Y\in\mathcal V_{[X]}}, \quad \widehat \alpha_Y \ddd \int \alpha(x)\Psi_Y(X;x)d\rho_{[X]}(x),
\end{equation}
is a unitary map from $L^2(\rho_{[X]})$ onto $\mathfrak C_{[X]}^{(X)}$. In particular,  it holds that
\begin{equation}
\label{864-1}
\|\alpha\|_{L^2(\rho_{[X]})}^2 = \|\widehat\alpha\|_{\ell^2(\mathcal V_{[X]})}^2 \quad \text{and} \quad \mathfrak C_{[X]}^{(X)} = \left\{\widehat\alpha:~\alpha\in L^2(\rho_{[X]})\right\}.
\end{equation}
Thus, the formula
\begin{equation}
\label{864-3}
\alpha(\jox)\delta^{(X)} \ddd \widehat\alpha = \int \alpha(x)\Psi(X;x)d\rho_{[X]}(x) 
\end{equation}
extends the definition of \(\alpha(\jox)\delta^{(X)} \) from continuous functions \( \alpha(x) \) to those in \( L^2(\rho_{[X]}) \). We also have that
\begin{equation}
\label{864-2}
x\alpha(x) \mapsto \jox\widehat\alpha, \quad \alpha\in L^2(\rho_{[X]}).
\end{equation}
\end{Prop}
\begin{proof} 
The following argument is standard and we adduce it solely for completeness. Let \( \alpha(x) \) be a continuous function on \( \Delta_1\cup\Delta_2 \). It follows from the Spectral Theorem that
\begin{eqnarray}\nonumber
\|\alpha(\jox)\delta^{(X)}\|^2_{\ell^2(\mathcal V_{[X]})} &=& \big\langle \alpha(\jox)\delta^{(X)},\alpha(\jox)\delta^{(X)} \big\rangle = \big\langle |\alpha(\jox)|^2\delta^{(X)},\delta^{(X)} \big\rangle \\
& = & \int |\alpha(\lambda)|^2d\langle P_{[X],\lambda} \delta^{(X)},\delta^{(X)}\rangle=\int|\alpha(x)|^2d\rho_{[X]}(x) = \|\alpha\|_{L^2(\rho_{[X]})}^2\,\label{5sd1}
\end{eqnarray}
since $\rho_{[X]}$ is the spectral measure for $\delta^{(X)}$ in \( \ell^2(\mathcal V_{[X]}) \). Take any $\alpha\in L^2(\rho_{[X]})$ and approximate it in $L^2(\rho_{[X]})$ by a sequence $\big\{\alpha^{(n)}\big\}$ of polynomials. Recall that each \( \alpha^{(n)}(\jox) \) is compactly supported and therefore is in \( \ell^2(\mathcal V_{[X]}) \). Because \( \Psi_Y(X;x) \) is continuous on \( \Delta_1\cup\Delta_2 \), it holds that \( \widehat\alpha_Y=\lim_{n\to\infty} \widehat\alpha^{(n)}_Y \) for every $Y$. Thus, 
\[
\sum_{|\Pi(Y)|<N}|\widehat\alpha_Y|^2\le \int|\alpha|^2d\rho_{[X]}\, \, \text{for any } N\in \mathbb{N} \quad \Rightarrow \quad \|\widehat\alpha\|^2_{\ell^2(\mathcal V_{[X]})}\le \int|\alpha|^2d\rho_{[X]}\,
\]
and therefore \( \big\{\widehat\alpha:~\alpha\in L^2(\rho_{[X]})\big\} \subseteq \mathfrak C_{[X]}^{(X)} \). Furthermore, let \(\Phi\in\mathfrak C_{[X]}^{(X)} \) and \( \widehat\alpha^{(n)}\to\Phi \) as \( n\to\infty \) in \( \ell^2(\mathcal V_{[X]}) \) for some sequence $\big\{\alpha^{(n)}\big\}$ of polynomials. By \eqref{5sd1}, we have $\sup_n \|\alpha^{(n)}\|_{L^2(\rho_{[X]})}<\infty$
and, according to Banach-Alaoglu, there exists \( \phi\in L^2(\rho_{[X]}) \) such that \( \alpha^{(n_k)}\to\phi \) weakly in \( L^2(\rho_{[X]}) \) as \( k\to\infty \). Therefore, evaluating at each $Y\in \mathcal{V}$, we get
\[
\widehat\phi_Y = \int\phi(x)\Psi_Y(X;x)d\rho_{[X]}(x) \leftarrow \int \alpha^{(n_k)}(x)\Psi_Y(X;x)d\rho_{[X]}(x) = \widehat\alpha^{(n_k)}_Y \rightarrow \Phi_Y
\]
as \( k\to\infty \). Hence, \( \big\{\widehat\alpha:~\alpha\in L^2(\rho_{[X]})\big\} =\mathfrak C_{[X]}^{(X)} \). That is, the map $\alpha \mapsto \widehat\alpha$ is onto as well as isometric on the dense subset so it is isometric everywhere. Thus, the considered map  $\alpha \mapsto \widehat\alpha$ is actually unitary, which finishes the proof of \eqref{864-1}. Finally, one can readily see that
\[
\jox\widehat\alpha = \jox\int \Psi(X;x)\alpha(x)d\rho_{[X]}(x) = \int \jox\Psi(X;x)\alpha(x)d\rho_{[X]}(x) =  \int x\Psi(X;x)\alpha(x)d\rho_{[X]}(x)
\]
by \eqref{378-3}, which shows \eqref{864-2}.
\end{proof}

\subsection{Non-trivial cyclic subspaces}

Fix \( X\in\mathcal V \) and let \( X_i=X_{(ch),i} \), \( i\in\{1,2\} \). 
Put
\begin{equation}
\label{hatrhoX}
\widetilde\rho_X \ddd \omega_X + \sum_{E\in E_X}\mu^\star(\{E\})\delta_E,
\end{equation}
where \( \omega_X \) is the reference measure from Proposition~\ref{prop:RefM}, \( E_X \) is the set of zeroes of \( \Lambda_X^{(1)}(x)\Lambda_X^{(2)}(x) \), and \( \mu^\star \) is the concatenated measure from \eqref{mustar}. It readily follows from \eqref{378-1} that
\begin{equation}
\label{nuXi}
d\rho_{[X_i]}(x) = \nu_{X_i}(x)d\widetilde\rho_X(x),
\end{equation}
where \( \nu_{X_i}(x) = S_{X_i}(x) \) for \( x\in(\Delta_1\cup\Delta_2)\setminus E_X \) and \( \nu_{X_i}(E) = A_{X_i}^{(k)}(E) L_{X+}^{-1}(E) \) for \( E\in E_X\cap\Delta_k \). Most importantly for us there exists \( c_X>1 \) such that
\begin{equation}\label{emph1}
c_X^{-1}\leq \nu_{X_i}(x)  \leq c_X, \quad x\in\Delta_1\cup\Delta_2,
\end{equation}
 according to Lemmas~\ref{lem:Commut} and~\ref{lem:interlacing} (it is also continuous on \( (\Delta_1\cup\Delta_2) \setminus E_X\)). Let \( \Psi(X_i;x)\) be the generalized eigenfunction from Proposition~\ref{prop:rhoX}. Recall that \( W_{X_i}>0 \). Let
\begin{equation}
\label{hatPsi}
\widehat \Psi_Y(X;x) \ddd (-1)^iW_{X_i}^{-1/2}\Psi_Y(X_i;x),~~Y\in\mathcal V_{[X_i]}, \quad \text{and} \quad \widehat \Psi_Y(X;x) \ddd 0, ~~\text{otherwise}.
\end{equation}
We stress that \( \widehat \Psi(X;x) \) is a function on $\mathcal V$ that is supported by $\mathcal V_{[X]}$ with value zero at \( X \) itself. Define
\begin{equation}
\label{Chat}
\widehat{\mathfrak C}^{(X)} \ddd \left\{\int \alpha(x)\widehat\Psi(X;x)d\widetilde\rho_X(x):~~\alpha\in L^2(\widetilde\rho_X)\right\}.
\end{equation}
Let \( \chi_i \) be the restriction operator that sends \( f\in \widehat{\mathfrak C}^{(X)}\) to its restriction to \( \mathcal V_{[X_i]} \), \( i\in\{1,2\} \). It readily follows from Proposition \ref{prop:triv-c},  \eqref{nuXi}, and property \eqref{emph1} that
\[
\left\{\chi_if:~~f\in\widehat{\mathfrak C}^{(X)} \right\} = \mathfrak C_{[X_i]}^{(X_i)} .
\]
Observe that  \( \chi_i:\widehat{\mathfrak C}^{(X)} \to \mathfrak C_{[X_i]}^{(X_i)} \) is a bijection and the composition \( \chi_2\circ\chi_1^{-1} \) is a bijection between \( \mathfrak C_{[X_1]}^{(X_1)} \) and \( \mathfrak C_{[X_2]}^{(X_2)} \). Altogether, we can say that
\begin{equation}
\label{Chat-char}
f\in \widehat{\mathfrak C}^{(X)} \quad \Leftrightarrow \quad \supp\,f\subseteq\mathcal V_{[X_1]}\cup\mathcal V_{[X_2]}, \quad f_i \in \mathfrak C_{[X_i]}^{(X_i)},~~i\in\{1,2\}, \quad \text{and} \quad \chi_1^{-1}f_1 = \chi_2^{-1}f_2,
\end{equation}
where \( f_i \) is the restriction of \( f \) to \( \mathcal V_{[X_i]} \).

\begin{Prop}
\label{prop:Chat}
Fix \( X\in\mathcal V \). The function \( \widehat\Psi(X;x) \) is a generalized eigenfunction of \( \jack \), that is, it holds that
\begin{equation}
\label{093-0}
\jack \widehat\Psi(X;x) = x\widehat\Psi(X;x).
\end{equation}
Moreover, let the function \( g_i^{(X)}\in\widehat{\mathfrak C}^{(X)} \) be given by
\begin{equation}
\label{093-2}
g_i^{(X)} \ddd \int\alpha(X_i;x)\widehat\Psi(X;x)d\widetilde\rho_X(x), \quad \alpha(X_i;x) \ddd (-1)^i W_{X_i}^{1/2}\nu_{X_i}(x).
\end{equation}
Then, it holds that \( \chi_ig_i^{(X)} = \chi_i\delta^{(X_i)} \), \( i\in\{1,2\} \), and
\begin{equation}
\label{093-1}
\widehat{\mathfrak C}^{(X)} = \overline{ \mathrm{span}\left\{\jack^n g_i^{(X)}:~n\in\Z_+\right\} }.
\end{equation}
That is, each $g_i^{(X)}$ is a generator of the cyclic subspace $\widehat{\mathfrak C}^{(X)}$. In particular, the formula
\begin{equation}
\label{093-3}
\alpha(\jack) g_i^{(X)} \ddd \int \alpha(x)\alpha(X_i;x)\widehat\Psi(X;x)d\widetilde\rho_X(x) 
\end{equation}
extends the definition of \( \alpha(\jack) g_i^{(X)} \) from continuous functions \( \alpha(x) \) to those in \( L^2(\widetilde\rho_X) \). Furthermore, it holds that
\begin{equation}
\label{093-4}
d\rho_{X,i}(x) = \sum_{k=1}^2 \frac{W_{X_i}}{W_{X_k}}\frac{\nu_{X_i}^2(x)}{\nu_{X_k}(x)} d\widetilde\rho_X(x),
\end{equation}
where \( \rho_{X,i} = \rho_{g_i^{(X)}} \) is the spectral measure of \( g_i^{(X)} \) with respect to the operator $\jack$.
\end{Prop}

\begin{proof}
If \( Y\not\in\mathcal V_{[X]} \), it clearly holds that \( (\jack \widehat\Psi(X;x))_Y = 0 = x\widehat\Psi_Y(X;x) \). Further, we get straight from \eqref{hatPsi} that
\[
(\jack \widehat\Psi(X;x))_X = W_{X_1}^{1/2}\widehat\Psi_{X_1}(X;x) + W_{X_2}^{1/2}\widehat\Psi_{X_2}(X;x) = - \Psi_{X_1}(X_1;x) + \Psi_{X_2}(X_2;x) = 0 = x\widehat\Psi_X(X;x) 
\]
since \( \Psi_{X_i}(X_i;x) = 1 \) according to their definition, see remark after Proposition~\ref{prop:rhoX}. Moreover, if \( Y\in\mathcal V_{X_i} \), then we get from \eqref{hatPsi} and \eqref{378-3} that
\[
(\jack \widehat\Psi(X;x))_Y = (-1)^iW_{X_i}^{-1/2}(\jo_{[X_i]} \Psi(X_i;x))_Y = (-1)^iW_{X_i}^{-1/2}x\Psi_Y(X_i;x) = x\widehat\Psi_Y(X;x),
\]
which proves \eqref{093-0}. Further, it holds that \( \chi_ig_i^{(X)} = \chi_i\delta^{(X_i)} \) since
\[
\big(g_i^{(X)}\big)_Y = (-1)^iW_{X_i}^{-1/2}\int\alpha(X_i;x)\Psi_Y(X_i;x)d\widetilde\rho_X(x) = \int\Psi_Y(X_i;x)d\rho_{[X_i]}(x) = \delta_Y^{(X_i)}, \quad Y\in\mathcal V_{[X_i]},
\]
where we used \eqref{hatPsi}, \eqref{nuXi}, and \eqref{378-3}. Now, according to \eqref{Chat-char}, to prove \eqref{093-1} it is enough to show that the closure of the span of \( \chi_i\jack^n g_i^{(X)} \) is equal to \( \mathfrak C_{[X_i]}^{(X_i)} \). As \( \chi_i \) and \( \jack \) commute by \eqref{Chat} and \eqref{093-0} (or, put differently, 
$
\chi_i \jack^n g_i^{(X)}=\mathcal{J}_{[X_i]}^n (\chi_i \delta^{(X_i)})\,)
$
the latter claim follows. Formula \eqref{093-3} can be obtained through approximation by polynomials exactly as an analogous formula of Proposition~\ref{prop:triv-c} was proved. Finally, to get \eqref{093-4}, observe that
\begin{eqnarray*}
\left\langle (\jack-z)^{-1}g_i^{(X)},g_i^{(X)} \right\rangle &=& \sum_{k=1}^2 \frac1{W_k}   \left\langle  \int 
\frac{\alpha(X_i;x)}{x-z} \Psi(X_k;x)d\widetilde\rho_X(x), \int \alpha(X_i;x) \Psi(X_k;x)d\widetilde\rho_X(x) \right\rangle \\
&=& \sum_{k=1}^2 \frac{W_i}{W_k}   \left\langle \int \frac{\nu_{X_i}(x)}{(x-z)\nu_{X_k}(x)} \Psi(X_k;x)d\rho_{[X_k]}(x), \int \frac{\nu_{X_i}(x)}{\nu_{X_k}(x)}\Psi(X_k;x)d\rho_{[X_k]}(x) \right\rangle
\end{eqnarray*}
where we used \eqref{nuXi} and \eqref{093-3}. Now it follows from \eqref{spectral-m}, \eqref{378-3}, \eqref{864-3}, and \eqref{nuXi} that
\[
\int\frac{d\rho_{X,i}(x)}{x-z} = \sum_{k=1}^2 \frac{W_i}{W_k}\int \frac{\nu_{X_i}^2(x)}{\nu_{X_k}(x)} \frac{d\widetilde\rho_X(x)}{x-z}.
\]
Since Markov functions are uniquely determined by their defining measures, \eqref{093-4} follows.
\end{proof}

\subsection{Decomposition into an orthogonal sum of cyclic subspaces}

In this subsection, we will prove a theorem that, in the view of Theorem~\ref{thm:spectrum}, constitutes the central result of this paper.

\begin{Thm} 
\label{mt1} 
The Hilbert space \( \ell^2(\mathcal V) \) decomposes into an orthogonal sum of cyclic subspaces of \( \jack \) as follows:
\begin{equation}
\label{g2}
\ell^2(\cal{V})=\mathfrak{C}^{(O)}\oplus \mathcal{L}, \quad \mathcal{L}= \oplus_{Z\in \cal{V}} \widehat{\mathfrak{C}}^{(Z)}\,.
\end{equation}
\end{Thm}
\begin{proof}
First, we need to show that the subspaces on the right-hand side of \eqref{g2} are orthogonal to each other. Recall that \( \widehat{\mathfrak C}^{(Y)} \) is supported by the subtree \( \mathcal T_{[Y]} \). Let \( Z,X\in\mathcal V \), \( Z\neq X \).  If the subtrees \( \mathcal T_{[X]} \) and \( \mathcal T_{[Z]} \) are disjoint, the subspaces \( \widehat{\mathfrak C}^{(X)} \) and \( \widehat{\mathfrak C}^{(Z)} \) are naturally orthogonal. If they are not disjoint, one is a subtree of another. Assume for definiteness that \( \mathcal T_{[Z]} \) is a (proper) subtree of \( \mathcal T_{[X]} \). That is, \( Z \) is a descendant of \( X \). Let \( i\in\{1,2\} \) be such that \( Z \) is equal to or is a descendant of \( X_{(ch),i} \). Let \( \alpha(x) \) be a polynomial and \( f\in\widehat{\mathfrak C}^{(Z)} \). Then
\[
\left\langle f,\alpha(\jack) g_i^{(X)}\right\rangle = \left\langle \overline\alpha(\jack)f, g_i^{(X)}\right\rangle = \left\langle \overline\alpha(\jack)f, \delta^{(X_i)}\right\rangle = \big( \overline\alpha(\jack)f \big)_{X_i} = 0
\]
since \( \overline\alpha(\jack)f\in\widehat{\mathfrak C}^{(Z)} \) and \( X_i \) does not belong to the support of any \( h\in\widehat{\mathfrak C}^{(Z)} \). Because functions \( \alpha(\jack) g_i^{(X)} \)  are dense in  \( \widehat{\mathfrak C}^{(X)} \) by \eqref{093-1}, we get that $\widehat{\mathfrak C}^{(X)} \perp \widehat{\mathfrak C}^{(Z)}$ as claimed. When the subspace \( \widehat{\mathfrak C}^{(X)} \) is replaced by \( \mathfrak C^{(O)} \),  the proof remains absolutely the same except that we need to consider functions \( \alpha(\jack) \delta^{(O)} \) instead of \( \alpha(\jack) g_i^{(X)} \).

Since all cyclic subspaces are orthogonal to each other, to prove the theorem, it is enough to show that finite sums of the above cyclic subspaces contain all the functions with compact support. As the latter are linear combinations of delta functions, it is sufficient to show that all delta functions belong to such finite sums. Trivially, it holds that \( \delta^{(O)}\in\mathfrak C^{(O)} \). By going down the tree \( \mathcal T \), we shall inductively show  that
\[
\delta^{(X)} \in \mathfrak{C}^{(O)}\oplus \mathcal{L}_X, \quad \mathcal{L}_X= \oplus_{Y\in\mathrm{path}(X_{(p)},O)} \widehat{\mathfrak{C}}^{(Y)}\,,
\]
for any \( X\in\mathcal V \), \( X\neq O \), where \( \mathrm{path}(X_{(p)},O) \) is the same as \eqref{gev1}. Take such \( X \) and assume the claim is true for \( X_{(p)} \) and \( X_{(g)} \), where \( X_{(g)} \) is parent of \( X_{(p)} \). Let \( Z \) be the sibling of \( X \). It follows from \eqref{093-2} that
\[
\left(g_{\iota_Z}^{(X_{(p)})}\right)_Z=1 \quad \text{and} \quad \left(g_{\iota_Z}^{(X_{(p)})}\right)_X = (-1)^{\iota_X}W_X^{-1/2}\int\alpha(Z;x)d\widetilde\rho_X(x) = -(W_Z/W_X)^{1/2}.
\]
We further get from the very definition of \( \jack \) in \eqref{jack2} that
\[
\big(\jack\delta^{(X_{(p)})} \big)_X = W_X^{1/2}, \quad \big(\jack\delta^{(X_{(p)})} \big)_Z = W_Z^{1/2}, \quad \big(\jack\delta^{(X_{(p)})} \big)_{X_{(p)}} = V_{X_{(p)}}, \quad \text{and} \quad \big(\jack\delta^{(X_{(p)})} \big)_{X_{(g)}} = U_{X_{(p)}},
\]
where \( U_{X_{(p)}} =0 \) if \( X_{(p)} = O \) and \( U_{X_{(p)}} = W_{X_{(p)}}^{1/2} \) otherwise (all other values of \( \jack\delta^{(X_{(p)})} \) are equal to zero). Extend \( \Psi(X;x) \) from \( \mathcal V_{[X]} \) to the whole set \( \mathcal V \) by zero. Then
\begin{multline}
\label{beta-hat}
\sqrt{W_XW_Z} \left[W_Z^{-1/2}\left( \jack\delta^{(X_{(p)})} - V_{X_{(p)}}\delta^{(X_{(p)})} - U_{X_{(p)}}\delta^{(X_{(g)})} \right) - g_{\iota_Z}^{(X_{(p)})}\right] \\ = W_X\delta^{(X)} + W_Z \int \frac{\nu_Z(x)}{\nu_X(x)}\Psi(X;x)d\rho_{[X]}(x) = \int\beta(X;x) \Psi(X;x)d\rho_{[X]}(x) \ddd \widehat\beta(X),
\end{multline}
where we used \eqref{378-3} for the last equality. By \eqref{emph1}, the function  \( \beta(X;x) \ddd W_X+ W_Z(\nu_Z/\nu_X)(x) \) is  strictly positive on the support of \( \rho_{[X]} \). Observe that \( \widehat\beta(X) \) is supported on \( \mathcal V_{[X]} \) and has value \( W_X+W_Z>0 \) at \( X \). It follows from the properties of \( \beta(X;x) \) that
\[
\overline{\big\{\alpha(x)\beta(X;x):~~\alpha\text{ is a polynomial}\big\}} = L^2(\rho_{[X]}),
\]
where the closure is taken in \( L^2(\rho_{[X]}) \)-norm. Thus, there exists a sequence of polynomials \( \{\alpha^{(n)}(x)\} \) such that \(  \alpha^{(n)}(x)\beta(X;x) \to 1 \) as \( n\to\infty \) in  \( L^2(\rho_{[X]}) \)-norm and therefore
\begin{equation}
\label{alphan-deltaX}
\alpha^{(n)}(\jox)\widehat\beta(X) = \int\alpha^{(n)}(x)\beta(X;x)\Psi(X;x)d\rho_{[X]}(x) \to \delta^{(X)} 
\end{equation}
as \( n\to\infty \) in \( \ell^2(\mathcal V_{[X]})\) by \eqref{864-2} and since \( \widehat\beta(X)\in\mathfrak C_{[X]}^{(X)} \), where we extend \( \alpha^{(n)}(\jox)\widehat\beta(X) \) from \( \mathcal V_{[X]} \) to \( \mathcal V \) by zero. On the other hand, it follows from \eqref{beta-hat} that
\[
\jox\widehat\beta(X) = \jack\widehat\beta(X) - W_X^{1/2}\widehat \beta_X(X)\delta^{(X_{(p)})} = \gamma(\jack)\delta^{(X_{(p)})} + \jack\left(c_1\delta^{(X_{(g)})} + c_2 g_{\iota_Z}^{(X_{(p)})}\right) \in \mathfrak C^{(O)}\oplus \mathcal L_X,
\]
where \( \gamma(x) \) is a certain quadratic polynomial and \( c_1,c_2 \) are certain constants (all explicitly expressible using \eqref{beta-hat}) and the last conclusion follows from the inductive hypothesis and the nature of cyclic subspaces, see \eqref{093-1}. By iterating the above relation we get that
\[
\alpha^{(n)}(\jox)\widehat\beta(X)\in \mathfrak C^{(O)}\oplus \mathcal L_X \quad \Rightarrow \quad \delta^{(X)}\in \mathfrak C^{(O)}\oplus \mathcal L_X,
\]
where the last conclusion is a consequence of \eqref{alphan-deltaX} and \( \mathfrak C^{(O)}\oplus \mathcal L_X \) being closed. That finishes the proof of the theorem.
\end{proof}

\section{Spectral analysis} 
\label{s:sa}

In this section, we will apply Theorem \ref{mt1} to the analysis of the spectral type of $\jack$. 

\begin{Thm}
\label{c_1}
Let  $E_{\vec\varkappa}$ be as in Lemma~\ref{lem:Lkappa}. It holds that
\begin{equation}
\label{s-jack}
\sigma(\jack)\subseteq \Delta_1\cup\Delta_2\cup E_{\vec\varkappa}.
\end{equation}
Furthermore, if \( \supp\mu_k=\Delta_k \) for each \( k\in\{1,2\} \), then inclusion in \eqref{s-jack} becomes equality.
\end{Thm}
\begin{proof}
It follows from Theorems~\ref{thm:spectrum} and~\ref{mt1}, and Proposition~\ref{prop:Chat} that
\[
\sigma(\jack) = \supp\,\rho_O \cup \bigcup_{Z\in\mathcal V} \supp\,\rho_{Z,1}
\]
where \( \rho_{Z,1} \) is the spectral measure of \( g_1^{(Z)} \). As stated in Proposition~\ref{prop:rhoO}, we have that
\[
\supp\, \rho_O\subseteq \Delta_1\cup\Delta_2\cup E_{\vec\varkappa},
\]
where inclusion becomes equality when \( \supp\mu_k=\Delta_k \) for each \( k\in\{1,2\} \) as can be seen from \eqref{rhoO} and \eqref{mO}. We further get from \eqref{093-4}  that \( \rho_{Z,1} \) is absolutely continuous with respect to \( \widehat\rho_Z \). Since \( \supp\,\widehat\rho_Z\subseteq\Delta_1\cup\Delta_2 \) by \eqref{hatrhoX}, Proposition~\ref{prop:RefM}, and Lemma~\ref{lem:interlacing}, the claim of the theorem follows.
\end{proof}

This result complements  characterization of the essential spectrum of $\jack$ obtained in the recent preprint \cite{ady2} where all right limits of \( \jack \) for \( \vec\kappa=\vec e_i \) were computed for the case where the measures $\mu_1,\mu_2$ are absolutely continuous with analytic and non-vanishing densities.

As the following example shows, in general, $\sigma(\jack)\neq  \supp \mu_1\cup \supp\mu_2$ even when \( E_{\vec\varkappa}=\varnothing \). Thus, equality \eqref{sd_21} does not hold for the case of multiple orthogonality. 

\begin{Exa} 
Consider any probability measures $\mu_1,\mu_2$ for which \( \supp\,\mu_1=[-1,0] \) and \( \supp \, \mu_2 = \{1,2\}\cup [3,4] \), i.e., $1$ and $2$ are  isolated atoms of \( \mu_2 \).  Clearly, $\Delta_1=[-1,0]$ and \( \Delta_2 = [1,4] \). Consider \( \jo_{\vec e_1} \). Formulae \eqref{spectral-Markov} and \eqref{mO} become
\[
\widehat\rho_O(z) = \frac1{\Xi_{\vec\mu}}\frac{\widehat \mu_2(z) - \widehat \mu_1(z)}{\widehat \mu_2(z)}.
\]
Since \( \widehat \mu_2(z) \) necessarily has a zero on \( (1,2) \), \( \rho_O \) has a point mass there and therefore its support is clearly not a subset of \( \supp\,\mu_1\cup\supp\,\mu_2\).
\end{Exa}

It is standard in the multidimensional scattering theory to deal with operators that have purely absolutely continuous spectrum (see \cite{rs1} for basics of Spectral Theory). In the next theorem, we provide simple conditions for  $\jack$  to have such a spectrum. 

\begin{Thm}
\label{c_3}
Suppose that $d\mu_k(x)=\mu_k^\prime(x)dx$  and $(\mu_k^\prime)^{-1} \in L^\infty(\Delta_k)$ for each $k\in \{1,2\}$. Then, the spectrum of $\jo_{\vec e_i}$ is purely absolutely continuous for each \( i\in\{1,2\} \).
\end{Thm}
\begin{proof}
We need to show that the spectral measures \( \rho_O \) and \( \{\rho_{Z,1}\} \), \( Z\in\mathcal V \), are all absolutely continuous. It follows from \eqref{093-4} that \( \rho_{Z,1} \) is absolutely continuous with respect to \( \widetilde\rho_Z \). Since measures \( \mu_1,\mu_2 \) have no mass points, we get from \eqref{hatrhoX} that \( \widetilde\rho_Z \) is equal to the reference measure \( \omega_\n \), \( \n=\Pi(Z) \). To show that the latter has no singular part, it is enough to prove that
\[
\limsup_{y\to0^+}\Im\left(\big(D_{\n,\xi}(x+\ic y)L_\n(x+\ic y)\big)^{-1}\right) < \infty \quad \text{for every} \quad x\in (\Delta_1\cup\Delta_2)\setminus E_\n,
\]
according to \eqref{DnLnOmn} and Proposition~\ref{prop:Markov}(3), where \( E_\n\) is the set of zeroes of \( A_\n^{(1)}(z)A_\n^{(2)}(z) \) and \( D_{\n,\xi}(z) \) is given by \eqref{ur3} with \( \xi\in(\beta_1,\alpha_2) \). It clearly holds that
\[
\lim_{y\to0^+}\Im\left(\big(D_{\n,\xi}(x+\ic y)L_\n(x+\ic y)\big)^{-1}\right) \leq \lim_{y\to0^+} \left(-\Im\big(D_{\n,\xi}(x+\ic y)L_\n(x+\ic y)\big)\right) ^{-1}.
\]
Fix \( k\in\{1,2\} \) and a closed subinterval \( \Delta \) of \( \Delta_k\setminus E_\n \). By the conditions of the theorem and the definition of \( E_\n \) there exists \( \epsilon>0 \) such that
\[
|x-\xi||A_\n^{(3-k)}(x)|A_\n^{(k)}(x)^2\mu_k^\prime(x)\geq\epsilon
\]
almost everywhere on \( \Delta \). Then, it follows from Lemma~\ref{lem:L}(1,2) that
\[
-\Im\big(D_{\n,\xi}(x+\ic y)L_\n(x+\ic y)\big) = y \int \frac{Q_\n(s)D_{\n,\xi}(s)}{(x-s)^2+y^2} \geq \epsilon\,\int_\Delta\frac{y\,ds}{(x-s)^2+y^2}.
\]
Therefore, for every \( x\in\Delta \) it holds that
\[
\limsup_{y\to0^+}\Im\left(\big(D_{\n,\xi}(x+\ic y)L_\n(x+\ic y)\big)^{-1}\right) \leq 2/(\epsilon\,\pi).
\]
As \( \Delta \) was arbitrary closed subinterval of \( (\Delta_1\cup\Delta_2)\setminus E_\n \) and \( \omega_\n \) has no mass points at the elements of \( E_\n \) by its very definition, \( \omega_\n \) is indeed absolutely continuous. The absolute continuity of \( \rho_O \) can be shown analogously using \eqref{spectral-Markov}, \eqref{rhoO}, and \eqref{mO}. 
\end{proof}

\section{Appendix to Part~\ref{part3}}
\label{s:ap3}

In this appendix we collected some results that were used in the main text.

\subsection{Some properties of \( A_\n^{(k)}(x) \)}
\label{ss:An}

Recall that \( A_{(1,1)}^{(1)}(x) \) and \( A_{(1,1)}^{(2)}(x) \) have degree \( 0 \) and therefore are constants. 
\begin{Lem}
\label{lem:1.4.1}
It holds that
\begin{equation}
\label{378_1}
A_{(1,1)}^{(1)} = - \Xi_{\vec\mu}^{-1}\|\mu_1\|^{-1} \quad \text{and} \quad A_{(1,1)}^{(2)} = \Xi_{\vec\mu}^{-1}\|\mu_2\|^{-1}\,,
\end{equation}
where \( \Xi_{\vec\mu} \) was defined in \eqref{3sd6}. In particular, $A_{(1,1)}^{(1)}<0$ and $A_{(1,1)}^{(2)}>0$\,.
\end{Lem}
\begin{proof}
The claim is a consequence of the fact that \( A_{(1,1)}^{(1)},A_{(1,1)}^{(2)} \) solve the system of equations
\[
\int \bigg(A_{(1,1)}^{(1)}d\mu_1(x)+A_{(1,1)}^{(2)}d\mu_2(x)\bigg) =0 \quad \text{and} \quad \int x\bigg(A_{(1,1)}^{(1)}d\mu_1(x)+A_{(1,1)}^{(1)}d\mu_2(x)\bigg)=1. \qedhere
\]
\end{proof}

Recall that we assumed $\Delta_1<\Delta_2$. Let 
\[
\lambda_{\n,1}\ddd \text{coeff}_{n_1-1}\, A_{\n}^{(1)} \quad \text{and} \quad \lambda_{\n,2}\ddd \text{coeff}_{n_2-1}\, A_{\n}^{(2)}\,.
\]
\begin{Lem} 
\label{lem:1.4.2}
We have that
\[
\sgn \lambda_{\n,1}=(-1)^{n_2} \quad \text{and} \quad \sgn \lambda_{\n,2}=1\,.
\]
\end{Lem}
\begin{proof}
Comparing the leading coefficients in recursion relations \eqref{3sd1} gives \( \lambda_{\n,j}=a_{\n,j}\lambda_{\n+\vec e_j,j} \). By taking into account that  $a_{\n,j}>0$,  we get
\begin{equation}\label{dd1}
\sgn \lambda_{(n_1,n_2),1}=\sgn \lambda_{(1,n_2),1} \quad \text{and} \quad \sgn \lambda_{(n_1,n_2),2}=\sgn \lambda_{(n_1,1),2}\,,\quad \n\in \mathbb{N}^2\,.
\end{equation}
It follows from Lemma~\ref{lem:1.4.1} that $\lambda_{(1,1),1}=A_{(1,1)}^{(1)}<0$ and $\lambda_{(1,1),2}=A_{(1,1)}^{(2)}>0$. Therefore,
\[
\sgn \lambda_{(n_1,1),1}=-1\, \quad \text{and} \quad \sgn \lambda_{(1,n_2),2}=1\,.
\]
It follows from orthogonality conditions \eqref{sad15} for the multi-index $(1,n_2)$ that
\[
\int q(x)\bigg(A_{(1,n_2)}^{(1)}(x)d\mu_1(x)+A_{(1,n_2)}^{(2)}(x)d\mu_2(x)\bigg)=0
\]
for all polynomials $q(x)$ of degree at most $n_2-1$. By taking $q(x)= A_{(1,n_2)}^{(2)}(x)$, we get
\[
-\int \big(A_{(1,n_2)}^{(2)}(x)\big)^2d\mu_2(x) = \int A_{(1,n_2)}^{(1)}(x)A_{(1,n_2)}^{(2)}(x) d\mu_1(x).
\]
Since all the zeroes of $A_{(1,n_2)}^{(2)}(x)$ are on $\Delta_2$ and \( A_{(1,n_2)}^{(1)}=\lambda_{(1,n_2),1}\) is a constant, we get that
\[
-1 = \sgn \lambda_{(1,n_2),1}\cdot  \sgn \lambda_{(1,n_2),2} \cdot (-1)^{n_2-1} = (-1)^{n_2-1}\cdot\sgn \lambda_{(1,n_2),1}
\]
and therefore \( \sgn \lambda_{\n,1}= \sgn \lambda_{(1,n_2),1} = (-1)^{n_2}\) by \eqref{dd1}. That proves the first statement. The second one can be proved similarly.
\end{proof}

Let \( E_{\n,k} \) be the set of zeroes of \( A_\n^{(k)}(x) \), \( k\in\{1,2\} \), and \( E_\n=E_{\n,1}\cup E_{\n,2} \).
 
 \begin{Lem}
 \label{lem:interlacing}
It holds that \( E_{\n,k} \subset \Delta_k \) and \( \#E_{\n,k}=n_k-1 \). That is, all the zeroes of \( A_\n^{(k)}(x) \) are simple and belong to \( \Delta_k \). Write \( E_{\n,k} = \big\{ x_1^{(\n,k)},\ldots,x_{n_k-1}^{(\n,k)} \big\} \), where the zeroes are labeled in the increasing order. The sets \( E_{\n+\vec e_l,k} \) and \( E_{\n,k} \) interlace for any \( k,l\in\{1,2\} \) and
 \begin{equation}
\label{interl1}
x_1^{(\n,2)}< x_1^{(\n+\vec{e}_1,2)}<x_2^{(\n,2)}<\ldots<x_{n_2-1}^{(\n,2)}<x_{n_2-1}^{(\n+\vec{e}_1,2)}
\end{equation}
while
\begin{equation}
\label{interl2}
x_1^{(\n+\vec{e}_2,1)}< x_1^{(\n,1)}<x_2^{(\n+\vec{e}_2,1)}<\ldots<x_{n_1}^{(\n+\vec{e}_2,1)}<x_{n_1}^{(\n,1)}\,
\end{equation}
(in the other two situations the order is uniquely induced by the fact that \( \#E_{\n+\vec e_k,k} = \#E_{\n,k} + 1 \)).
 \end{Lem}
 \begin{proof}
 The statements about location of zeroes and interlacing can be proved in the standard way (see, e.g., \cite[Proposition~2.2 and Theorem~5]{fp} for the proofs). We only need to show \eqref{interl1} and \eqref{interl2}. Let us prove \eqref{interl1}, the argument for \eqref{interl2} is identical. By
 \eqref{3sd1}, we have two identities
 \[
 xA_{\n}^{(2)}(x)=A_{\n-\vec{e}_i}^{(2)}(x) + b_{\n-\vec{e}_i,i}A_{\n}^{(2)}(x) + a_{\n,1}A_{\n+\vec{e}_1}^{(2)}(x) + a_{\n,2}A_{\n+\vec{e}_2}^{(2)}(x)\,,\,\, i\in \{1,2\}\,.
 \] 
 Subtracting one from another, we get
 \[
 A^{(2)}_{\vec n-\vec e_1}(x)-A^{(2)}_{\vec n-\vec e_2}(x)=(b_{\n-\vec{e}_2,2}-b_{\n-\vec{e}_1,1})A_{\n}^{(2)}(x)\,.
 \]
 Taking $x=x^{(\vec n,2)}_{n_2-1}$, the largest zero of $A^{(2)}_{\vec{n}}(x)$, in the previous identity yields
 \begin{equation}\label{3sd2}
 A^{(2)}_{\vec n-\vec e_1}(x^{(\vec n,2)}_{n_2-1})=A^{(2)}_{\vec n-\vec e_2}(x^{(\vec n,2)}_{n_2-1})\,.
 \end{equation}
 The leading coefficients of $\{A^{(2)}_{\vec{m}}(x)\}$ are all positive by Lemma \ref{lem:1.4.2} and  the zeroes of $A^{(2)}_{\vec n-\vec e_2}(x)$ and  $A^{(2)}_{\vec n}(x)$ interlace, so $A^{(2)}_{\vec n-\vec e_2}(x^{(\vec n,2)}_{n_2-1})>0$. Thus, $A^{(2)}_{\vec n-\vec e_1}(x^{(\vec n,2)}_{n_2-1})>0$ by \eqref{3sd2}.  Since the zeroes of $A^{(2)}_{\vec n-\vec e_2}(x)$ and  $A^{(2)}_{\vec n}(x)$ also interlace, we conclude that the zeroes of  $A^{(2)}_{\vec n}(x)$ dominate those of $A^{(2)}_{\vec n-\vec e_1}(x)$.
 \end{proof}
 
Define the polynomials $\{T_{\n,l}(x)\}$ by
\begin{equation}\label{sdpt}
T_{\n,l}(x) \ddd \big(A_{\n+\vec e_l}^{(2)}A_\n^{(1)}-A_{\n+\vec e_l}^{(1)}A_\n^{(2)}\big)(x),\quad l\in \{1,2\}\,.
\end{equation}

\begin{proof}[Proof of Lemma~\ref{lem:Commut}]
It holds by the very definition \eqref{Qn} that
\[
T_{\n,l}(x)d\mu_1(x) = T_{\n,l}(x)d\mu_1(x) \pm A_{\n+\vec e_l}^{(2)}(x)A_\n^{(2)}(x)d\mu_2(x) = A_{\n+\vec e_l}^{(2)}(x)Q_\n(x) -A_\n^{(2)}(x)Q_{\n+\vec e_l}(x).
\]
Since the degree of \( A_{\n+\vec e_l}^{(2)}(x) \) is \( n_2+l-2 \), we get from \eqref{sad15} that
\[
\int x^kT_{\n,l}(x)d\mu_1(x) = 0, \quad k\in\{0,\ldots,n_1-l \}.
\]
Thus, polynomial \( T_{\n,l}(x) \) has at least \( n_1-l+1 \) zeroes on \( \Delta_1 \). Similarly, we can show that $T_{\n,l}(x)$ satisfies  $n_2+l-2$ orthogonality conditions with respect to $\mu_2$ and therefore it has at least \( n_2+l-2 \) zeroes on \( \Delta_2 \). Because its degree is  \( n_1+n_2-1 \), all its zeroes are accounted for and are simple. We can write this polynomials as a product of its leading coefficient and monic polynomials \( T_{\n,l,1}(x) \) and \( T_{\n,l,2}(x) \) that have their zeroes on \( \Delta_1 \) and \( \Delta_2 \), respectively.

Without loss of generality we assume that \( \vec\mu \) satisfies the conditions of Lemma~\ref{lem:Ln} as the general case can be  obtained via weak$^*$ approximation of measures. First, we undo the transformations leading to the definition of \( S_X(x) \). Let \( \n=\Pi(X_{(p)}) \) and \( l=\iota_X \).  It follows from \eqref{commut} that
\[
S_X(x) = \left(\big(A_{\n+\vec e_l}^{(0)}A_\n^{(k)} - A_{\n+\vec e_l}^{(k)}A_\n^{(0)}\big)(x) + (-1)^k\widehat\mu_{3-k}(x) \big(A_{\n+\vec e_l}^{(2)}A_\n^{(1)}-A_{\n+\vec e_l}^{(1)}A_\n^{(2)}\big)(x) \right)
\]
for \( x\in\Delta_k \).  Taking the formulae \eqref{GYX} and \eqref{GYX+} with $Y=X$, we get
\[
\pi S_X(x)\mu_k^\prime(x) = -\Im\big( L_{\n+\vec e_l+}(x) L_{\n-}(x) \big).
\]
On the other hand, it follows from Plemelj-Sokhotski formulae, see \cite[Section~I.4.2]{Gakhov}, that
\[
\pi S_X(x)\mu_k^\prime(x)  = -\Im\left(\left(\mathrm{p.v.}\int_\R\frac{Q_{\n+\vec e_l}(s)}{x-s} - \pi\ic \,\mu_k^\prime(x) \, A_{\n+\vec e_l}^{(k)}(x) \right)\left(\mathrm{p.v.}\int_\R\frac{Q_\n(s)}{x-s} + \pi\ic \,\mu_k^\prime(x) \, A_\n^{(k)}(x) \right)\right)
\]
for \( x\in \Delta_k \), where ``p.v.'' stands for ``principal value''. Notice that it follows form \eqref{sad15} that
\[
P^{-1}(x)\,\mathrm{p.v.}\int_\R\frac{P(s)Q_{\vec m}(s)}{x-s} = P^{-1}(x)\int \frac{P(s)-P(x)}{x-s}Q_{\vec m}(x) + \mathrm{p.v.}\int_\R\frac{Q_{\vec m}(s)}{x-s} = \mathrm{p.v.}\int_\R\frac{Q_{\vec m}(s)}{x-s}
\]
for any polynomial \( P(x) \) of degree at most \( |\vec m|-1 \). In particular, if \( Y=X \) and we let \( l=\iota_X \), in which case \( \vec m=\n+\vec e_l \), then it holds that
\begin{multline*}
A_{\n+\vec e_l}^{(k)}(x)\,\mathrm{p.v.}\int_\R\frac{Q_\n(s)}{x-s} - A_\n^{(k)}(x)\,\mathrm{p.v.}\int_\R\frac{Q_{\n+\vec e_l}(s)}{x-s} \\ = \frac1{T(x)} \mathrm{p.v.}\int_\R\frac{T(s)\big(A_{\n+\vec e_l}^{(k)}(s)Q_\n(s) - A_\n^{(k)}(s)Q_{\n+\vec e_l}(s)\big)}{x-s}  \\ = \frac{(-1)^k}{T(x)}\int_\R\frac{T(s)T_{\n,l}(s)}{x-s}d\mu_{3-k}(s), \quad x\in \Delta_k,
\end{multline*}
for any polynomial \( T(x) \) with real coefficients and of degree at most \( n_2 + l -2 \) if \( k=1 \) and of degree at most \( n_1 - l+1 \) when \( k=2 \). Hence,  taking \( T(x)=T_{\n,l,3-k}(x) \),  we have shown that
\begin{equation}
\label{SX}
S_{\n,l,k}(x) \ddd S_X(x) = \frac{(-1)^k}{T_{\n,l,3-k}(x)}\int_\R\frac{T_{\n,l,3-k}(s)T_{\n,l}(s)}{x-s}d\mu_{3-k}(s), \quad x\in \Delta_k,
\end{equation}
which is clearly a non-vanishing function. To prove positivity, take \( k=1 \). Polynomial \( T_{\n,l,2}(x) \) is monic and has all of its \( n_2+l-2 \) zeroes on \( \Delta_2 \). Thus, its sign on \( \Delta_1 \) is equal to \( (-1)^{n_2+l} \). Polynomial \( T_{\n,l,2}(x)T_{\n,l}(x) \) has double zeroes on \( \Delta_2 \) and the same leading coefficient as \( (-1)^l A_{\n+\vec e_l}^{(l)}(x)A_\n^{(3-l)}(x) \). The latter has the same sign as \( (-1)^{n_2+l} \) by Lemma~\ref{lem:1.4.2}, and therefore,
\[
S_X(x) = \frac{-1}{|T_{\n,l,2}(x)|}\,\int_\R\frac{|T_{\n,l,2}(s)T_{\n,l}(s)|}{x-s}d\mu_2(s) > 0, \quad x<\beta_2,  
\]
as claimed. The case of \( k=2 \) can be considered similarly.
\end{proof}

\subsection{Properties of $L_{\n}(z)$.} 
\label{ss:pLn}

Recall the definitions of \( D_{\n,\xi}(z) \) in \eqref{ur3}, the measure \( \nu_{\n,E} \) in \eqref{nunE}, the polynomials \( T_{\n,k}(x) \) in \eqref{sdpt}, and the functions \( S_{\n,l,k}(x) \) in \eqref{SX}. The set $E_n$ is the set of zeroes of the polynomial $A_{\vec n}^{(1)}(z)A_{\vec n}^{(2)}(z)$.

\begin{Lem}
\label{lem:L}
It holds that
\begin{itemize}
\item[(1)] If \( D(x) \) is a polynomial of degree at most $|\n|-1$, then
 \begin{equation}\label{3sd10}
 L_{\n}(z)=D^{-1}(z)\int_\R \frac{\big(Q_{\n}D\big)(x)}{z-x}.
 \end{equation}
 \item[(2)] The measure \( D_{\n,\xi}(x) Q_\n (x) \) is non-negative on \( \Delta_1\cup \Delta_2 \) for every \( \xi\in(\beta_1,\alpha_2) \).  In particular, $\nu_{\vec n,E}$ is a positive measure. \medskip
 \item[(3)] The function $L_{\n}(z)$ has no zeroes outside $\Delta_1\cup\Delta_2$ and its restriction to $\R\backslash (\Delta_1\cup\Delta_2)$ has well-defined nonzero limits at the endpoints of $\Delta_1$ and $\Delta_2$. \medskip
 \item[(4)] If $E\in E_\n$, then  \( -D_{\n,\xi}^\prime(E)\lim_{\epsilon\to 0^+}L_{\n}(E+\ic\epsilon) = \|\nu_{\n,E}\| - \nu_{\n,E}(\{E\}) >0 \). \medskip
\item[(5)] If \( E\in E_\n\cap\Delta_k \), then \( \|\nu_{\n,E}\| = -D_{\n,\xi}^\prime(E)S_{\n,l,k}(E)/A_{\n+\vec e_l}^{(k)}(E) \) for either \( l\in\{1,2\} \).
\end{itemize}
\end{Lem}
\begin{proof}
(1) The claim follows form orthogonality condition \eqref{sad15},  \eqref{Qn}, and \eqref{Ln} since
\[
0=\int_\R \frac{Q_{\n}(x)(D(x)-D(z))}{x-z}=\int_\R \frac{\big(Q_{\n}D\big)(x)}{x-z}+\big(DL_{\n}\big)(z).
\]
(2) Since $A_{\n}^{(k)}(x)$ has all its zeroes localized to $\Delta_k$, it follows from Lemma~\ref{lem:1.4.2} that
\[
(-1)^{n_2}(x-\xi)A_{\n}^{(2)}(x)>0, ~~ x\in\Delta_1, \quad \text{and} \quad (-1)^{n_2}(x-\xi)A_{\n}^{(1)}(x)>0, ~~x\in\Delta_2,
\]
which yields positivity of \( D_{\n,\xi}(x) Q_\n (x) \).

\smallskip 

\noindent
(3) It follows from claims (2) and (1), applied with \( D(x) = D_{\n,\xi}(x) \),  that \( (D_{\n,\xi} L_\n)(y)<0 \) for \( y\in(-\infty,\alpha_1] \) and \( (D_{\n,\xi} L_\n)(y)>0 \) for \( y\in[\beta_2,\infty) \) (the limits at \( \alpha_1 \) and \( \beta_1 \) might be infinite, but they always exist since Markov functions are decreasing on the real line away from the support of the defining measure). Hence, \( L_\n(x) \) is non-vanishing there. To show that $L_{\n}(x)$ has no zeroes in the lacuna $[\beta_1,\alpha_2]$, take \( D(x)=D_{\n,\eta}(x) \) with \( \eta<\alpha_1 \) and $D_{\n,\eta}$ is defined by $\eqref{ur3}$. Observe that in this case \( (Q_{\n}D_{\n,\eta})(x) \) is non-positive on \( \Delta_1 \) and still non-negative on \( \Delta_2 \). Hence, \( (D_{\n,\eta} L_\n)(y)<\zeta<0 \) for all \( y\in(\beta_1,\alpha_2) \), where \( \zeta = \int_{\Delta_1}(\alpha_2-x)^{-1}(Q_\n D_{\n,\eta})(x) \), which finishes the proof of the desires statement.

\smallskip 

\noindent
(4) Notice that
\[
(x-E)^2d\nu_{\n,E}(x) = (x-E)^2d\widetilde\nu_{\n,E}(x), \quad \widetilde \nu_{\n,E} \ddd \nu_{\n,E} - \nu_{\n,E}(\{E\})\delta_E.
\]
Then, it follow from the dominated convergence theorem (the integrands below are bounded by \( 1 \) in absolute value) that
\begin{equation}
\label{ur2}
\lim_{\epsilon\to 0^+} \int_\R \frac{(x-E)d\nu_{\n,E}(x)}{x-(E+\ic\epsilon)}= \lim_{\epsilon\to 0^+}\int_\R\frac{(x-E)^2d\widetilde \nu_{\n,E}(x)}{(x-E)^2+\epsilon^2} + \ic\lim_{\epsilon\to0^+} \int_\R \frac{\epsilon(x-E)d\widetilde \nu_{\n,E}(x)}{(x-E)^2+\epsilon^2} = \|\widetilde \nu_{\n,E}\| > 0,
\end{equation}
where the last conclusion holds since the measures \( \mu_1,\mu_2 \) have supports of infinite cardinality. Thus, claim (4) follows from claim (1) applied with \( D(x) = D_{\n,\xi}(x)/(x-E) \).

\smallskip

\noindent
(5) For a polynomial \( P(x) \) vanishing at \( E \), let us set \( P(E;x) \ddd P(x)/(x-E) \). Clearly, \( P(E;E) = P^\prime(E) \). Recall that  \( \deg(T_{\n,l,1}) = n_1 - l+1  \) and \( \deg(T_{\n,l,2}) = n_2 + l -2 \). It holds that
\[
\|\nu_{n,E}\| = \int_\R\frac{D_{\n,\xi}(E;x)Q_\n(x)}{x-E} = D_{\n,\xi}^\prime(E)\int_\R\frac{Q_\n(x)}{x-E} = \frac{D_{\n,\xi}^\prime(E)}{T_{\n,l,3-k}(E)A_{\n+\vec e_l}^{(k)}(E)}\int_\R \frac{T_{\n,l,3-k}(x)A_{\n+\vec e_l}^{(k)}(x)Q_\n(x)}{x-E},
\]
where we used the fact that \( Q_\n(x) \) is divisible by \( (x-E) \), orthogonality relations \eqref{sad15}  twice, and Lemma~\ref{lem:interlacing} to observe that \( A_{\n+\vec e_l}^{(k)}(E)\neq0 \). Assume that \( k \in\{1,2\} \) is such that \( E\in\Delta_k \), that is, it is a zero of \( A_\n^{(k)}(x) \). Then
\[
\int_\R \frac{T_{\n,l,3-k}(x)A_\n^{(k)}(x)Q_{\n+\vec e_l}(x)}{x-E} = \int_\R T_{\n,l,3-k}(x)A_\n^{(k)}(E;x)Q_{\n+\vec e_l}(x) = 0,
\]
again, due to orthogonality relations \eqref{sad15}. Therefore, it holds by \eqref{SX} that
\begin{eqnarray*}
\|\nu_{n,E}\| &=& \frac{D_{\n,\xi}^\prime(E)}{T_{\n,l,3-k}(E)A_{\n+\vec e_l}^{(k)}(E)} \int_\R \frac{T_{\n,l,3-k}(s)\big(A_{\n+\vec e_l}^{(k)}(s)Q_\n(s)-A_\n^{(k)}(s)Q_{\n+\vec e_l}(s)\big)}{s-E} \\
&=& -\frac{(-1)^kD_{\n,\xi}^\prime(E)}{T_{\n,l,3-k}(E)A_{\n+\vec e_l}^{(k)}(E)} \int_\R\frac{T_{\n,l,3-k}(s)T_{\n,l}(s)}{E-s}d\mu_{3-k}(s) = -\frac{D_{\n,\xi}^\prime(E)S_{\n,l,k}(E)}{A_{\n+\vec e_l}^{(k)}(E)}
\end{eqnarray*}
as claimed.
\end{proof}

\begin{Lem}
\label{lem:boot-strap}
Assume that \( \vec \mu \) satisfies the conditions of Lemma~\ref{lem:Ln} and that \( (\mu_k^\prime(x))^{-1}\in L^p(\Delta_k)\)  for some \( p>1 \) and each \( k\in\{1,2\} \). Suppose further that there exists \( \gamma\in\{\alpha_1,\beta_1,\alpha_2,\beta_2 \} \) such that \( |L_{\vec\varkappa}(\gamma)|=0 \). Then, \( |L_{\vec\varkappa}(x)|^{-2}\mu_k^\prime(x) \in L^p(\Delta_k) \)  for each \( k\in\{1,2\} \) and
\[
\lim_{y\to 0^+} \ic y \frac{L_{\vec 1}(\gamma +\ic y)}{L_{\vec\varkappa}(\gamma+\ic y)} = 0.
\]
\end{Lem}
\begin{proof}
Clearly, the first claim is obvious unless \( |L_{\vec\varkappa}(x)| \) vanishes at the endpoint of \( \Delta_k \). In the latter situation it follows from Proposition~\ref{prop:Markov}(2-4) that
\[
|L_{\vec\varkappa}(x)|^2 \geq \Im (L_{\vec\varkappa+}(x))^2 = \varkappa_k^2 \Im (\sigma_{k+}(x))^2 = \pi (\varkappa_k /\|\mu_k\|)^2 (\mu_k^\prime(x))^2,
\]
where we used the notation  \( \sigma_k = \|\mu_k\|^{-1}\mu_k \). That yields the desired claim \( |L_{\vec\varkappa}(x)|^{-2}\mu_k^\prime(x) \in L^p(\Delta_k) \).

To prove the limit, assume for definiteness that \( \gamma\in\{\alpha_1,\beta_1\} \). Then, we get that
\[
\lim_{y\to 0^+} \ic y \frac{L_{\vec 1}(\gamma +\ic y)}{L_{\vec\varkappa}(\gamma+\ic y)}= L_{\vec 1}(\gamma)\lim_{y\to 0^+} \left( \varkappa_2\widehat\sigma_2^\prime(\gamma) + \varkappa_1\frac{\widehat\sigma_1(\gamma +\ic y) - \widehat\sigma_1(\gamma) }{\ic y}\right)^{-1},
\]
recall that by Lemma~\ref{lem:Lkappa} the value \( \widehat\sigma_1(\gamma) \) is well-defined. The fraction above can be rewritten as
\[
\frac{\widehat\sigma_1(\gamma +\ic y) - \widehat\sigma_1(\gamma) }{\ic y} =-\left( \int_\R\frac{d\sigma_1(x)}{(\gamma-x)^2+y^2} - \ic y \int_\R\frac{d\sigma_1(x)}{(\gamma-x)((\gamma-x)^2+y^2)}\right),
\]
where the first integral is a strictly decreasing function of \( y \in (0,\infty)\). 

Notice that $\varkappa_1\neq 0$ since otherwise $L_{\vec\varkappa}=L_{\vec{e}_1}$ which  has no zeroes on $\R$.
Then, it only remains to show that \( (\gamma-x)^{-2}\mu_1^\prime(x) \) is not \( L^1 \)-integrable on \( \Delta_1 \). Let \( \Delta_\epsilon = [\alpha_1+\epsilon,\beta_1-\epsilon] \) and \( d\nu(x) = \mu_1^{-1}(x)dx \), which is a finite measure on \( \Delta_1 \). Hence, we get from Cauchy-Schwarz inequality that
\[
\left(\int_{\Delta_\epsilon}\frac{dx}{|x-\gamma|}\right)^2 = \left(\int_{\Delta_\epsilon}\frac{\mu_1^\prime(x)d\nu(x)}{|x-\gamma|}\right)^2 \leq \|\nu_{|\Delta_\epsilon}\| \int_{\Delta_\epsilon} \frac{\mu_1^\prime(x)^2d\nu(x)}{(x-\gamma)^2} = \|\nu_{|\Delta_\epsilon}\|\int_{\Delta_\epsilon} \frac{d\mu_1(x)}{(x-\gamma)^2}
\]
and the desired claim follows by letting \( \epsilon\to 0 \).
\end{proof}

\part{Periodic Jacobi operators on rooted trees and Angelesco systems}
 \label{part4}

In Part~\ref{part1}, we introduced operators $\jackn$, see \eqref{jackn}, defined on finite trees $\mathcal T_{\vec N}$, \( \vec N\in \mathbb N^2 \), see Section~\ref{ss:TN}, and studied their spectra and spectral decompositions. In this part of the paper, we  consider Angelesco system, as in Part~\ref{part3}, see~\eqref{1.10}, in the case when \( \supp\,\mu_i=\Delta_i \), \( d\mu_i(x)=\mu_i^\prime(x)dx \), \( \mu_i(x)>0 \), \( x\in\Delta _i \), and \( \mu'_i(x) \) is a restriction of an analytic function defined around $\Delta_i$.   This situation was studied in great detail in  \cite{ady1} and \cite{ady2}, see also \cite{ya1}.   In particular, it was proved that $\jackn$ converges to a limiting operator $\cal{L}_c^{(i)}$ when $\vec{N}$ goes to infinity along the ray 
\begin{equation}
\label{multi-indices}
\mathcal N_c=\big\{\n:~n_i = c_i|\n| + o(|\n|),~~i\in\{1,2\}\big\}, \quad  (c_1,c_2)=(c,1-c),~~c\in[0,1].
\end{equation}
Hereafter, \(\lim_{\mathcal N_{c}} \) stands for the limit as \( |\n|\to\infty \) and  \( \n\in\mathcal N_{c} \). \smallskip

\section{Definitions}

It was shown in the work of Gonchar and Rakhmanov \cite{gr} that for Angelesco systems with two measures there exists a family of vector equilibrium problems, depending on a parameter \( c\in[0,1] \), whose solutions describe the limiting asymptotics of zeroes of the polynomials  $P_{\n} (z)$, see \eqref{1.2}, along all ray sequences \( \mathcal N_c \). In particular, if an Angelesco system \( \vec\mu \) is as described before \eqref{multi-indices}, then the support of the vector equilibrium measure corresponding to \( c \) is a union of two intervals $\Delta_{c,1}\cup\Delta_{c,2}$ where $\Delta_{c,i}\subseteq \Delta_i$, see, e.g., \cite{gr,ady1} for detail. 

\subsection{Riemann surface}

To define operators $\cal{L}_c^{(i)}$ rigorously, we need  the following Riemann surfaces. Let $\RS_c$ be a 3-sheeted Riemann surface
realized as follows: cut a copy of \( \overline \C \) along  \(
\Delta_{c,1}\cup \Delta_{c,2} \), which henceforth is denoted by
$\RS_c^{(0)}$, the second copy of \( \overline\C \) is cut along \(
\Delta_{c,1} \) and is denoted by $\RS_c^{(1)}$, while the third
copy is cut along \( \Delta_{c,2} \) and is denoted by \(
\RS_c^{(2)} \). These copies are then glued to each other crosswise
along the corresponding cuts. It can
be easily verified that thus constructed Riemann surface has genus
0. We denote by $\pi$ the natural projection from $\RS_c$ to
$\overline\C$ and employ the notation \( \z \) for a generic point
on \( \RS_c \) with \( \pi(\z)=z \) as well as $z^{(i)}$ for a
point on $\RS_c^{(i)}$ with $\pi(z^{(i)})=z$.

Since $\RS_c$ has genus zero, one can arbitrarily prescribe zero/pole divisors of rational functions on $\RS_c$ as long as the degree of the divisor is zero. Clearly, a rational function with a given divisor is unique up to multiplication by a constant. Let  \( \chi_c(\z) \) be the conformal map of \( \RS_c \) onto \( \overline\C \) defined uniquely by the condition
\begin{equation}
\label{chi}
\chi_c\big( z^{(0)} \big) = z + \mathcal O\big(z^{-1} \big),\quad z\to\infty.
\end{equation}
The following constants are going to be central to our investigations in this part of the paper. Let \( A_{c,1},A_{c,2},B_{c,1},B_{c,2} \) be determined by
\begin{equation}
\label{AngPar1}
\chi_c\big( z^{(i)} \big) = B_{c,i} + A_{c,i}z^{-1} + \mathcal O\big(z^{-2} \big), \,\, z\to\infty, \,\, i\in\{1,2\}.
\end{equation}
It was shown in \cite[Proposition~2.1]{ady2} that these constants continuously depend on the parameter \( c \) and have well-defined limits as \( c\to 0^+ \) and \( c\to1^- \), which we denote by \( A_{0,i},B_{0,i} \) and \( A_{1,i},B_{1,i} \), respectively. Moreover, constants \( A_{c,1}>0 \) for all \( c\in[0,1) \) while \( A_{1,1}=0 \) and \( A_{c,2}>0 \) for all \( c\in(0,1] \) while \( A_{0,2}=0 \). 

\subsection{Periodic Jacobi operators on rooted trees}

Let $\mathcal{T},\mathcal V$, and \( O \) be as in Section~\ref{ss:rct}. There are two edges meeting at the root $O$. We label one of them type 1 and the other one -- type 2. Next, consider the children of  $O$. Each of them is coincident with exactly three edges, one of which has already been labeled. We label the remaining two as an edge of type 1 and an edge of type 2. We continue in a similar fashion going down the tree generation by generation  and calling one of the unlabelled edges type 1 and the other one type 2. After assigning types to all the edges, we continue by labeling the vertices. If a vertex $Y$ meets two edges of type 1 and one edge of type 2, we call it a vertex of type 1; otherwise, if it is incident with two edges of type 2 and one edge of type 1, we call it type 2. We do not need to assign any type to the root $O$. Given a vertex $Y\neq O$, we denote its type by $\ell_Y$ (this is similar to the index function introduced in \eqref{imath}). 

Both operators  $\mathcal{L}_c^{(1)}$ and $\mathcal{L}_c^{(2)}$ are Jacobi matrices defined on $\cal{T}$. At a vertex \( Y\neq O \) of type \( \ell_Y \), we  define them by the same formula:
\begin{equation}
\label{Ll1}
(\mathcal{L}_c^{(l)}\psi)_{Y}=\sum_{j\in \{1,2\},{Y'}\sim {Y},\,{\rm type\, of\, edge}\, ({Y},{Y'})=j}\sqrt{A_{c,j}}\psi_{{Y}'}+{B_{c,\ell_Y}}\psi_{{Y}}, \quad l\in \{1,2\};
\end{equation}
and at the root $O$ we define the operators  $\mathcal{L}_c^{(1)}$ and $\mathcal{L}_c^{(2)}$ differently  by writing
\begin{equation}
\label{Ll2}
(\mathcal{L}_c^{(l)}\psi)_{O}=\sum_{j\in \{1,2\},{Y'}\sim {O},\,{\rm type\, of\, edge}\, ({O},{Y'})=j}\sqrt{A_{c,j}}\psi_{{Y}'}+{B_{c,l}}\psi_{{O}}\,, \quad l\in\{1,2\}.
\end{equation}
Recall that $A_{c,j}>0$ when $c\in (0,1)$, but either $A_{c,1}$ or $A_{c,2}$ becomes zero when $c\in \{0,1\}$. The latter cases are trivial and we do not study them, see \cite[Appendix~A]{ady2}.

Our operators $\cal{L}_c^{(l)}$  have  ``periodic coefficients'' and ``self-similar structure''. They are defined on the binary tree and should not be confused with a similar class of Jacobi matrices defined on trees associated with the universal cover of finite connected graphs.  The latter class was studied in several papers, see, e.g.,  \cite{ao1,ao2,sb1}. In the rest of this part, we will apply the arguments from Section~\ref{s:cs} to obtain the spectral decomposition of $\cal{L}_c^{(l)}$ using their generalized eigenfunctions.
 
 The following theorem provides the connection between operators $\cal{L}_c^{(l)}$ and \( \jackn \). It is stated in \cite{ady1} for \( c\in(0,1) \) and is a simple consequence of the results of \cite{ya1}. Its extension to \( c\in\{0,1\} \) was obtained in \cite{ady2}.
 
\begin{Thm}
\label{thm:recurrenceOld}
Let $\vec\mu $ be an Angelesco system \eqref{1.10} such that \( \supp\,\mu_i=\Delta_i \), \( d\mu_i(x)=\mu^\prime(x)dx \), \( \mu_i(x)>0 \), \( x\in\Delta _i \), and \( \mu'_i(x) \) is a restriction of a function analytic around $\Delta_i$ for each \( i\in\{1,2\} \). Further, let the constants \( A_{c,i},B_{c,i} \), \( c\in[0,1] \) and \( i\in\{1,2\} \), be given by \eqref{AngPar1}. Then, the ray limits \eqref{multi-indices} of  coefficients $\big\{a_{\n,i},b_{\n,i}\big\}$ from \eqref{1.5}--\eqref{1.7} exist for any  $c\in (0,1) $ and
\begin{equation}
\label{limit}
\lim_{\mathcal N_c}a_{\n,i} =A_{c,i} \quad \text{and} \quad \lim_{\mathcal N_c}b_{\n,i} =B_{c,i}, \quad i\in\{1,2\}.
\end{equation}
\end{Thm}

In \cite[Section~4.5]{ady1}, this theorem was used to prove that \(  \jo_{\vec e_l,\vec{N}}\to \mathcal{L}_c^{(l)} \), \( l\in \{1,2\} \), when $\vec{N}\in \mathcal N_c$ converges to infinity. This convergence can be understood as the strong operator convergence on the same Hilbert space $\ell^2(\cal{T})$ when $\jo_{\vec e_l,\vec{N}}$ is properly extended to this space. 

\subsection{Green's functions}

 In \cite[Appendix A]{ady2}, it was proved that $\sigma(\mathcal{L}_c^{(j)})=\Delta_{c,1}\cup\Delta_{c,2}$ and the spectrum  is purely absolutely continuous. Moreover, if we denote Green's functions of \( \mathcal{L}_c^{(l)} \) corresponding to the root \( O \) by
\begin{equation}
\label{g-lo}
G_c^{(l)}(Y,O;z) \ddd \big\langle(\mathcal{L}_c^{(l)}-z)^{-1}\delta^{(O)},\delta^{(Y)}\big\rangle,
\end{equation}
then it was shown in \cite[Section 4.5]{ady1} that
\begin{equation}
\label{spm}
G_c^{(l)}(O,O;z) = M_{c}^{(l)}(z^{(0)}), \quad z\notin \Delta_{c,1}\cup\Delta_{c,2}\,,
\end{equation}
where $M_{c}^{(l)}(\z)$ is a function on \( \RS_c \) given by
\begin{equation}
\label{mki}
M_{c}^{(l)}(\z) \ddd \frac{1}{B_{c,l}-\chi_{c}(\z)},\quad  l\in \{1,2\}\,.
\end{equation}
Clearly, \( M_{c}^{(l)}(\z) \) is an analytic function on \( \RS_c \) apart from a single pole at \( \infty^{(l)} \), which is simple. Therefore, the traces \( G_c^{(l)}(O,O;x)_\pm \) exists and are continuous on $\Delta_{c,1}\cup\Delta_{c,2}$. Moreover, they are complex conjugates of each other. In particular, \( |G_c^{(l)}(O,O;x)|\) is well-defined for all \( x\in\Delta_{c,1}\cup\Delta_{c,2} \). 

\begin{Lem}
The identity
\begin{equation}
\label{ap1}
A_{c,1}|G_c^{(1)}(O,O;x)|^2+A_{c,2}|G_c^{(2)}(O,O;x)|^2 = 1
\end{equation}
holds for each $x\in  \Delta_{c,1}\cup\Delta_{c,2} $. Moreover, 
\begin{equation}
\label{ap11}
A_{c,1}|G_c^{(1)}(O,O;z)|^2 + A_{c,2}|G_c^{(2)}(O,O;z)|^2 < 1
\end{equation}
for  $z\notin \Delta_{c,1}\cup\Delta_{c,2}$.
\end{Lem}
\begin{proof} From \cite[formula (4.27)]{ady1}, we get that
\begin{equation}
\label{fo2}
z=-1/M_c^{(l)}(\z) + B_{c,j} - A_{c,1}M_c^{(1)}(\z) -  A_{c,2}M_c^{(2)}(\z)
\end{equation}
for each \( l\in\{1,2\} \) and \( \z\in\RS_c \). Formula \eqref{fo2}, in particular, implies that
\[
B_{c,1} - 1/M_c^{(1)}(\z) = B_{c,2} - 1/M_c^{(2)}(\z)
\]
for all $\z\in \RS_c$. Fix \( i\in\{1,2\} \). Using the above relation with \( \z = z^{(3-i)} \) gives us
\begin{equation}
\label{fo1}
\frac{1}{M_c^{(1)}(z^{(3-i)})}-\frac{1}{M_c^{(2)}(z^{(3-i)})}=B_{c,1}-B_{c,2}.
\end{equation}
Since the product of all the branches of an algebraic function is a polynomial, behavior at infinity yields that
\[
M_c^{(l)}(z^{(0)})M_c^{(l)}(z^{(1)})M_c^{(l)}(z^{(2)})=(-1)^l(A_{c,l}(B_{c,2}-B_{c,1}))^{-1}.
\]
By plugging the above relations into \eqref{fo1} we get
\[
A_{c,1}M_{c}^{(1)}{(z^{(0)})}M_{c}^{(1)}{(z^{(i)})} + A_{c,2} M_{c}^{(2)}{(z^{(0)})} M_{c}^{(2)}{(z^{(i)})}=1
\]
for all $z \in \C\setminus(\Delta_{c,1}\cup \Delta_{c,2})$. Taking the boundary values on \( \Delta_i \) from the upper half-plane, we obtain
\[
A_{c,1}M_{c+}^{(1)}{(x^{(0)})}M_{c+}^{(1)}{(x^{(i)})} + A_{c,2} M_{c+}^{(2)}{(x^{(0)})} M_{c+}^{(2)}{(x^{(i)})}=1,
\]
for \( x\in\Delta_i \). To prove \eqref{ap1}, it only remains to observe that
\[
G_c^{(l)}(O,O;x)_\pm = M_{c\pm}^{(l)}\big(x^{(0)}\big) = M_{c\mp}^{(l)}\big(x^{(i)}\big)
\]
for \( x\in\Delta_i \) in view of \eqref{spm}. To show \eqref{ap11} observe that its right-hand side is subharmonic, decays at infinity, and equals $1$ on the cuts. Thus, the maximum principle gives the claimed bound.
\end{proof}

\begin{Rem}
Identity \eqref{ap1} gives a simple description of the image of the cuts $\Delta_{c,1}$ and $\Delta_{c,2}$ under the conformal map $\chi_{c}(\z)$. Namely, this image is a contour in the plane described by the equation
\begin{equation}
\frac{A_{c,1}}{|\chi-B_{c,1}|^2}+\frac{A_{c,2}}{|\chi-B_{c,2}|^2}=1\,, \quad \chi\in \C\,.
\end{equation}
\end{Rem}
The self-similar nature of the operators \( \mathcal L^{(l)}_c \) and \eqref{spm} make it possible to compute their Green's functions.

\begin{Prop}
\label{prop:GlXO}
For \( z\notin \Delta_{c,1}\cup\Delta_{c,2} \) and \( X\neq O \), it holds that
\begin{equation}
\label{GlXO}
G_c^{(l)}(X,O;z) = M_c^{(l)}(z^{(0)}) \cdot \prod_{Y\in \mathrm{path}^*(X,O)}\Bigl(-A_{c,\ell_Y}^{1/2}\Bigr)M_{c}^{(\ell_Y)}(z^{(0)})\,,
\end{equation}
where $\mathrm{path}^*(X,O)$ is the path that connects $O$ to $X$, it includes $X$, but excludes $O$. Moreover, 
\begin{equation}
\label{GlOn}
\Bigl\| G_c^{(l)}(\cdot,O;z)  \Bigr\|_{\ell^2(\cal{V})}^2 = \frac{|M_c^{(l)}(z^{(0)})|^2}{1-(A_{c,1}|M_{c}^{(1)}(z^{(0)})|^2+A_{c,2}|M_{c}^{(2)}(z^{(0)})|^2)}
\end{equation}
for all \( z\notin \Delta_{c,1}\cup\Delta_{c,2} \), where we consider \( \{G_c^{(l)}(Y,O;z) \}\) as a function of \( Y \) on \( \mathcal V \).
\end{Prop}
\begin{proof}
Let \( g(z) \) be a function on \( \mathcal V \) given by the right-hand side of \eqref{GlXO} with \( g_O(z) \ddd M_c^{(l)}(z^{(0)}) \). By induction in $n\in \mathbb{N}$, one gets that
\[
\sum_{|Y|=n} |g_Y(z)|^2 = |M_c^{(l)}(z^{(0)})|^2\left(A_{c,1}|M_{c}^{(1)}(z^{(0)})|^2+A_{c,2}|M_{c}^{(2)}(z^{(0)})|^2\right)^n,
\]
where \( |Y| \) stands for the distance from \( Y \) to the root \( O \). Therefore, it follows from \eqref{ap11} that \( \|g(z)\|^2_{\ell^2(\mathcal V)}\) is finite and is equal to the right-hand side of \eqref{GlOn} for all \( z\notin \Delta_{c,1}\cup\Delta_{c,2} \). Thus, to prove the lemma we only need to show that \( (\mathcal L_c^{(l)}-z) g(z) =\delta^{(O)} \). The latter is a straightforward application of \eqref{Ll1} and \eqref{Ll2}. Indeed, let \( Y\neq O \) be of type \( i \) and \( Y_1 \) and \( Y_2 \) be the children of \( Y \) of types \( 1 \) and \( 2 \), respectively. Then
\begin{eqnarray*}
\big((\mathcal L_c^{(l)}-z) g(z)\big)_Y &=& (B_{c,i}-z)g_Y(z) +  \sqrt{A_{c,i}}g_{Y_{(p)}}(z) + \sqrt{A_{c,1}}g_{Y_1}(z) + \sqrt{A_{c,2}}g_{Y_2}(z) \\
&=& g_Y(z) \left(B_{c,i}-z - M_c^{(i)}\big(z^{(0)}\big)^{-1} - A_{c,1}M_c^{(1)}\big(z^{(0)}\big) - A_{c,2}M_c^{(2)}\big(z^{(0)}\big) \right) = 0,
\end{eqnarray*}
where the last equality follows from \eqref{fo2}. Similarly, it holds that
\begin{eqnarray*}
\big((\mathcal L_c^{(l)}-z) g(z)\big)_O &=& (B_{c,l}-z)g_O(z) + \sqrt{A_{c,1}}g_{O_1}(z) + \sqrt{A_{c,2}}g_{O_2}(z) \\
&=& M_c^{(l)}(z^{(0)}) \left(B_{c,i}- z - A_{c,1}M_c^{(1)}\big(z^{(0)}\big) - A_{c,2}M_c^{(2)}\big(z^{(0)}\big) \right) = 1,
\end{eqnarray*}
where \( O_1 \) and \( O_2 \) be the children of \( O \) of types \( 1 \) and \( 2 \), respectively.
\end{proof}

\begin{Rem}
Direct algebraic proof  of \eqref{GlXO}, rather than a posteriori computation given above, can be found in \cite[Remark~4.15]{ady1}.
\end{Rem}

\section{Spectral analysis}

To carry our spectral analysis of the operators $\mathcal{L}_c^{(l)}$ we follow the blueprint of Sections~\ref{s:gf}--\ref{s:sa}.

\subsection{Trivial cyclic subspaces of $\mathcal{L}_c^{(l)}$ generated by $\delta^{(O)}$}

From \eqref{chi} and the symmetries of the surface \( \RS_c \), one can deduce that \( \chi_c(z^{(0)}) \)  has positive imaginary part when \( z\in \C_+ \), i.e.,  that \( \chi_c(z^{(0)})\in \HN \). That is consistent with \( G_c(O,O;\cdot)\in\HN \) due to  \eqref{spm} and \eqref{mki}. It is indeed a negative of a Markov function of the spectral measure of $\mathcal{L}_c^{(l)}$ with respect to \( \delta^{(O)} \). Let us denote this spectral measure by \( \rho_O^{(c,l)} \). Since functions \( M_c^{(l)}(\z) \) map the surface \( \RS_c \) conformally onto the Riemann sphere, it follows from Proposition~\ref{prop:Markov}(1-3) and \eqref{spm} that
\[
d\rho_O^{(c,l)}(x) = \Im\big(M_c^{(l)}\big(x^{(0)}_+\big)\big)dx, \quad x\in\Delta_{c,1}\cup\Delta_{c,2},
\]
where \( x^{(0)}_+ \ddd \lim_{y\to0^+}z^{(0)} \), \( z=x+\ic y \). 
Define the reference measure $\omega^{(c)}$ as the sum of two
\[
d\omega^{(c)}(x) \ddd \sqrt{|(x-\alpha_{c,1})(x-\beta_{c,1})(x-\alpha_{c,2})(x-\beta_{c,2})|}dx,
\]
where we write \( \Delta_{c,i}=[\alpha_{c,i},\beta_{c,i}] \) (in fact, it always holds that \( \alpha_{c,1}=\alpha_1 \) and \( \beta_{c,2}=\beta_2\)).
The analysis of the conformal map $\chi(\z)$ at the endpoints of $\Delta_{c,i}$ reveals that the densities of both spectral measures $\rho_O^{(c,l)}$ satisfy
\[
C_1 (\omega^{(c)})'(x)<(\rho_O^{(c,l)})'(x)<C_2 (\omega^{(c)})'(x)
\]
for $x\in \Delta_{c,i}$ and some positive constants \( C_1,C_2 \) that might depend on $c$ but do not depend on $x$.  In particular, if we define $\nu^{(c,l)}(x)\ddd (\rho_O^{(c,l)})'(x)/(\omega^{(c)})'(x)$, then
\begin{equation}\label{emph2}
\nu^{(c,l)}\in L^\infty(\Delta_{c,1}\cup\Delta_{c,2}), \quad
(\nu^{(c,l)})^{-1}\in L^\infty(\Delta_{c,1}\cup\Delta_{c,2})
\end{equation}
 for each $l\in \{1,2\}$. Similarly to \eqref{987-2}, we can then define 
\begin{eqnarray}
\nonumber
\Psi_X^{(c,l)}(x) & \ddd & \frac{d\Im G_c^{(l)}(X,O)^+(x)}{d\rho_O^{(c,l)}(x)} \\
\label{Psicl}
&=& \Im\big(M_c^{(l)}\big(x^{(0)}_+\big)\big)^{-1}\Im\left(M_c^{(l)}\big(x^{(0)}_+\big)\prod_{Y\in {\rm path}^*(X,O)}\Bigl(-A_{c,\ell_Y}^{1/2}\Bigr)M_{c}^{(\ell_Y)}\big(x_+^{(0)}\big)\right)
\end{eqnarray}
for \( X\in\mathcal V \) and \( x\in\Delta_{c,1}\cup\Delta_{c,2} \), where the second equality follows from \eqref{GlXO}. Notice that the same computation as in the second part of the proof of Proposition~\ref{prop:GlXO} shows that $\Psi^{(c,l)}(x)$ is a formal generalized eigenvector for $\mathcal{L}_c^{(l)}$ corresponding to $x\in\Delta_{c,1}\cup\Delta_{c,2}$ that satisfies $\Psi^{(c,l)}_O(x)=1$.

Denote by $\mathfrak{C}_{c,l}^{(O)}$ the cyclic space generated by $\delta^{(O)}$ and $\mathcal{L}_c^{(l)}$. 
Recall that the operator $\alpha(\mathcal{L}_c^{(l)})$ can be defined for every continuous function $\alpha$ using the Spectral Theorem for self-adjoint operators. The proof of the next proposition repeats the proof of Proposition~\ref{prop:triv-c}.

\begin{Prop}
\label{sadok1}
The map
\[
\alpha(x) \mapsto \widehat \alpha^{(c,l)} = \left\{\widehat\alpha_Y^{(c,l)}\right\}_{Y\in\mathcal V}, \quad \widehat \alpha_Y^{(c,l)} \ddd \int \alpha(x)\Psi_Y^{(c,l)}(x)d\rho_O^{(c,l)}(x),
\]
is a unitary map from $L^2(\rho_O^{(c,l)})$ onto $\mathfrak C_{(c,l)}^{(O)}$. In particular,  it holds that
\[
\|\alpha\|_{L^2(\rho_O^{(c,l)})}^2 = \big\|\widehat\alpha^{(c,l)}\big\|_{\ell^2(\mathcal V)}^2 \quad \text{and} \quad \mathfrak C_{(c,l)}^{(O)} = \left\{\widehat\alpha^{(c,l)}:~\alpha\in L^2(\rho_O^{(c,l)})\right\}.
\]
Thus, the formula
\[
\alpha(\mathcal L_c^{(l)})\delta^{(O)} \ddd \widehat\alpha^{(c,l)} = \int \alpha(x)\Psi^{(c,l)}(x)d\rho_O^{(c,l)}(x) 
\]
extends the definition of \(  \alpha(\mathcal L_c^{(l)})\delta^{(O)}  \) to all  \( \alpha\in L^2(\rho_O^{(c,l)}) \). We also have that
\[
x\alpha(x) \mapsto \mathcal L_c^{(l)}\widehat\alpha, \quad \alpha\in L^2(\rho_O^{(c,l)}).
\]
\end{Prop}

\subsection{Nontrivial cyclic subspaces of $\mathcal{L}_c^{(l)}$}

Let \( X\in\mathcal V \) and \( X_1, X_2 \) be children of \( X \) of types \( 1 \) and \( 2 \), respectively. Observe that the restriction of \( \mathcal{L}_c^{(l)} \) to \( \mathcal T_{[X_i]} \) is equal to \( \mathcal{L}_c^{(i)} \), where, as before, \( \mathcal T_{[X_i]} \) is the subtree of \( \mathcal T \) with root at \( X_i \). Here, we can use the self-similar structure to naturally identify \( \mathcal T_{[X_i]} \) with \( \mathcal T \) when talking about the operator \( \mathcal{L}_c^{(i)} \) on \( \mathcal T_{[X_i]} \). Let us further denote by \( \Psi^{(c)}(X_i;x) \) the function \( \Psi^{(c,i)}(x) \), defined in \eqref{Psicl}, carried to \( \mathcal V_{[X_i]} \) from \( \mathcal V \) by using this natural identification. Similarly to \eqref{hatPsi} define
\[
\widehat \Psi_Y^{(c)}(X;x) \ddd (-1)^iA_{c,i}^{-1/2}\Psi_Y^{(c)}(X_i;x),~~Y\in\mathcal V_{[X_i]}, \quad \text{and} \quad \widehat \Psi_Y(X;x) \ddd 0, ~~\text{otherwise}.
\]
Observe that \( \widehat \Psi^{(c)}(X;x) \) does not depend on \( l \) and it follows from \eqref{Ll1} and \eqref{Ll2} that
\[
\left(\big(\mathcal L_c^{(l)}-x\big)\widehat \Psi^{(c)}(X;x)\right)_X = A_{c,1}^{1/2}\widehat \Psi_{X_1}^{(c)}(X;x) + A_{c,2}^{1/2}\widehat \Psi_{X_2}^{(c)}(X;x) = 0.
\]
Similarly to \eqref{Chat}, define
\[
\widehat{\mathfrak C}_c^{(X)} \ddd \left\{\int \alpha(x)\widehat\Psi^{(c)}(X;x)d\omega^{(c)}(x):~~\alpha\in L_{\omega^{(c)}}^2(\Delta_{c,1}\cup\Delta_{c,2})\right\}.
\]
The following proposition is analogous to Proposition~\ref{prop:Chat} and can be proven similarly using \eqref{emph2} and Proposition \ref{sadok1}.

\begin{Prop}
Fix \( X\in\mathcal V \) and let \( X_1, X_2 \) be children of \( X \) of types \( 1 \) and \( 2 \), respectively. The function \( \widehat\Psi^{(c)}(X;x) \) is a generalized eigenvector of \( \mathcal L_c^{(l)} \), that is, it holds that
\[
\mathcal L_c^{(l)} \widehat\Psi^{(c)}(X;x) = x\widehat\Psi^{(c)}(X;x).
\]
Moreover, let the function \( g_{c,i}^{(X)}\in\widehat{\mathfrak C}_c^{(X)} \), \( i\in\{1,2\} \), be given by
\[
g_{c,i}^{(X)} \ddd \int\alpha^{(c)}(X_i;x)\widehat\Psi^{(c)}(X;x)d\omega^{(c)}(x), \quad \alpha^{(c)}(X_i;x) \ddd (-1)^i A_{c,i}^{1/2}\nu^{(c,i)}(x).
\]
Then, it holds that \( \chi_ig_{c,i}^{(X)} = \chi_i\delta^{(X_i)} \), where \( \chi_i \) is the restriction operator that sends \( f\in \widehat{\mathfrak C}_c^{(X)}\) to its restriction to \( \mathcal V_{[X_i]} \), and
\[
\widehat{\mathfrak C}_c^{(X)} = \overline{ \mathrm{span}\left\{\big(\mathcal L_c^{(l)}\big)^n g_{c,i}^{(X)}:~n\in\Z_+\right\} }.
\]
That is, each $g_{c,i}^{(X)}$ is a generator of the cyclic subspace $\widehat{\mathfrak C}_c^{(X)}$. In particular, the formula
\[
\alpha\big(\mathcal L_c^{(l)}\big) g_{c,i}^{(X)} \ddd \int \alpha(x)\alpha^{(c)}(X_i;x)\widehat\Psi^{(c)}(X;x)d\omega^{(c)}(x)
\]
extends the definition of \( \alpha\big(\mathcal L_c^{(l)}\big) g_{c,i}^{(X)} \) to all \( \alpha\in L_{\omega^{(c)}}^2(\Delta_{c,1}\cup\Delta_{c,2}) \). 
Furthermore, it holds that
\[
d\rho_{g_{c,i}^{(X)}}(x) = \sum_{k=1}^2 \frac{A_{c,i}}{A_{c,k}}\frac{\nu^{(c,i)}(x)^2}{\nu^{(c,k)}(x)} d\omega^{(c)}(x),
\]
where \( \rho_{g_{c,i}^{(X)}} \) is the spectral measure of \( g_{c,i}^{(X)} \).
\end{Prop}

\subsection{Orthogonal decomposition}

The proof of the following statement repeats that of Theorem~\ref{mt1}.

\begin{Thm} 
\label{mt2}
The Hilbert space \( \ell^2(\mathcal V) \) decomposes into an orthogonal sum of cyclic subspaces of \( \mathcal{L}_c^{(l)} \) as follows:
\begin{equation}
\label{g44}
\ell^2(\cal{V})=\mathfrak{C}_{c,l}^{(O)}\oplus \mathcal{L}, \quad \mathcal{L}= \oplus_{Z\in \cal{V}} \widehat{\mathfrak{C}}_c^{(Z)}, \quad l\in \{1,2\}\,.
\end{equation}
\end{Thm}

\begin{Rem} 
This theorem implies immediately that $\sigma(\mathcal{L}_c^{(l)})=\Delta_{c,1}\cup\Delta_{c,2}$, that the spectrum is purely absolutely continuous, and that it has infinite multiplicity.
\end{Rem}

\end{document}